\documentclass[11pt,leqno]{article}
\normalfont
\renewcommand{\baselinestretch}{1.25}

\usepackage{graphicx}

\newcommand{\IntroSection}{1}
\newcommand{\IntroFig}{Figure 1.1}
\newcommand{\IntroFigTwo}{Figure 1.2}
\newcommand{\OpenQuestionOne}{Open Question 1.1}
\newcommand{\OpenQuestionTwo}{Open Question 1.2}

\newcommand{\SetupSectionOne}{2}
\newcommand{\SpecializationTheorem}{Proposition 2.1}
\newcommand{\WInvariantLemma}{Lemma 2.2}

\newcommand{\SetupSectionTwo}{3}
\newcommand{\SupportingGraphProp}{Proposition 3.1}
\newcommand{\RepDiagramTheorem}{Proposition 3.2}

\newcommand{\CombinatorialFeatures}{4}
\newcommand{\MainCombinatorialTheorem}{Proposition 4.1}
\newcommand{\TypeAFormulas}{Proposition 4.2}
\newcommand{\CombinatorialPropositions}{Propositions 4.1 and 4.2}

\newcommand{\GTLatticesSection}{5}
\newcommand{\GTParallelogramFigure}{Figure 5.1}
\newcommand{\GTTriangleFigure}{Figure 5.2}
\newcommand{\ParallelogramFig}{Figure 5.3}
\newcommand{\FirstGTLatticeResult}{Proposition 5.1}
\newcommand{\LRCorollaryOne}{Theorem 5.2}

\newcommand{\GTStatementsSection}{6}
\newcommand{\ObservationalLemma}{Lemma 6.1}
\newcommand{\ConsistencyLemma}{Proposition 6.2}
\newcommand{\GTBasisTheorem}{Theorem 6.3}
\newcommand{\GTGeneralTheorem}{Theorem 6.4}

\newcommand{\StandardizeProp}{Proposition 6.5}
\newcommand{\GTCombinatorics}{Corollary 6.6}

\newcommand{\OrthogonalSection}{7}
\newcommand{\OddOrthPatternFigure}{Figure 7.1}
\newcommand{\EvenOrthPatternFigure}{Figure 7.2}
\newcommand{\OrthogonalBijectionFig}{Figure 7.3}
\newcommand{\NewOrthogonalConstructions}{Theorem 7.1}
\newcommand{\NewOrthogonalCorollary}{Corollary 7.2}
\newcommand{\NewOrthogonalCombinatorics}{Corollary 7.3}

\newcommand{\SolitarySection}{8}
\newcommand{\StarProposition}{Proposition 8.1}
\newcommand{\StarIsomorphisms}{Proposition 8.2}
\newcommand{\StarProps}{Propositions 8.1/8.2}
\newcommand{\GTSolitaryMinimal}{Theorem 8.3}
\newcommand{\SolitarySigmaFig}{Figure 8.1}
\newcommand{\SolitaryDualFig}{Figure 8.2}
\newcommand{\SolitaryBowtieFig}{Figure 8.3}
\newcommand{\SolitaryFigs}{Figures 8.1-3}
\newcommand{\OrthogonalSolitaryTheorem}{Theorem 8.4}

\usepackage{amsfonts}
\usepackage{amssymb}
\usepackage{mathrsfs}
\usepackage{rotating}
\usepackage{mathdots}

\newfont{\myscbolditalics}{ecoc0500 at 11pt}

\newfont{\mybolditalics}{ecbi0500 at 11pt}

\newcommand{\myqx}{\mbox{\mybolditalics x}}
\newcommand{\myqy}{\mbox{\mybolditalics y}}
\newcommand{\myqX}{\mbox{\mybolditalics X}}
\newcommand{\myqY}{\mbox{\mybolditalics Y}}
\newcommand{\myqh}{\mbox{\mybolditalics h}}
\newcommand{\myqP}{\mbox{\mybolditalics P}}
\newcommand{\myqQ}{\mbox{\mybolditalics Q}}

\newfont{\eulercursive}{eurm10 at 11pt}
\newcommand{\myd}{\mbox{\eulercursive d}}

\newcommand{\mys}{\mbox{\eulercursive s}}
\newcommand{\mym}{\mbox{\eulercursive m}}

\newfont{\smalleulercursive}{eurm10 at 9pt}

\newfont{\smallereulercursive}{eurm10 at 7pt}

\newfont{\myslantcyrillic}{wncyi10 at 11pt}

\newcommand{\QED}{\raisebox{0.5mm}{\fbox{\rule{0mm}{1.5mm}\ }}}

\newcounter{myfn}[page]
\renewcommand{\thefootnote}{\fnsymbol{footnote}}
\newcommand{\myfootnote}[1]{\setcounter{footnote}{\value{myfn}}%
    \footnote{#1}\stepcounter{myfn}}
    
\newcounter{rone}
\newcounter{rtwo}
\newcounter{rthree}
\newcounter{rfour}
\newcounter{rfive}
\newcounter{rsix}
\newcounter{rseven}
\newcommand{\myA}{\mbox{\sffamily A}}

\newcommand{\mytinyA}{\mbox{\tiny \sffamily A}}
\newcommand{\myB}{\mbox{\sffamily B}}

\newcommand{\mytinyB}{\mbox{\tiny \sffamily B}}
\newcommand{\myC}{\mbox{\sffamily C}}

\newcommand{\myD}{\mbox{\sffamily D}}

\newcommand{\mytinyD}{\mbox{\tiny \sffamily D}}
\newcommand{\myE}{\mbox{\sffamily E}}

\newcommand{\myF}{\mbox{\sffamily F}}

\newcommand{\myG}{\mbox{\sffamily G}}

\newcommand{\mysmallM}{\mbox{\footnotesize \sffamily M}}
\newcommand{\mytinyM}{\mbox{\tiny \sffamily M}}

\newcommand{\mysmallN}{\mbox{\footnotesize \sffamily N}}
\newcommand{\mytinyN}{\mbox{\tiny \sffamily N}}

\newcommand{\mysmallP}{\mbox{\footnotesize \sffamily P}}
\newcommand{\myscriptsizeP}{\mbox{\scriptsize \sffamily P}}
\newcommand{\mytinyP}{\mbox{\tiny \sffamily P}}

\newcommand{\mysmallQ}{\mbox{\footnotesize \sffamily Q}}
\newcommand{\myscriptsizeQ}{\mbox{\scriptsize \sffamily Q}}
\newcommand{\mytinyQ}{\mbox{\tiny \sffamily Q}}

\newcommand{\mysmallS}{\mbox{\footnotesize \sffamily S}}

\newcommand{\myX}{\mbox{\sffamily X}}

\newcommand{\mytinyX}{\mbox{\tiny \sffamily X}}

\newcommand{\myvarZ}{\mbox{\scriptsize \sffamily Z}}

\newcommand{\melt}{\mathbf{m}} \newcommand{\nelt}{\mathbf{n}}
 
 \newcommand{\relt}{\mathbf{r}}
\newcommand{\selt}{\mathbf{s}} \newcommand{\telt}{\mathbf{t}}
\newcommand{\uelt}{\mathbf{u}} 
 \newcommand{\xelt}{\mathbf{x}}
\newcommand{\yelt}{\mathbf{y}} \newcommand{\zelt}{\mathbf{z}}

\newcommand{\wt}{\mbox{\sffamily wt}}
\newcommand{\smallwt}{\mbox{\scriptsize \sffamily wt}}
\newcommand{\tinywt}{\mbox{\tiny \sffamily wt}}

\newcommand{\mychar}{\mbox{\sffamily char}}

\newcommand{\WGF}{\mbox{\sffamily WGF}}
\newcommand{\RGF}{\mbox{\sffamily RGF}}
\newcommand{\myCARD}{\mbox{\sffamily CARD}}
\newcommand{\myLENGTH}{\mbox{\sffamily LENGTH}}

\newcommand{\comp}{\mbox{\sffamily comp}}

\newcommand{\myarrow}[1]{\stackrel{#1}{\rightarrow}}

\newcommand{\mylongarrow}[1]{\stackrel{#1}{\longrightarrow}}

\newcommand{\NEEdgeLabelForLatticeI}[1]{
\setlength{\unitlength}{1.5cm}
\begin{picture}(0,0)
\put(-0.25,0){
\begin{picture}(0,0)
\put(0.4,0.4){\footnotesize #1} 
\end{picture}
}
\end{picture}
}

\newcommand{\NWEdgeLabelForLatticeI}[1]{
\setlength{\unitlength}{1.5cm}
\begin{picture}(0,0)
\put(-0.25,0){
\begin{picture}(0,0)
\put(-0.525,0.4){\footnotesize #1} 
\end{picture}
}
\end{picture}
}

\newcommand{\VerticalEdgeLabelForLatticeI}[1]{
\setlength{\unitlength}{1.5cm}
\begin{picture}(0,0)
\put(-0.25,0){
\begin{picture}(0,0)
\put(-0.05,0.4){\footnotesize #1} 
\end{picture}
}
\end{picture}
}

\newcommand{\VertexTableau}[6]{
\setlength{\unitlength}{1.5cm}
\begin{picture}(0,0)
\put(-0.25,0){
\begin{picture}(0,0)
\put(0,0){\circle*{0.1}} 
\put(#5,#6){\setlength{\unitlength}{0.3cm}\begin{picture}(0,0)\put(0,0){\line(0,1){1}} \put(1,0){\line(0,1){1}} \put(2,0){\line(0,1){2}} \put(3,0){\line(0,1){2}} \put(0,0){\line(1,0){3}} \put(0,1){\line(1,0){3}} \put(2,2){\line(1,0){1}} \put(0.25,0.25){\tiny #1} \put(1.25,0.25){\tiny #2} \put(2.25,1.25){\tiny #3} \put(2.25,0.25){\tiny #4} \end{picture}}
\end{picture}
}
\end{picture}
}

\newcommand{\VertexTableauFib}[7]{
\setlength{\unitlength}{1.5cm}
\begin{picture}(0,0)
\put(-0.25,0){
\begin{picture}(0,0)
\put(0,0){\circle*{0.1}} 
\put(#6,#7){\setlength{\unitlength}{0.25cm}\begin{picture}(0,0)\put(0,-1){\line(0,1){2}} \put(1,-1){\line(0,1){2}} \put(2,0){\line(0,1){2}} \put(3,0){\line(0,1){2}} \put(0,-1){\line(1,0){1}} \put(0,0){\line(1,0){3}} \put(0,1){\line(1,0){3}} \put(2,2){\line(1,0){1}} \put(0.25,0.25){\tiny #1} \put(0.25,-0.75){\tiny #2} \put(1.25,0.25){\tiny #3} \put(2.25,1.25){\tiny #4} \put(2.25,0.25){\tiny #5} \end{picture}}
\end{picture}
}
\end{picture}
}

\newcommand{\VertexTableauTwo}[7]{
\setlength{\unitlength}{1.5cm}
\begin{picture}(0,0)
\put(-0.25,0){
\begin{picture}(0,0)
\put(0,0){\circle*{0.1}} 
\put(#6,#7){\setlength{\unitlength}{0.3cm}\begin{picture}(0,0)\put(0,0){\line(0,1){2}} \put(1,0){\line(0,1){2}} \put(2,0){\line(0,1){2}} \put(3,1){\line(0,1){1}} \put(0,0){\line(1,0){2}} \put(0,1){\line(1,0){3}} \put(0,2){\line(1,0){3}} \put(0.25,1.25){\tiny #1} \put(1.25,1.25){\tiny #2} \put(2.25,1.25){\tiny #3} \put(0.25,0.25){\tiny #4} \put(1.25,0.25){\tiny #5}\end{picture}}
\end{picture}
}
\end{picture}
}

\newcommand{\VertexTableauThree}[7]{
\setlength{\unitlength}{1.5cm}
\begin{picture}(0,0)
\put(-0.25,0){
\begin{picture}(0,0)
\put(0,0){\circle*{0.1}} 
\put(#6,#7){\setlength{\unitlength}{0.3cm}\begin{picture}(0,0)\put(0,0){\line(0,1){1}} \put(1,0){\line(0,1){2}} \put(2,0){\line(0,1){2}} \put(3,0){\line(0,1){2}} \put(0,0){\line(1,0){3}} \put(0,1){\line(1,0){3}} \put(1,2){\line(1,0){2}} \put(0.25,0.25){\tiny #1} \put(1.25,1.25){\tiny #2} \put(1.25,0.25){\tiny #3} \put(2.25,1.25){\tiny #4} \put(2.25,0.25){\tiny #5}\end{picture}}
\end{picture}
}
\end{picture}
}

\newcommand{\VertexTableauFour}[6]{
\setlength{\unitlength}{1.5cm}
\begin{picture}(0,0)
\put(-0.25,0){
\begin{picture}(0,0)
\put(0,0){\circle*{0.1}} 
\put(#5,#6){\setlength{\unitlength}{0.3cm}\begin{picture}(0,0)\put(0,0){\line(0,1){2}} \put(1,0){\line(0,1){2}} \put(2,1){\line(0,1){1}} \put(3,1){\line(0,1){1}} \put(0,0){\line(1,0){1}} \put(0,1){\line(1,0){3}} \put(0,2){\line(1,0){3}} \put(0.25,1.25){\tiny #1} \put(0.25,0.25){\tiny #2} \put(1.25,1.25){\tiny #3} \put(2.25,1.25){\tiny #4}\end{picture}}
\end{picture}
}
\end{picture}
}

\newcommand{\VertexTableauForText}[5]{
\setlength{\unitlength}{1.5cm}
\begin{picture}(0.675,0)
\put(-0.05,-0.15){
\begin{picture}(0,0)
\put(0,0){\setlength{\unitlength}{0.3cm}\begin{picture}(0,0)\put(0,0){\line(0,1){2}} \put(1,0){\line(0,1){2}} \put(2,0){\line(0,1){2}} \put(3,1){\line(0,1){1}} \put(0,0){\line(1,0){2}} \put(0,1){\line(1,0){3}} \put(0,2){\line(1,0){3}} \put(0.25,1.25){\tiny #1} \put(1.25,1.25){\tiny #2} \put(2.25,1.25){\tiny #3} \put(0.25,0.25){\tiny #4} \put(1.25,0.25){\tiny #5}\end{picture}}
\end{picture}
}
\end{picture}
}

\begin{document}
\pagenumbering{arabic}
\thispagestyle{empty}%
\vspace*{-0.7in}
\noindent {\scriptsize \bf \em To appear in Advances in Applied Mathematics} \hfill {\scriptsize April 28, 2022}

\begin{center}
{\large \bf Gelfand--Tsetlin-type weight bases for all special linear Lie algebra\\ 
representations corresponding to skew Schur functions} 

\vspace*{0.05in}
\renewcommand{\thefootnote}{1}
Robert G.\ Donnelly\footnote{Department of Mathematics and Statistics, Murray State
University, Murray, KY 42071\\ 
\hspace*{0.25in}Email: {\tt rob.donnelly@murraystate.edu}} 
\renewcommand{\thefootnote}{2} 
\hspace*{-0.07in}and Molly W.\ Dunkum\footnote{Department of Mathematics, Western Kentucky University, Bowling Green, KY 42101\\ 
\hspace*{0.25in}Email: {\tt molly.dunkum@wku.edu}}

{\small \em Dedicated to Robert A.~(Bob) Proctor, in honor of his 65$^{\mbox{\tiny th}}$ birthday.}
\end{center} 

\begin{abstract}
We generalize the famous weight basis constructions of the finite-dimensional irreducible representations of $\mathfrak{sl}(n,\mathbb{C})$ obtained by Gelfand and Tsetlin in 1950. 
Using combinatorial methods, we construct one such basis for each finite-dimensional representation of $\mathfrak{sl}(n,\mathbb{C})$ associated to a given skew Schur function. 
Our constructions use diamond-colored distributive lattices of skew-shaped semistandard tableaux that generalize some classical Gelfand--Tsetlin (GT) lattices. 
Our constructions take place within the context of a certain programmatic study of poset models for semisimple Lie algebra representations and Weyl group symmetric functions undertaken by the first-named author and others. 
Some key aspects of the methodology of that program are recapitulated here. 

Combinatorial and representation-theoretic applications of our constructions are pursued here and elsewhere. 
Here, we extend combinatorial results about classical GT lattices to our more general lattices. 
We also use classical GT lattices to construct new and combinatorially distinctive weight bases for certain families of irreducible representations of the orthogonal Lie algebras. 
In another paper, via an entirely similar approach using the non-classical lattices of this paper, we obtain explicit weight bases for some infinite families of irreducible representations of the exceptional simple Lie algebras of types $\myE_{6}$ and $\myE_{7}$. 
Other combinatorial applications and generalizations are pursued in companion papers. 

\
\vspace*{0.05in}

\begin{center}

\end{center} 
\end{abstract}

\noindent {\bf \S \IntroSection\ Introduction}. 
The origin of what is now called Gelfand--Tsetlin theory is the explicit construction, presented in 1950 by I.\ Gelfand and M.\ Tsetlin \cite{GT}, of weight bases for all finite-dimensional irreducible representations of the general linear Lie algebras. 
Their originating idea has since been generalized in a number of ways: To the combinatorial study of weight bases for irreducible simple Lie algebra representations, e.g.\ \cite{PrGZ}, \cite{DonSupp}; to the discovery, by advanced algebraic methods, of GT-type weight bases for irreducible representations of classical simple Lie algebras in types $\myB$, $\myC$, and $\myD$, see \cite{MoLast}, \cite{MoSurvey}; to the study of $q$-deformations of simple Lie algebras, called quantum groups \cite{Cherednik}; and to the study of infinite-dimensional $\mathfrak{gl}_{n}$-modules --- called Gelfand--Tsetlin modules --- admitting GT-type weight bases, which originated in \cite{DOF} and \cite{DFO} and continues as part of the more general theory of Galois orders, see e.g.\ \cite{FGRZ}, \cite{Har}, and references therein. 

\begin{figure}[t]
\begin{center}
\IntroFig: The skew-tabular lattice $L_{\mytinyA_{2}}^{\mbox{\tiny skew}}({\setlength{\unitlength}{0.125cm}\begin{picture}(3,0)\put(0,0){\line(0,1){1}} \put(1,0){\line(0,1){1}} \put(2,0){\line(0,1){2}} \put(3,0){\line(0,1){2}} \put(0,0){\line(1,0){3}} \put(0,1){\line(1,0){3}} \put(2,2){\line(1,0){1}}\end{picture}})$, built from skew-shaped tableaux.  

\setlength{\unitlength}{1.5cm}
\begin{picture}(4,6.5)
\put(1,0){\line(-1,1){1}}
\put(1,0){\line(1,1){2}}
\put(0,1){\line(1,1){2}}
\put(2,1){\line(-1,1){1}}
\put(2,1){\line(0,1){1}}
\put(1,2){\line(0,1){1}}
\put(2,2){\line(-1,1){2}}
\put(2,2){\line(1,1){2}}
\put(3,2){\line(-1,1){1}}
\put(3,2){\line(0,1){1}}
\put(1,3){\line(1,1){2}}
\put(2,3){\line(0,1){1}}
\put(3,3){\line(-1,1){2}}
\put(0,4){\line(1,1){2}}
\put(4,4){\line(-1,1){2}}
\put(2,6){\VertexTableau{1}{1}{1}{2}{-0.7}{0,05}}
\put(1,5){\VertexTableau{1}{2}{1}{2}{-0.7}{0.05}}
\put(3,5){\VertexTableau{1}{1}{1}{3}{0.15}{-0.05}}
\put(0,4){\VertexTableau{2}{2}{1}{2}{-0.75}{-0.05}}
\put(2,4){\VertexTableau{1}{2}{1}{3}{0.25}{-0.15}}
\put(4,4){\VertexTableau{1}{1}{2}{3}{0.15}{-0.05}}
\put(1,3){\VertexTableau{2}{2}{1}{3}{-0.8}{-0.3}}
\put(2,3){\VertexTableau{1}{3}{1}{3}{-0.75}{-0.1}}
\put(3,3){\VertexTableau{1}{2}{2}{3}{0.15}{-0.3}}
\put(1,2){\VertexTableau{2}{3}{1}{3}{-0.8}{-0.1}}
\put(2,2){\VertexTableau{2}{2}{2}{3}{0.15}{-0.2}}
\put(3,2){\VertexTableau{1}{3}{2}{3}{0.15}{-0.2}}
\put(0,1){\VertexTableau{3}{3}{1}{3}{-0.75}{-0.25}}
\put(2,1){\VertexTableau{2}{3}{2}{3}{0.15}{-0.25}}
\put(1,0){\VertexTableau{3}{3}{2}{3}{0.15}{-0.25}}
\put(1,5){\NEEdgeLabelForLatticeI{{\em 1}}}
\put(3,5){\NWEdgeLabelForLatticeI{{\em 2}}}
\put(0,4){\NEEdgeLabelForLatticeI{{\em 1}}}
\put(2,4){\NWEdgeLabelForLatticeI{{\em 2}}}
\put(2,4){\NEEdgeLabelForLatticeI{{\em 1}}}
\put(4,4){\NWEdgeLabelForLatticeI{{\em 1}}}
\put(1,3){\NEEdgeLabelForLatticeI{{\em 1}}}
\put(1,3){\NWEdgeLabelForLatticeI{{\em 2}}}
\put(2,3){\VerticalEdgeLabelForLatticeI{{\em 2}}}
\put(3,3){\NWEdgeLabelForLatticeI{{\em 1}}}
\put(3,3){\NEEdgeLabelForLatticeI{{\em 1}}}
\put(1,2){\VerticalEdgeLabelForLatticeI{{\em 2}}}
\put(1.25,2.25){\NEEdgeLabelForLatticeI{{\em 1}}}
\put(2.2,1.8){\NWEdgeLabelForLatticeI{{\em 1}}}
\put(1.8,1.8){\NEEdgeLabelForLatticeI{{\em 1}}}
\put(3,2){\VerticalEdgeLabelForLatticeI{{\em 2}}}
\put(2.75,2.25){\NWEdgeLabelForLatticeI{{\em 1}}}
\put(0,1){\NEEdgeLabelForLatticeI{{\em 2}}}
\put(2,1){\VerticalEdgeLabelForLatticeI{{\em 2}}}
\put(2,1){\NWEdgeLabelForLatticeI{{\em 1}}}
\put(2,1){\NEEdgeLabelForLatticeI{{\em 1}}}
\put(1,0){\NWEdgeLabelForLatticeI{{\em 1}}}
\put(1,0){\NEEdgeLabelForLatticeI{{\em 2}}}
\end{picture}
\end{center}
\end{figure}

Here, however, we will return to the original setting of Gelfand and Tsetlin and produce explicit weight bases for representations of the special linear Lie algebras $\mathfrak{sl}(n,\mathbb{C})$ associated with the skew Schur functions, which are symmetric functions in the traditional sense cf.\ \cite{StanText2} or \cite{vL}. 
Specifically, for each skew-shape $\mysmallP/\mysmallQ$ obtained from `skew-compatible' 
partitions\myfootnote{Partitions $\mysmallP$ and $\mysmallQ$ are to be viewed as weakly decreasing nonnegative integer $m$-tuples $\mysmallP = (\myscriptsizeP_{1},\ldots,\myscriptsizeP_{m})$ and $\mysmallQ = (\myscriptsizeQ_{1},\ldots,\myscriptsizeQ_{m})$ where $m \geq n$ and such that the skew shape $\mysmallP/\mysmallQ$ has no column with more than $n$ boxes. This latter constraint guarantees that there are no unrealizable columns when we consider semistandard (and therefore column-strict) tableaux with skew shape $\mysmallP/\mysmallQ$ and entries from the set $\{1,2,\ldots,n\}$.} 
$\mysmallP$ and $\mysmallQ$, we explicitly construct a weight basis for a representation of $\mathfrak{sl}(n,\mathbb{C})$ whose associated character --- what we refer to as a Weyl symmetric function --- is the skew Schur function here denoted $\vartheta_{_{\mytinyP/\mytinyQ}}$. 
Each such weight basis is obtained from a particular diamond-colored distributive lattice naturally built from semistandard tableaux of the appropriate skew shape, which leads to some interesting combinatorics. 
The distributive lattice associated to $\vartheta_{_{\mytinyP/\mytinyQ}}$ is here denoted 
$L_{\mytinyA_{n-1}}^{\mbox{\tiny skew}}(\mysmallP/\mysmallQ)$ and called a `skew-tabular lattice'  (see \IntroFig).  
When $\mysmallQ = (0,\ldots,0)$, so the shape $\mysmallP/\mysmallQ$ is not strictly skewed, the `classical Gelfand--Tsetlin (GT) lattice' $L_{\mytinyA_{n-1}}^{\mbox{\tiny GT}}(\mysmallP)$ is known to realize the irreducible $\mathfrak{sl}(n,\mathbb{C})$-representation whose character is the Schur function $\vartheta_{_{\mytinyP}}$; all irreducible $\mathfrak{sl}(n,\mathbb{C})$-representations can be realized in this way \cite{PrGZ}, \cite{DonSupp}, \cite{HL}. 
Henceforth, the phrase `skew-tabular lattice' will generically refer to strictly skew-tabular or classical Gelfand--Tsetlin lattices. 

\begin{figure}[t]
\begin{center}
\IntroFigTwo:  The skew-tabular lattice $L_{\mytinyA_{2}}^{\mbox{\tiny skew}}({\setlength{\unitlength}{0.125cm}\begin{picture}(3,0)\put(0,-0.75){\line(0,1){2}} \put(1,-0.75){\line(0,1){2}} \put(2,0.25){\line(0,1){2}} \put(3,0.25){\line(0,1){2}} \put(0,-0.75){\line(1,0){1}} \put(0,0.25){\line(1,0){3}} \put(0,1.25){\line(1,0){3}} \put(2,2.25){\line(1,0){1}}\end{picture}})$.  

\setlength{\unitlength}{1.5cm}
\begin{picture}(4,6.5)
\put(0,0){\qbezier(2,6)(3.25,5)(4.5,4)}
\put(0,0){\qbezier(2,6)(1.5,5.5)(1,5)}
\put(0,0){\qbezier(2,6)(2,5.5)(2,5)}
\put(0,0){\qbezier(3.25,5)(2.75,4.5)(2.25,4)}
\put(0,0){\qbezier(3.25,5)(3.25,4)(3.25,3)}
\put(0,0){\qbezier(4.5,4)(4,3.5)(3.5,3)}
\put(0,0){\qbezier(4.5,4)(4.5,3)(4.5,2)}
\put(0,0){\qbezier(1,5)(2.25,4)(3.5,3)}
\put(0,0){\qbezier(1,5)(1,4.5)(1,4)}
\put(0,0){\qbezier(2,5)(3.25,4)(4.5,3)}
\put(0,0){\qbezier(2,5)(1,4)(0,3)}
\put(0,0){\qbezier(1,4)(2.25,3)(3.5,2)}
\put(0,0){\qbezier(2.25,4)(2.25,3)(2.25,2)}
\put(0,0){\qbezier(3.5,3)(3.5,2)(3.5,1)}
\put(0,0){\qbezier(3.25,4)(2.25,3)(1.25,2)}
\put(0,0){\qbezier(4.5,3)(3.5,2)(2.5,1)}
\put(0,0){\qbezier(0,3)(1.25,2)(2.5,1)}
\put(0,0){\qbezier(1.25,2)(1.25,1.5)(1.25,1)}
\put(0,0){\qbezier(2.5,1)(2.5,0.5)(2.5,0)}
\put(0,0){\qbezier(3.25,3)(2.25,2)(1.25,1)}
\put(0,0){\qbezier(3.25,3)(3.875,2.5)(4.5,2)}
\put(0,0){\qbezier(4.5,2)(3.5,1)(2.5,0)}
\put(0,0){\qbezier(2.25,2)(2.875,1.5)(3.5,1)}
\put(0,0){\qbezier(1.25,1)(1.875,0.5)(2.5,0)}
\put(2,6){\VertexTableauFib{1}{2}{1}{1}{2}{-0.7}{0.05}}
\put(3.25,5){\VertexTableauFib{1}{2}{1}{1}{3}{0.15}{0.1}}
\put(4.5,4){\VertexTableauFib{1}{2}{1}{2}{3}{0.15}{-0.05}}
\put(2,5){\VertexTableauFib{1}{2}{2}{1}{2}{-0.6}{0.05}}
\put(3.25,4){\VertexTableauFib{1}{2}{2}{1}{3}{0.1}{0.1}}
\put(4.5,3){\VertexTableauFib{1}{2}{2}{2}{3}{0.15}{0.05}}
\put(3.25,3){\VertexTableauFib{1}{2}{3}{1}{3}{-0.625}{-0.05}}
\put(4.5,2){\VertexTableauFib{1}{2}{3}{2}{3}{0.15}{0.05}}
\put(1,5){\VertexTableauFib{1}{3}{1}{1}{2}{-0.6}{0.05}}
\put(2.25,4){\VertexTableauFib{1}{3}{1}{1}{3}{-0.6}{-0.3}}
\put(3.5,3){\VertexTableauFib{1}{3}{1}{2}{3}{0.12}{-0.1}}
\put(1,4){\VertexTableauFib{1}{3}{2}{1}{2}{-0.6}{0.05}}
\put(2.25,3){\VertexTableauFib{1}{3}{2}{1}{3}{-0.775}{-0.15}}
\put(3.5,2){\VertexTableauFib{1}{3}{2}{2}{3}{0.125}{-0.2}}
\put(2.25,2){\VertexTableauFib{1}{3}{3}{1}{3}{-0.625}{-0.05}}
\put(3.5,1){\VertexTableauFib{1}{3}{3}{2}{3}{0.15}{-0.1}}
\put(0,3){\VertexTableauFib{2}{3}{2}{1}{2}{-0.6}{0.05}}
\put(1.25,2){\VertexTableauFib{2}{3}{2}{1}{3}{-0.625}{-0.3}}
\put(2.5,1){\VertexTableauFib{2}{3}{2}{2}{3}{-0.625}{-0.3}}
\put(1.25,1){\VertexTableauFib{2}{3}{3}{1}{3}{-0.625}{-0.05}}
\put(2.5,0){\VertexTableauFib{2}{3}{3}{2}{3}{0.15}{-0.05}}
\put(1,5){\NEEdgeLabelForLatticeI{{\em 2}}}
\put(2,3.75){\NEEdgeLabelForLatticeI{{\em 2}}}
\put(2.1,2.85){\NEEdgeLabelForLatticeI{{\em 2}}}
\put(2.1,1.85){\NEEdgeLabelForLatticeI{{\em 2}}}
\put(3.75,3.25){\NEEdgeLabelForLatticeI{{\em 2}}}
\put(3.75,2.25){\NEEdgeLabelForLatticeI{{\em 2}}}
\put(3.6,1.1){\NEEdgeLabelForLatticeI{{\em 2}}}
\put(1.3,4.3){\NEEdgeLabelForLatticeI{{\em 2}}}
\put(2,5){\VerticalEdgeLabelForLatticeI{{\em 1}}}
\put(1,4){\VerticalEdgeLabelForLatticeI{{\em 1}}}
\put(3.25,4){\VerticalEdgeLabelForLatticeI{{\em 1}}}
\put(2.25,3){\VerticalEdgeLabelForLatticeI{{\em 1}}}
\put(4.5,3){\VerticalEdgeLabelForLatticeI{{\em 1}}}
\put(3.5,2){\VerticalEdgeLabelForLatticeI{{\em 1}}}
\put(3.25,3){\VerticalEdgeLabelForLatticeI{{\em 2}}}
\put(2.25,2){\VerticalEdgeLabelForLatticeI{{\em 2}}}
\put(1.25,1){\VerticalEdgeLabelForLatticeI{{\em 2}}}
\put(4.5,2){\VerticalEdgeLabelForLatticeI{{\em 2}}}
\put(3.5,1){\VerticalEdgeLabelForLatticeI{{\em 2}}}
\put(2.5,0){\VerticalEdgeLabelForLatticeI{{\em 2}}}
\put(3.15,5){\NWEdgeLabelForLatticeI{{\em 2}}}
\put(2.4,3.8){\NWEdgeLabelForLatticeI{{\em 2}}}
\put(2.9,4.2){\NWEdgeLabelForLatticeI{{\em 2}}}
\put(2.025,3.1){\NWEdgeLabelForLatticeI{{\em 2}}}
\put(1.025,2.1){\NWEdgeLabelForLatticeI{{\em 2}}}
\put(4.525,3.9){\NWEdgeLabelForLatticeI{{\em 1}}}
\put(4.15,3.2){\NWEdgeLabelForLatticeI{{\em 1}}}
\put(4.65,1.8){\NWEdgeLabelForLatticeI{{\em 1}}}
\put(3.525,1.9){\NWEdgeLabelForLatticeI{{\em 1}}}
\put(2.1,3.3){\NEEdgeLabelForLatticeI{{\em 1}}}
\put(2.1,1.3){\NEEdgeLabelForLatticeI{{\em 1}}}
\put(1.6,0.9){\NEEdgeLabelForLatticeI{{\em 1}}}
\put(1.225,0.2){\NEEdgeLabelForLatticeI{{\em 1}}}
\put(0.1,3.1){\NEEdgeLabelForLatticeI{{\em 1}}}
\put(1.25,2){\NEEdgeLabelForLatticeI{{\em 1}}}
\put(2.45,1.3){\NWEdgeLabelForLatticeI{{\em 1}}}
\put(2.3,0.8){\NEEdgeLabelForLatticeI{{\em 1}}}
\put(3.5,0.1){\NWEdgeLabelForLatticeI{{\em 1}}}
\end{picture}
\end{center}
\end{figure}

Before we discuss some general features of our weight basis constructions, we mention some overall motivation. 
A primary goal of our related work in papers such as \cite{DonSymp}, \cite{DonAdjoint}, \cite{DLP1}, and \cite{ADLPOne}/\cite{ADLMPPW} has been 
to find ranked posets, and in particular diamond-colored modular and distributive lattices, that can serve as models for (weight bases of) semisimple Lie algebra representations and/or their associated Weyl symmetric functions. 
Such poset models for semisimple Lie algebra representations are herein called `supporting graphs' and for Weyl symmetric functions are called `splitting posets'. 
(For some general discourses on these ideas, see \cite{DonPosetModels}, \cite{DonDiamond}.) 
Often, the representations of interest are irreducible, and their associated Weyl symmetric functions can be expressed as quotients of certain alternants. 
So, these latter `Weyl bialternants' are general and type-independent versions of the Schur functions associated with the irreducible representations of the type $\myA_{n-1}$ simple Lie algebra $\mathfrak{g}(\myA_{n-1}) \cong \mathfrak{sl}(n,\mathbb{C})$. 

Splitting posets for Weyl bialternants and supporting graphs for irreducible semisimple Lie algebra representations have many nice combinatorial features (see e.g.\ \MainCombinatorialTheorem\ below). 
An interesting historical example is that Stanley demonstrated in \cite{StanGLn} some lovely combinatorial properties of the classical GT lattices essentially by understanding them as splitting posets for Schur functions; in \cite{PrGZ}, Proctor recognized that the lattices Stanley had considered could be naturally viewed as supporting graphs for the GT bases, a fact that inspires our use of the name `classical GT lattices'. 
So, such observed instances of supporting graphs and splitting posets afford combinatorial connections and results. 
Likewise, the consideration of certain nice ordered structures can afford algebraic results and insights.  
Starting with well-chosen ranked posets and lattices, we have, on occasion, obtained new explicit irreducible representation constructions and gained otherwise difficult-to-discern insights into their combinatorial structure (see e.g.\ \cite{DonSupp}, \cite{DonAdjoint}, \cite{DLP2}, \cite{ADLPOne}). 

Reducible representations and Weyl symmetric functions are an important part of this overall project. 
At times, we have discovered distributive lattice supporting graphs/splitting posets of irreducible representations/Weyl bialternants as sublattices of larger lattices, the latter of which naturally realize reducible representations/Weyl symmetric functions, e.g.\ \cite{DonSymp} or \cite{DLP1}. 
Also, components of supporting graphs and splitting posets for irreducible representations/Weyl bialternants are often not irreducible, and analyzing such components can be crucial (\cite{DonSymp}, \cite{DDW}). 
Supporting graphs for irreducible $\mathfrak{g}(\myA_{1})$-modules are simply chains, but, arguably, the most interesting and most frequently occurring supporting graphs for $\mathfrak{g}(\myA_{1})$-modules are for reducible representations: Consider that any connected supporting graph for a semisimple Lie algebra representation realizes a (typically) reducible representation of $\mathfrak{g}(\myA_{1})$ via the induced action of a principal three-dimensional subalgebra, an idea that is largely responsible for \MainCombinatorialTheorem\ below via Proctor's `Peck Poset Theorem' \cite{PrPeck}. 

Diamond-colored modular/distributive lattice supporting graphs provide an answer to a problem posed by Stanley (see Problem 3 of \cite{StanUnim}) to find distributive lattices, or other ranked posets, whose rank-generating functions are of the form presented in \MainCombinatorialTheorem\ of this paper. 
This is one motivation for identifying weight bases that afford combinatorially distinctive supporting graphs. 
Supporting graphs can also be used to define several other properties of interest concerning weight bases of semisimple Lie algebra representations.\footnote{These and other combinatorial notions related to semisimple Lie algebra representations, weight bases, and the like are developed more precisely in \S \SetupSectionTwo\ below.} 
We say a weight basis is edge-minimal if no proper subgraph of its supporting graph supports another weight basis for the same representation, and the weight basis is edge-minimizing if no weight basis has a supporting graph with fewer edges. 
The weight basis is solitary if the only other weight bases which share the same supporting graph are those obtained by diagonal change-of-basis matrices. 
The adjectives edge-minimal, edge-minimizing, and solitary apply to supporting graphs as well. 
The fact that, in many special cases we have studied, diamond-colored modular/distributive lattice supporting graphs often evince solitarity or edge-minimality and seem finely-tuned to the problem of explicitly constructing weight bases for semisimple Lie algebra representations or serving as splitting posets for Weyl symmetric functions is, for now, a phenomenon that has no fully general theoretical foundation. 
Ultimately, we hope the present paper might serve as a useful landmark on the journey toward better understanding this phenomenon. 

A typical skew-tabular lattice is depicted in \IntroFigTwo. 
We will see in \GTGeneralTheorem\ that the 21-element lattice $L_{\mytinyA_{2}}^{\mbox{\tiny skew}}({\setlength{\unitlength}{0.125cm}\begin{picture}(3,0)\put(0,-0.75){\line(0,1){2}} \put(1,-0.75){\line(0,1){2}} \put(2,0.25){\line(0,1){2}} \put(3,0.25){\line(0,1){2}} \put(0,-0.75){\line(1,0){1}} \put(0,0.25){\line(1,0){3}} \put(0,1.25){\line(1,0){3}} \put(2,2.25){\line(1,0){1}}\end{picture}})$ is a supporting graph for some $\mathfrak{g}(\myA_{2})$-module we notate as `$V[L_{\mytinyA_{2}}^{\mbox{\tiny skew}}({\setlength{\unitlength}{0.125cm}\begin{picture}(3,0)\put(0,-0.75){\line(0,1){2}} \put(1,-0.75){\line(0,1){2}} \put(2,0.25){\line(0,1){2}} \put(3,0.25){\line(0,1){2}} \put(0,-0.75){\line(1,0){1}} \put(0,0.25){\line(1,0){3}} \put(0,1.25){\line(1,0){3}} \put(2,2.25){\line(1,0){1}}\end{picture}})]$'. 
By the Littlewood--Richrdson Rule (as recounted in, say, Theorem 1.4.4 of \cite{vL} or in \LRCorollaryOne\ below), the associated skew Schur function is expressible as $\vartheta_{\setlength{\unitlength}{0.075cm}\begin{picture}(3,0)\put(0,-0.75){\line(0,1){2}} \put(1,-0.75){\line(0,1){2}} \put(2,0.25){\line(0,1){2}} \put(3,0.25){\line(0,1){2}} \put(0,-0.75){\line(1,0){1}} \put(0,0.25){\line(1,0){3}} \put(0,1.25){\line(1,0){3}} \put(2,2.25){\line(1,0){1}}\end{picture}} = \vartheta_{\setlength{\unitlength}{0.075cm}\begin{picture}(3,0)\put(0,-0.25){\line(0,1){2}} \put(1,-0.25){\line(0,1){2}} \put(2,-0.25){\line(0,1){2}} \put(3,0.75){\line(0,1){1}} \put(0,-0.25){\line(1,0){2}} \put(0,0.75){\line(1,0){3}} \put(0,1.75){\line(1,0){3}}\end{picture}} + \vartheta_{\setlength{\unitlength}{0.075cm}\begin{picture}(2,0)\put(0,0){\line(0,1){1}} \put(1,0){\line(0,1){1}} \put(2,0){\line(0,1){1}} \put(0,0){\line(1,0){2}} \put(0,1){\line(1,0){2}}\end{picture}}$, and therefore $V[L_{\mytinyA_{2}}^{\mbox{\tiny skew}}({\setlength{\unitlength}{0.125cm}\begin{picture}(3,0)\put(0,-0.75){\line(0,1){2}} \put(1,-0.75){\line(0,1){2}} \put(2,0.25){\line(0,1){2}} \put(3,0.25){\line(0,1){2}} \put(0,-0.75){\line(1,0){1}} \put(0,0.25){\line(1,0){3}} \put(0,1.25){\line(1,0){3}} \put(2,2.25){\line(1,0){1}}\end{picture}})]$ decomposes as the direct sum $V[L_{\mytinyA_{2}}^{\mbox{\tiny GT}}({\setlength{\unitlength}{0.125cm}\begin{picture}(3,0)\put(0,-0.25){\line(0,1){2}} \put(1,-0.25){\line(0,1){2}} \put(2,-0.25){\line(0,1){2}} \put(3,0.75){\line(0,1){1}} \put(0,-0.25){\line(1,0){2}} \put(0,0.75){\line(1,0){3}} \put(0,1.75){\line(1,0){3}}\end{picture}})] \oplus V[L_{\mytinyA_{2}}^{\mbox{\tiny GT}}({\setlength{\unitlength}{0.125cm}\begin{picture}(2,0)\put(0,0.25){\line(0,1){1}} \put(1,0.25){\line(0,1){1}} \put(2,0.25){\line(0,1){1}} \put(0,0.25){\line(1,0){2}} \put(0,1.25){\line(1,0){2}}\end{picture}})]$ of irreducible $\mathfrak{g}(\myA_{2})$-modules, where the disjoint sum $L_{\mytinyA_{2}}^{\mbox{\tiny GT}}({\setlength{\unitlength}{0.125cm}\begin{picture}(3,0)\put(0,-0.25){\line(0,1){2}} \put(1,-0.25){\line(0,1){2}} \put(2,-0.25){\line(0,1){2}} \put(3,0.75){\line(0,1){1}} \put(0,-0.25){\line(1,0){2}} \put(0,0.75){\line(1,0){3}} \put(0,1.75){\line(1,0){3}}\end{picture}}) \oplus L_{\mytinyA_{2}}^{\mbox{\tiny GT}}({\setlength{\unitlength}{0.125cm}\begin{picture}(2,0)\put(0,0.25){\line(0,1){1}} \put(1,0.25){\line(0,1){1}} \put(2,0.25){\line(0,1){1}} \put(0,0.25){\line(1,0){2}} \put(0,1.25){\line(1,0){2}}\end{picture}})$ of classical GT lattices  
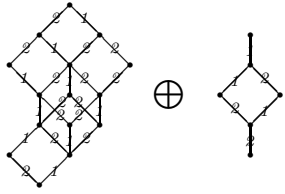
\begin{center}
{\setlength{\unitlength}{0.4cm}\begin{picture}(9,6)
\put(1,0){\circle*{0.15}} \put(0,1){\circle*{0.15}} \put(2,1){\circle*{0.15}} \put(1,2){\circle*{0.15}} 
\put(2,2){\circle*{0.15}} \put(3,2){\circle*{0.15}} \put(1,3){\circle*{0.15}} \put(2,3){\circle*{0.15}}
\put(3,3){\circle*{0.15}} \put(0,4){\circle*{0.15}} \put(2,4){\circle*{0.15}} \put(4,4){\circle*{0.15}}
\put(1,5){\circle*{0.15}} \put(3,5){\circle*{0.15}} \put(2,6){\circle*{0.15}} 
\put(1,0){\line(1,1){2}} \put(0,1){\line(1,1){2}} 
\put(0,4){\line(1,1){2}} \put(1,3){\line(1,1){2}} \put(2,2){\line(1,1){2}} 
\put(1,0){\line(-1,1){1}} \put(2,1){\line(-1,1){1}} \put(3,2){\line(-1,1){1}} 
\put(2,2){\line(-1,1){2}} \put(3,3){\line(-1,1){2}} \put(4,4){\line(-1,1){2}}  
\put(2,1){\line(0,1){1}} \put(1,2){\line(0,1){1}} \put(3,2){\line(0,1){1}} \put(2,3){\line(0,1){1}} 
\put(1.3,0.3){\tiny \em 1} 
\put(0.35,1.4){\tiny \em 1} 
\put(0.35,0.3){\tiny \em 2} 
\put(1.3,1.4){\tiny \em 2} 
\put(2.35,1.4){\tiny \em 2} 
\put(1.35,3.4){\tiny \em 2} 
\put(0.35,4.4){\tiny \em 2} 
\put(3.35,3.4){\tiny \em 2} 
\put(2.35,4.4){\tiny \em 2} 
\put(1.35,5.4){\tiny \em 2} 
\put(2.3,5.4){\tiny \em 1} 
\put(1.3,4.4){\tiny \em 1} 
\put(0.3,3.4){\tiny \em 1} 
\put(2.1,2.15){\tiny \em 2} 
\put(1.55,2.6){\tiny \em 2} 
\put(2.1,2.6){\tiny \em 2} 
\put(1.5,2.15){\tiny \em 2} 
\put(2.3,3.4){\tiny \em 2} 
\put(3.3,4.4){\tiny \em 2} 
\put(0.8,2.3){\tiny \em 1} 
\put(2.8,2.3){\tiny \em 1} 
\put(1.8,1.3){\tiny \em 1} 
\put(1.8,3.3){\tiny \em 1} 
\put(4.75,2.75){$\bigoplus$}
\put(8,1){\circle*{0.15}} \put(8,2){\circle*{0.15}} \put(7,3){\circle*{0.15}} \put(9,3){\circle*{0.15}} \put(8,4){\circle*{0.15}} \put(8,5){\circle*{0.15}} 
\put(8,1){\line(0,1){1}} \put(8,4){\line(0,1){1}} 
\put(8,2){\line(-1,1){1}} \put(8,2){\line(1,1){1}} \put(7,3){\line(1,1){1}} \put(9,3){\line(-1,1){1}} 
\put(7.8,4.3){\tiny \em 1} 
\put(7.8,1.3){\tiny \em 2} 
\put(8.3,3.4){\tiny \em 2} 
\put(7.3,2.4){\tiny \em 2} 
\put(8.3,2.3){\tiny \em 1} 
\put(7.3,3.3){\tiny \em 1} 
\end{picture}}
\end{center}
is the supporting graph for the weight basis obtained by joining together the GT bases within each of the two irreducible components. 
By Theorem 4.4 of \cite{DonSupp} and Theorem 5.5 of \cite{HL}, classical GT lattices are solitary and edge-minimal. 
Since the decomposition $\vartheta_{\setlength{\unitlength}{0.075cm}\begin{picture}(3,0)\put(0,-0.75){\line(0,1){2}} \put(1,-0.75){\line(0,1){2}} \put(2,0.25){\line(0,1){2}} \put(3,0.25){\line(0,1){2}} \put(0,-0.75){\line(1,0){1}} \put(0,0.25){\line(1,0){3}} \put(0,1.25){\line(1,0){3}} \put(2,2.25){\line(1,0){1}}\end{picture}} = \vartheta_{\setlength{\unitlength}{0.075cm}\begin{picture}(3,0)\put(0,-0.25){\line(0,1){2}} \put(1,-0.25){\line(0,1){2}} \put(2,-0.25){\line(0,1){2}} \put(3,0.75){\line(0,1){1}} \put(0,-0.25){\line(1,0){2}} \put(0,0.75){\line(1,0){3}} \put(0,1.75){\line(1,0){3}}\end{picture}} + \vartheta_{\setlength{\unitlength}{0.075cm}\begin{picture}(2,0)\put(0,0){\line(0,1){1}} \put(1,0){\line(0,1){1}} \put(2,0){\line(0,1){1}} \put(0,0){\line(1,0){2}} \put(0,1){\line(1,0){2}}\end{picture}}$ of the skew Schur function $\vartheta_{\setlength{\unitlength}{0.075cm}\begin{picture}(3,0)\put(0,-0.75){\line(0,1){2}} \put(1,-0.75){\line(0,1){2}} \put(2,0.25){\line(0,1){2}} \put(3,0.25){\line(0,1){2}} \put(0,-0.75){\line(1,0){1}} \put(0,0.25){\line(1,0){3}} \put(0,1.25){\line(1,0){3}} \put(2,2.25){\line(1,0){1}}\end{picture}}$ is `multiplicity-free' in the sense of \cite{GHO}, it follows that the supporting graph $L_{\mytinyA_{2}}^{\mbox{\tiny GT}}({\setlength{\unitlength}{0.125cm}\begin{picture}(3,0)\put(0,-0.25){\line(0,1){2}} \put(1,-0.25){\line(0,1){2}} \put(2,-0.25){\line(0,1){2}} \put(3,0.75){\line(0,1){1}} \put(0,-0.25){\line(1,0){2}} \put(0,0.75){\line(1,0){3}} \put(0,1.75){\line(1,0){3}}\end{picture}}) \oplus L_{\mytinyA_{2}}^{\mbox{\tiny GT}}({\setlength{\unitlength}{0.125cm}\begin{picture}(2,0)\put(0,0.25){\line(0,1){1}} \put(1,0.25){\line(0,1){1}} \put(2,0.25){\line(0,1){1}} \put(0,0.25){\line(1,0){2}} \put(0,1.25){\line(1,0){2}}\end{picture}})$ depicted immediately above is solitary and edge-minimal. 
Now, the number of edges in $L_{\mytinyA_{2}}^{\mbox{\tiny skew}}({\setlength{\unitlength}{0.125cm}\begin{picture}(3,0)\put(0,-0.75){\line(0,1){2}} \put(1,-0.75){\line(0,1){2}} \put(2,0.25){\line(0,1){2}} \put(3,0.25){\line(0,1){2}} \put(0,-0.75){\line(1,0){1}} \put(0,0.25){\line(1,0){3}} \put(0,1.25){\line(1,0){3}} \put(2,2.25){\line(1,0){1}}\end{picture}})$ exceeds the number of edges in $L_{\mytinyA_{2}}^{\mbox{\tiny GT}}({\setlength{\unitlength}{0.125cm}\begin{picture}(3,0)\put(0,-0.25){\line(0,1){2}} \put(1,-0.25){\line(0,1){2}} \put(2,-0.25){\line(0,1){2}} \put(3,0.75){\line(0,1){1}} \put(0,-0.25){\line(1,0){2}} \put(0,0.75){\line(1,0){3}} \put(0,1.75){\line(1,0){3}}\end{picture}}) \oplus L_{\mytinyA_{2}}^{\mbox{\tiny GT}}({\setlength{\unitlength}{0.125cm}\begin{picture}(2,0)\put(0,0.25){\line(0,1){1}} \put(1,0.25){\line(0,1){1}} \put(2,0.25){\line(0,1){1}} \put(0,0.25){\line(1,0){2}} \put(0,1.25){\line(1,0){2}}\end{picture}})$, so the former cannot be edge-minimizing. 
However, $L_{\mytinyA_{2}}^{\mbox{\tiny skew}}({\setlength{\unitlength}{0.125cm}\begin{picture}(3,0)\put(0,-0.75){\line(0,1){2}} \put(1,-0.75){\line(0,1){2}} \put(2,0.25){\line(0,1){2}} \put(3,0.25){\line(0,1){2}} \put(0,-0.75){\line(1,0){1}} \put(0,0.25){\line(1,0){3}} \put(0,1.25){\line(1,0){3}} \put(2,2.25){\line(1,0){1}}\end{picture}})$ is also solitary and edge-minimal; this can be checked by hand or seen as a special case of Proposition 6.1 of \cite{DDMN}. 
So, the skew-tabular lattice $L_{\mytinyA_{2}}^{\mbox{\tiny skew}}({\setlength{\unitlength}{0.125cm}\begin{picture}(3,0)\put(0,-0.75){\line(0,1){2}} \put(1,-0.75){\line(0,1){2}} \put(2,0.25){\line(0,1){2}} \put(3,0.25){\line(0,1){2}} \put(0,-0.75){\line(1,0){1}} \put(0,0.25){\line(1,0){3}} \put(0,1.25){\line(1,0){3}} \put(2,2.25){\line(1,0){1}}\end{picture}})$ achieves the distinctions of edge-minimality and solitarity while also being connected. 
For another example, one can check by hand that the 55-element lattice $L := L_{\mytinyA_{2}}^{\mbox{\tiny skew}}({\setlength{\unitlength}{0.125cm}\begin{picture}(4,0)\put(0,-0.75){\line(0,1){2}} \put(1,-0.75){\line(0,1){2}} \put(2,0.25){\line(0,1){2}} \put(3,0.25){\line(0,1){2}} \put(0,-0.75){\line(1,0){1}} \put(0,0.25){\line(1,0){3}} \put(0,1.25){\line(1,0){4}} \put(2,2.25){\line(1,0){2}} \put(4,1.25){\line(0,1){1}}\end{picture}})$ with rank generating function $\RGF(L,q) = q^{8}+4q^{7}+7q^{6}+10q^{5}+11q^{4}+10q^{3}+7q^{2}+4q+1$ is also a solitary and edge-minimal supporting graph, and its associated skew Schur function $\vartheta_{\setlength{\unitlength}{0.075cm}\begin{picture}(4,0)\put(0,-0.75){\line(0,1){2}} \put(1,-0.75){\line(0,1){2}} \put(2,0.25){\line(0,1){2}} \put(3,0.25){\line(0,1){2}} \put(0,-0.75){\line(1,0){1}} \put(0,0.25){\line(1,0){3}} \put(0,1.25){\line(1,0){4}} \put(2,2.25){\line(1,0){2}} \put(4,1.25){\line(0,1){1}}\end{picture}} = \vartheta_{(4,2)}+\vartheta_{(3,3)}+\vartheta_{(3,0)}+\vartheta_{(2,1)}$ is also multiplicity-free.

Consider now a generic skew Schur function $\vartheta_{_{\mytinyP/\mytinyQ}} = \vartheta_{_{\mytinyP_{1}}}+\cdots+\vartheta_{_{\mytinyP_{k}}}$ expressed as a sum of Schur functions. 
The previous paragraph demonstrates that, as a supporting graph, the (connected) skew-tabular lattice $L_{\mytinyA_{n-1}}^{\mbox{\tiny skew}}(\mysmallP/\mysmallQ)$ will not in general be more `efficient' -- as measured by total edges --  than the associated disjoint sum $L_{\mytinyA_{n-1}}^{\mbox{\tiny GT}}(\mysmallP_{1}) \oplus \cdots \oplus L_{\mytinyA_{n-1}}^{\mbox{\tiny GT}}(\mysmallP_{k})$ of classical GT lattices. 
However, when the skew Schur function $\vartheta_{_{\mytinyP/\mytinyQ}} = \vartheta_{_{\mytinyP_{1}}}+\cdots+\vartheta_{_{\mytinyP_{k}}}$ is not multiplicity-free in the sense of \cite{GHO}, then the companion disjoint sum $L_{\mytinyA_{n-1}}^{\mbox{\tiny GT}}(\mysmallP_{1}) \oplus \cdots \oplus L_{\mytinyA_{n-1}}^{\mbox{\tiny GT}}(\mysmallP_{k})$ will not be a solitary supporting graph, although it will be edge-minimal. 
By contrast, it can be checked that the 144-element skew-tabular lattice $L := L_{\mytinyA_{2}}^{\mbox{\tiny skew}}({\setlength{\unitlength}{0.125cm}\begin{picture}(5,0) \put(0,-0.75){\line(0,1){1}} \put(0,-0.75){\line(1,0){1}} \put(0,0.25){\line(1,0){1}} \put(1,-0.75){\line(0,1){2}} \put(2,-0.75){\line(0,1){2}} \put(3,0.25){\line(0,1){2}} \put(4,0.25){\line(0,1){2}} \put(1,-0.75){\line(1,0){1}} \put(1,0.25){\line(1,0){3}} \put(1,1.25){\line(1,0){3}} \put(3,2.25){\line(1,0){1}} \put(4,1.25){\line(1,0){1}} \put(5,1.25){\line(0,1){1}} \put(4,2.25){\line(1,0){1}} \end{picture}})$ with rank generating function $\RGF(L,q) = q^{10}+5q^{9}+11q^{8}+18q^{7}+24q^{6}+26q^{5}+24q^{4}+18q^{3}+11q^{2}+5q+1$ is solitary and edge-minimal even though the associated skew Schur function $\vartheta_{_{\setlength{\unitlength}{0.075cm}\begin{picture}(5,0) \put(0,-0.75){\line(0,1){1}} \put(0,-0.75){\line(1,0){1}} \put(0,0.25){\line(1,0){1}} \put(1,-0.75){\line(0,1){2}} \put(2,-0.75){\line(0,1){2}} \put(3,0.25){\line(0,1){2}} \put(4,0.25){\line(0,1){2}} \put(1,-0.75){\line(1,0){1}} \put(1,0.25){\line(1,0){3}} \put(1,1.25){\line(1,0){3}} \put(3,2.25){\line(1,0){1}} \put(4,1.25){\line(1,0){1}} \put(5,1.25){\line(0,1){1}} \put(4,2.25){\line(1,0){1}} \end{picture}}} = \vartheta_{(5,2)}+\vartheta_{(4,0)}+2\vartheta_{(4,3)}+2\vartheta_{(3,1)}+\vartheta_{(2,2)}+\vartheta_{(1,0)}$ is not multiplicity-free. 

Such observations lead to many questions about skew-tabular lattices, two of which we bring to the reader's attention now as open problems. 

\noindent
{\bf \OpenQuestionOne}\ \ For which skew shapes is the associated skew-tabular lattice solitary and/or edge-minimal? 

\noindent 
See \S \SolitarySection\ below for some further discussion of this and related concepts. 
McNamara and van Willigenburg proposed a conjecture (Conjecture 3.26 of \cite{MvW}) specifying precisely when two skew Schur functions, say $\vartheta_{_{\mytinyP/\mytinyQ}}$ and $\vartheta_{_{\mytinyM/\mytinyN}}$, are equal. 
It is easy to deduce from first principles that if the skew-tabular lattices $L_{\mytinyA_{n-1}}^{\mbox{\tiny skew}}(\mysmallP/\mysmallQ)$ and $L_{\mytinyA_{n-1}}^{\mbox{\tiny skew}}(\mysmallM/\mysmallN)$ are isomorphic as diamond-colored distributive lattices, then $\vartheta_{_{\mytinyP/\mytinyQ}} =\vartheta_{_{\mytinyM/\mytinyN}}$. 
It is natural to ask about the nature of any converse relationship between equal skew Schur functions and their associated skew-tabular lattices: 

\noindent
{\bf \OpenQuestionTwo}\ \ When two skew Schur functions are equal, what can be said about the relationship between their associated skew-tabular lattices?

We remark that constructions of certain representations of the Lie superalgebra $\mathfrak{gl}_{m|n}$ considered in \cite{Stoilova} and \cite{MoSuperLie} and extended in \cite{Lu} and \cite{FSZ} have bases indexed by objects similar to the skew-shaped semistandard tableaux and corresponding parallelogram-shaped arrays used here.  
In those papers, however, the natural restrictions to general linear Lie algebra modules do not, in general, yield the reducible modules we study here. 
Even so, one can find in \cite{PSV} an explicit verification of relations for actions of generators of the general linear Lie algebra on the GT bases for irreducible representations.  
That said, the combinatorics of supporting graphs is at most implicit in these various constructions.  

Before we close this introduction, we mention several applications, appearing here and elsewhere, of skew-tabular lattices. 
In this paper, we present new, explicit weight basis constructions of some orthogonal Lie algebra representations associated with multiples of the `spin-node' minuscule fundamental weights in types $\myB$ and $\myD$. 
These constructions are obtained using some well-known distributive lattices (cf.\ \cite{PrEur}) related to these type $\myB/\myD$ dominant weights by recognizing certain large components of these lattices as (non-skew) classical GT lattices. 
We use a very similar approach in \cite{DDW} to explicitly construct of some infinite families of irreducible representations of the simple Lie algebras of types $\myE_{6}$ and $\myE_{7}$, but in those cases non-classical skew-tabular lattices are needed. 
In \cite{DDMN}, we use non-classical skew-tabular lattices to study a certain Fibonacci sequence of distributive lattices from \cite{HH} as well as some generalizations that we more recently discovered and named `symmetric Fibonaccian lattices'.  (\IntroFigTwo\ depicts an example.) 
In fact, it was our work with these symmetric Fibonaccian lattices that motivated the investigations undertaken in this paper. 
That skew-tabular lattices are rank symmetric, rank unimodal, and strongly Sperner is an immediate consequence of our constructions. 
While rank symmetry and rank unimodality are already known for the distributive lattices used in our orthogonal representation constructions, our proof that these lattices possess the strong Sperner property seems to be new. 
In a future paper, we will extend ideas from this paper to define, in a general and type-independent way, analogs of skew-Schur functions for all simple Lie algebras.  
Within that completely general setting, we will showcase a new, simple, and combinatorial proof of a generalization of the Littlewood--Richardson Rule whose type $\myA$ version was first obtained by Zelevinsky \cite{Z} and later augmented by Stembridge \cite{StemLR} -- see \LRCorollaryOne\ below for a statement of this type $\myA$ `ZS Rule'. 

The paper is organized as follows. 
In \S \SetupSectionOne, we provide a brief backgrounder on splitting posets for Weyl symmetric functions.  
In \S \SetupSectionTwo, we present those definitions and basic concepts about supporting graphs and representation diagrams that are relevant for our subsequent discussion. 
In \S \CombinatorialFeatures, we highlight some prominent combinatorial features of splitting posets and supporting graphs. 
The main result of that section (\MainCombinatorialTheorem) is derivative but has not appeared elsewhere in this same general form. 
In \S \GTLatticesSection, we define the skew-tabular lattices and related combinatorial objects that afford our new weight basis constructions. 
In \S \GTStatementsSection, we formally present our new weight basis constructions (\GTGeneralTheorem). 
We also note several combinatorial consequences for skew-tabular lattices, some of which appear to be new. 
In \S \OrthogonalSection, we use classical GT lattices to help us construct new weight bases for certain orthogonal Lie algebra representations (\NewOrthogonalConstructions/\NewOrthogonalCorollary). 
In \S \SolitarySection, we demonstrate certain extremal properties of our special linear and orthogonal Lie algebra representation constructions.

\noindent 
{\bf Acknowledgment}\ \ We are very grateful to the anonymous referee for offering many insights that added clarity and depth to our paper.

\vspace*{0.1in}
\noindent {\bf \S \SetupSectionOne\ Splitting posets for Weyl symmetric functions}. 
A Weyl symmetric function is a Laurent polynomial that is symmetric under a certain action of the Weyl group associated to a finite root system. 
Classical symmetric functions provide one family of examples, but to view them as Weyl symmetric functions proper requires an easy change of variables. 
The Weyl symmetric function analogs of Schur functions are what we call Weyl bialternants.   
For a self-contained account of the foundations of Weyl symmetric function theory and Weyl bialternants, see \S 2 of \cite{DonPosetModels}, a tutorial that aims to synthesize and unify standard content from various classical sources (e.g.\ \cite{FH}, \cite{Hum}). 
The purpose of the present section is to offer a distilled account of these background ideas sufficient to make accessible some results concerning weight- and rank-generating functions for the type $\myA$ skew-tabular lattices and the types $\myB$ and $\myD$ orthogonal lattices presented later. 

Here, we make note of some key notions.  
Our starting point is a finite rank $n$ root system $\Phi$ residing in an $n$-dimensional Euclidean space 
with inner product $\langle \cdot,\cdot \rangle$. 
Let $I$ be an indexing set with $|I|=n$, where we generally take $I=\{1,2,\ldots,n\}$. 
The related objects 
\rule[1mm]{3.5mm}{0.35mm} coroot system $\Phi^{\vee} := \{\alpha^{\vee} = \frac{2}{\langle \alpha,\alpha \rangle}\alpha\, |\, \alpha \in \Phi\}$; 
choice of simple roots $\{\alpha_{i}\}_{i \in I}$; 
Cartan matrix $A_{\Phi} := \left(a_{ij}\right)_{i,j \in I}$ with $a_{ij} := \langle \alpha_{i},\alpha_{j}^{\vee} \rangle$; 
sets of positive and negative roots $\Phi^{+}$ and $\Phi^{-}$ respectively; 
fundamental weights $\{\omega_{i}\}$ dual to the simple coroots $\{\alpha_{j}^{\vee}\}$ via the relations $\langle \omega_{i},\alpha_{j}^{\vee} \rangle = \delta_{ij}$; 
the lattice of weights $\Lambda := \{\mu = \sum_{i=1}^{n}m_{i}\omega_{i}\, |\, m_{i} \in \mathbb{Z}\}$; 
dominant weights $\Lambda^{+} = \{\mu = \sum_{i=1}^{n}m_{i}\omega_{i}\, |\, m_{i} \in \mathbb{Z} \mbox{ with } m_{i} \geq 0\}$; 
finite Weyl group $W=W(\Phi)$ with generators $\{s_{i}\}_{i \in I}$ and relations $(s_{i}s_{j})^{m_{ij}} = \varepsilon$ where $m_{ij} = k \in \{1,2,3,4,6\}$ with $k$ as the smallest positive integer such that $a_{ij}a_{ji} = 4\cos^{2}(\pi/k)$; 
action of $W$ on $\Lambda$ given by $s_{i}.\mu = \mu - \langle \mu,\alpha_{i}^{\vee} \rangle \alpha_{i}$ for each $i \in I$ and each $\mu \in \Lambda$; 
longest Weyl group element $w_{0}$ and involution $\sigma_{0}: I \longrightarrow I$ wherein $w_{0}.\alpha_{i} = -\alpha_{\sigma_{0}(i)}$; 
special elements $\varrho := \sum \omega_{i} = \frac{1}{2}\sum_{\alpha \in \Phi^{+}}\alpha$ and $\varrho^{\vee} = \frac{1}{2}\sum_{\alpha \in \Phi^{+}}\alpha^{\vee}$; etc \rule[1mm]{3.5mm}{0.35mm} 
are obtained as usual.  
At times for clarity we use `$\Phi$' as a modifier/superscript/subscript, as, for example, in the phrase``$\Phi$-dominant weight $\sum_{i=1}^{n}m_{i}\omega_{i}^{\Phi} \in \Lambda_{\Phi}^{+}$.'' 
Any irreducible root system of rank $n$ is of classification type $\myX_{n}$, where $\myX \in \{\myA, \myB, \myC, \myD, \myE, \myF, \myG\}$, cf.\ \cite{Hum}.

The group ring $\mathbb{Z}[\Lambda]$ has as a $\mathbb{Z}$-basis the formal exponentials $\{e^{\mu}\}_{\mu \in \Lambda}$.  
The Weyl group $W$ acts on $\mathbb{Z}[\Lambda]$ via $s_{i}.e^{\mu} := e^{s_{i}.\mu}$. 
From here on, we identify each $e^{\omega_{i}}$ as the indeterminate $z_{i}$. 
Then for $\mu = \sum_{i \in I} \mu_{i}\omega_{i} \in \Lambda$, the quantity $e^{\mu} = e^{\mu_{1}\omega_{1} + \mu_{2}\omega_{2} + \cdots + \mu_{n}\omega_{n}}$ is hereafter the monomial $z_{1}^{\mu_{1}}z_{2}^{\mu_{2}}{\cdots}z_{n}^{\mu_{n}} =: \myvarZ^{\mu}$. 
That is, each $\chi \in \mathbb{Z}[\Lambda]$ is a Laurent polynomial in the variables $z_{1}, z_{2}, \ldots, z_{n}$. 
The {\em ring of Weyl symmetric functions} $\mathbb{Z}[\Lambda]^{W}$ is the subring of $W$-invariants in $\mathbb{Z}[\Lambda]$.  
Each $\chi \in \mathbb{Z}[\Lambda]^{W}$ is a {\em Weyl symmetric function} or {\em $W$-symmetric function}. 
The {\em subgroup of $W$-alternants} $\mathbb{Z}[\Lambda]^{\mbox{\tiny alt}}$ consists of those group ring elements $\varphi$ for which $\sigma.\varphi = \det(\sigma)\varphi$ for all Weyl group elements $\sigma$. 
Define a mapping $\mathcal{A} = \mathcal{A}^{\Phi}: \mathbb{Z}[\Lambda] \longrightarrow \mathbb{Z}[\Lambda]^{\mbox{\tiny alt}}$ by the rule $\mathcal{A}(\varphi) := \sum_{\sigma \in W}\det(\sigma) \sigma.\varphi$.  
The Weyl denominator is the alternant $\mathcal{A}(\myvarZ^{\varrho})$, which  factors as 
\[\mathcal{A}(\myvarZ^{\varrho}) = \myvarZ^{\varrho}\Big(\prod_{\alpha \in \Phi^{+}}(1-\myvarZ^{-\alpha})\Big) = \prod_{\alpha \in \Phi^{+}}(\myvarZ^{\alpha/2}-\myvarZ^{-\alpha/2}) = \myvarZ^{-\varrho}\Big(\prod_{\alpha \in \Phi^{+}}(\myvarZ^{\alpha}-1)\Big).\]  
The following is a sort of fundamental theorem for Weyl symmetric functions: 
For any dominant weight $\lambda$, there exists a unique $\chi_{_{\lambda}}^{\Phi} \in \mathbb{Z}[\Lambda]$ for which $\mathcal{A}(\myvarZ^{\varrho})\, \chi_{_{\lambda}}^{\Phi} = \mathcal{A}(\myvarZ^{\lambda+\varrho})$, and in fact the $\chi_{_{\lambda}}^{\Phi}$'s comprise a $\mathbb{Z}$-basis for $\mathbb{Z}[\Lambda]^{W}$. 
We call $\chi_{_{\lambda}}^{\Phi}$ a {\em Weyl bialternant}, as it is, ostensibly, a quotient `$\mathcal{A}^{\Phi}(\myvarZ^{\lambda+\varrho})/\mathcal{A}^{\Phi}(\myvarZ^{\varrho})$' of alternants. 

There are two well-known specializations of Weyl bialternants that are important for the combinatorics of splitting posets. 
If we set $z_{i} := q^{\langle \omega_{i},\varrho^{\vee} \rangle}$, then the Weyl bialternant $\chi_{_{\lambda}}^{\Phi}$ specializes as a Laurent polynomial in the variable $q$; now if we scale this resulting quantity by the factor $q^{\langle \lambda,\varrho^{\vee} \rangle}$, the result is a polynomial in $q$, sometimes called a `Dynkin polynomial' \cite{Pan}. 
This polynomial has the product-of-quotients form of the \SpecializationTheorem\ below. 
Some further language: A polynomial $a_{0} + a_{1}q + \cdots + a_{l}q^{l}$ of degree $l$ is {\em symmetric} if $a_{i} = a_{l-i}$ for each $i \in \{0,1,\ldots,l\}$ and {\em unimodal} if there is some $u \in \{0,1,\ldots,l\}$ such that $a_{0} \leq \cdots \leq a_{u-1} \leq a_{u} \geq a_{u+1} \geq \cdots \geq a_{l}$. 

\noindent 
{\bf \SpecializationTheorem}\ \ {\sl The Dynkin polynomial obtained by the $q$-specialization of the $W(\Phi)$-Weyl bialternant $\chi_{_{\lambda}}^{\Phi}$ described in the preceding paragraph is a symmetric and unimodal polynomial of degree $2\langle \lambda,\varrho^{\vee} \rangle$ and satisfies the following identity:}
\begin{eqnarray*}
q^{\langle \lambda,\varrho^{\vee} \rangle}\left(\chi_{_{\lambda}}^{\Phi}\rule[-3.25mm]{0.2mm}{7mm}_{\, z_{i} := q^{\langle \omega_{i},\varrho^{\vee} \rangle}}\right) & = &  \mbox{\Large $\displaystyle \prod_{\alpha \in \Phi^{+}}$} \frac{1-q^{\langle \lambda+\varrho,\alpha^{\vee} \rangle}}{1-q^{\langle \varrho,\alpha^{\vee} \rangle}}
\end{eqnarray*}

\noindent 
This result is due to Dynkin \cite{Dyn}. 
For a proof from first principles, see Theorem 2.17 of \cite{DonPosetModels}. 
The Dynkin polynomial identity that concludes the preceding proposition statement is sometimes called  the `quantum dimension formula' \cite{Ram}. 
The second of our aforementioned specializations is often called `Weyl's dimension formula' and is found by taking $q=1$ in the quantum dimension formula. 
This yields, among other things, the number of terms (counting multiplicities) in the Laurent polynomial $\chi_{_{\lambda}}^{\Phi}$. 

Next, we connect $W(\Phi)$-symmetric functions with certain types of ranked and edge-colored posets. 
Before we do so, we fix our language and notation concerning such posets and stipulate the following finiteness hypothesis: {\em From here on, all posets are taken to be finite}. 
A ranked poset is a poset $R$ together with a surjective function \rule[1mm]{3.5mm}{0.35mm} i.e.\ its rank function \rule[1mm]{3.5mm}{0.35mm} $\rho: R \longrightarrow \{0,1,\ldots,l\}$ such that $\rho(\selt) + 1 = \rho(\telt)$ whenever $\telt$ covers $\selt$ in $R$. 
When $\rho$ is understood, we just use `$R$' to refer to the ranked poset; if $R$ is connected, its rank function is unique. 
The rank generating function of $R$ is $\RGF(R) = \RGF(R,q) := \sum_{\telt \in R}q^{\rho(\telt)} = \sum_{r=0}^{l}|\rho^{-1}(r)|q^{r}$. 
For consistency of notation, we set $\myCARD(R) := |R| = \RGF(R)|_{q=1}$ and $\myLENGTH(R) := l = \deg(\RGF(R))$. 
We say the ranked poset $R$ is rank symmetric and/or rank unimodal if the $q$-polynomial $\RGF(R)$ is symmetric and/or unimodal;   
is strongly Sperner if for each positive integer $k$, the largest union of  $k$ antichains of $R$ is no larger than the largest union of $k$ ranks; 
and has a symmetric chain decomposition (SCD) if we can write $R$ as a setwise disjoint union $\coprod_{i=1}^{k} \mathcal{C}_{i}$ where each $\mathcal{C}_{i} \subseteq R$ is a totally ordered subset of $R$ with minimal and maximal elements $\selt_{i}$ and $\telt_{i}$ respectively satisfying $\rho(\telt_{i}) + \rho(\selt_{i}) = l$ and $\rho(\telt_{i}) - \rho(\selt_{i}) = |\mathcal{C}_{i}|-1$. 
One can see that if $R$ has an SCD, then $R$ is rank symmetric, rank unimodal, and strongly Sperner. 

It will be convenient to identify a given ranked poset $R$ with the directed graph of its covering relations, i.e.\ its {\em order diagram}.  
Suppose the order diagram edges of $R$ are `colored' by our simple root index set $I$. 
Write $\selt \myarrow{i} \telt$ when the directed edge $\selt \rightarrow \telt$ in the order diagram of $R$  is assigned color $i$. 
Say $R$ is {\em diamond-colored} if, on any `diamond' of edges \parbox{1.4cm}{\begin{center}
\setlength{\unitlength}{0.2cm}
\begin{picture}(5,3)
\put(2.5,0){\circle*{0.5}} \put(0.5,2){\circle*{0.5}}
\put(2.5,4){\circle*{0.5}} \put(4.5,2){\circle*{0.5}}
\put(0.5,2){\line(1,1){2}} \put(2.5,0){\line(-1,1){2}}
\put(4.5,2){\line(-1,1){2}} \put(2.5,0){\line(1,1){2}}
\put(1.25,0.55){\em \small k} \put(3.2,0.7){\em \small l}
\put(1.2,2.7){\em \small i} \put(3.25,2.55){\em \small j}
\put(3,-0.75){\footnotesize $\relt$} \put(5.25,1.75){\footnotesize $\telt$}
\put(3,4){\footnotesize $\uelt$} \put(-1,1.75){\footnotesize $\selt$}
\end{picture} \end{center}} in the order diagram, necessarily $k=j$ and $l=i$.
If $\sigma: I \longrightarrow J$ is a set mapping, then the $\sigma$-{\em recoloring of} $R$, denoted $R^{\sigma}$, is the edge-colored poset obtained by assigning edge color $\sigma(i)$ to each edge of color $i$ in $R$, for all $i \in I$.   
For $\xelt \in R$ and a subset $J$ of $I$, define the $J$-{\em component} of $\xelt$, denoted $\comp_{J}(\xelt)$, to be the set of all $\yelt$ in $R$ that can be reached from $\xelt$ via some undirected path whose edges only have colors from $J$, together with all $J$-colored edges incident with these elements. 
For any $i \in I$ and $\xelt \in R$, let $l_{i}(\xelt)$ be the length of $\comp_{i}(\xelt)$, $\rho_{i}(\xelt)$ the rank of $\xelt$ within its $i$-component, and $\delta_{i}(\xelt) = l_{i}(\xelt) - \rho_{i}(\xelt)$ the depth of $\xelt$ in its $i$-component; set $\mym_{i}(\xelt) := \rho_{i}(\xelt)-\delta_{i}(\xelt)$. 
Let the {\em weight} of $\xelt$ be given by $\wt(\xelt) := \sum_{i \in I}\mym_{i}(\xelt)\omega_{i}$. 
Set $\WGF(R) = \WGF(R;z_{1}, z_{2}, \ldots, z_{n}) := \sum_{\xelt \in R}\myvarZ^{\smallwt(\xelt)}$, the {\em weight-generating function} of $R$. 
Say $R$ is $\Phi$-structured or $A_{\Phi}$-structured if, when $\selt \myarrow{i} \telt$, we have $\wt(\selt)+ \alpha_{i} = \wt(\telt)$, i.e.\ $\mym_{j}(\selt) + a_{ij} = \mym_{j}(\telt)$ for all $j \ne i$ in $I$. 
Say $R$ is a {\em splitting poset} for a Weyl symmetric function $\chi$ if $\WGF(R) = \chi$. 
If $R$ is a modular (respectively distributive) lattice, then call $R$ a {\em splitting modular }{(resp.}{\em  distributive}) {\em lattice}.  

The following lemma is a simple observation from \cite{DonPosetModels} (Lemma 3.5) that allows us, in certain circumstances, to discern $W(\Phi)$-invariance of weight generating functions. 

\noindent
{\bf \WInvariantLemma}\ \ {\sl Let $R$ be a ranked poset with edges colored by our simple root index set $I$. 
Suppose $R$ is $A_{\Phi}$-structured and that for each $i \in I$ the color $i$ components of $R$ are rank symmetric. 
Then} $\WGF(R)$ {\sl is $W(\Phi)$-symmetric.}

\vspace*{0.1in}
\noindent {\bf \S \SetupSectionTwo\ Supporting graphs and representation diagrams as splitting posets.} 
A supporting graph is an edge-colored poset that provides a kind of picture of the actions of semisimple Lie algebra generators on a weight basis for a given representation. 
A representation diagram is such a supporting graph together with two coefficients assigned to each edge so that the generator actions on a weight basis can be fully recovered. 
In this section, we summarize some basic properties of supporting graphs and representation diagrams and recount a method for constructing/confirming that a given edge-colored poset with coefficients is indeed a representation diagram. 

This paragraph follows \cite{DonSupp} and \cite{Hum}. 
Associate to the root system $\Phi$ with its given choice of simple roots the rank $n$ complex semisimple Lie algebra $\mathfrak{g} = \mathfrak{g}(\Phi)$ with Chevalley generators $\{\myqx_{i},\myqy_{i},\myqh_{i}\}_{i \in I}$ satisfying the Serre relations, cf.\ \S 18 of \cite{Hum}.  
For the remainder of this discussion  of supporting graphs, $V$ denotes a finite-dimensional (f.d.\ for short) $\mathfrak{g}$-module.  
For any $\mu \in \Lambda$, $V_{\mu} = \{v \in V\, |\, \myqh_{i}.v = \langle \mu,\alpha_{i}^{\vee} \rangle v \mbox{ for all } i \in I\}$ is the $\mu$-weight space for 
$V$.  
We have $V = \bigoplus_{\mu \in \Lambda} V_{\mu}$, and a weight basis is any basis for $V$ that respects this decomposition.  
Finite-dimensional $\mathfrak{g}$-modules are completely reducible, and the irreducible f.d.\ modules are indexed by dominant weights.  
Let $V(\lambda)$ denote an irreducible f.d.\ $\mathfrak{g}$-module corresponding to dominant weight $\lambda$ and constructed, say, using Verma modules cf.\ \cite{Hum}; then there is a `highest' weight vector $v_{\lambda}$ (unique up to scalar multiple) such that $\myqx_{i}.v_{\lambda} = 0$ for all $i \in I$.  
For such $V(\lambda)$, we have $\mychar(V(\lambda)) := \sum_{\mu \in \Lambda}(\dim V(\lambda)_{\mu})\myvarZ^{\mu} = \chi_{_{\lambda}}$, the latter equality by the famous Weyl character formula. 
For the generic $\mathfrak{g}$-module $V$, $\mychar(V) := \sum_{\mu \in \Lambda}\, (\dim V_{\mu})\myvarZ^{\mu} \in \mathbb{Z}[\Lambda]^{W}$ is a Weyl symmetric function, and $\mychar(V) = \sum_{i=1}^{k}\chi_{_{\lambda^{(i)}}}$ if and only if $V \cong V(\lambda^{(1)}) \oplus \cdots \oplus V(\lambda^{(k)})$ for some irreducible f.d.\ $\mathfrak{g}$-modules $V(\lambda^{(i)})$.   

Given any weight basis $\{v_{\relt}\}_{\relt \in R}$ for a $\mathfrak{g}$-module $V$ (with basis vectors indexed by a set $R$), we build as follows an edge-colored directed graph using elements of $R$ as vertices.  
In particular, whenever $\selt, \telt \in R$ and $i \in I$, express each of $\myqx_{i}.v_{\selt}$ and $\myqy_{i}.v_{\telt}$ in the basis $\{v_{\relt}\}_{\relt \in R}$, so $\myqx_{i}.v_{\selt} = \sum_{\uelt \in R}\myqX_{\uelt,\selt}v_{\uelt}$ and $\myqy_{i}.v_{\telt} = \sum_{\relt \in R}\myqY_{\relt,\telt}v_{\relt}$ for some scalars $\{\myqX_{\uelt,\selt}\}_{\uelt \in R}$ and $\{\myqY_{\relt,\telt}\}_{\telt \in R}$. 
Then create an $i$-colored and directed edge $\selt \myarrow{i} \telt$ if $\myqX_{\telt,\selt} \not= 0$ or $\myqY_{\selt,\telt} \not= 0$.  
For simplicity, we use $R$ to denote the resulting edge-colored digraph and call it a {\em supporting graph} for $V$. 
The supporting graph $R$ together with the set of coefficient pairs $\{(\myqX_{\telt,\selt},\myqY_{\selt,\telt})\}_{\selt \rightarrow \telt \mbox{\scriptsize \ in } R}$ assigned to the edges of $R$ is a {\em representation diagram}. 
We say the supporting graph $R$ or any associated weight basis $\{v_{\relt}\}_{\relt \in R}$ is {\em solitary} if, whenever any other weight basis $\{w_{\relt}\}_{\relt \in R}$ has supporting graph $R$, then there exist nonzero scalars $\{a_{\relt}\}_{\relt \in R}$ such that $w_{\relt} = a_{\relt}v_{\relt}$ for all $\relt \in R$, i.e.\ up to a change of basis by some diagonal matrix, the weight basis $\{v_{\relt}\}_{\relt \in R}$ is uniquely specified by $R$.  
We say $R$ is {\em edge-minimal} if no supporting graph for $V$ is isomorphic (as an edge-colored directed graph) to a subgraph of $R$. 

The following proposition gathers together in one place some simple observations (mainly from \cite{DonSupp}) concerning supporting graphs. 

\noindent 
{\bf \SupportingGraphProp}\ \ {\sl Any supporting graph for an f.d.\ $\mathfrak{g}$-module $V$ is a splitting poset for the Weyl symmetric function} $\mychar(V)$. 
{\sl Any supporting graph for $V(\lambda)$ is a connected splitting poset for $\chi_{_{\lambda}}$.}

{\em Proof.} Any supporting graph $R$ for $V$ is $A_{\Phi}$-structured by Lemmas 3.1.A and 3.2.A of \cite{DonSupp}.  
Now $\displaystyle \WGF(R) = \sum_{\relt \in 
R}\myvarZ^{\smallwt(\relt)} = \sum_{\mu \in \Lambda}\left(\sum_{\relt: \smallwt(\relt) = 
\mu}\myvarZ^{\mu}\right) = \sum_{\mu \in \Lambda}(\dim V_{\mu})\myvarZ^{\mu} = 
\mychar(V)$, which is a Weyl symmetric function by the paragraph 
preceding the proposition statement.  
So $R$ is a splitting poset for $\mychar(V)$. Now say $V = 
V(\lambda)$.  Lemma 3.1.F of \cite{DonSupp} guarantees that the 
supporting graph $R$ is connected.\hfill\QED

Any poset $R$ with edges colored by our simple root index set $I$ and with a pair of scalars $(\myqX_{\telt,\selt},\myqY_{\selt,\telt})$ assigned to each of its edges $\selt \myarrow{i} \telt$ will be called an {\em edge-tagged} poset. 
Suppose now that we are given an edge-tagged diamond-colored modular or distributive lattice $L$, and regard each member of the pair $(\myqX_{\telt,\selt},\myqY_{\selt,\telt})$ assigned to edge $\selt \myarrow{i} \telt$ of $L$ to be a scalar variable. 
We wish to identify combinatorial criteria sufficient to guarantee that $L$ is a supporting graph for some $\mathfrak{g}(\Phi)$-module. 
For any diamond of edges \parbox{1.4cm}{\begin{center}
\setlength{\unitlength}{0.2cm}
\begin{picture}(5,3)
\put(2.5,0){\circle*{0.5}} \put(0.5,2){\circle*{0.5}}
\put(2.5,4){\circle*{0.5}} \put(4.5,2){\circle*{0.5}}
\put(0.5,2){\line(1,1){2}} \put(2.5,0){\line(-1,1){2}}
\put(4.5,2){\line(-1,1){2}} \put(2.5,0){\line(1,1){2}}
\put(1.25,0.55){\em \small j} \put(3.2,0.7){\em \small i}
\put(1.2,2.7){\em \small i} \put(3.25,2.55){\em \small j}
\put(3,-0.75){\footnotesize $\relt$} \put(5.25,1.75){\footnotesize $\telt$}
\put(3,4){\footnotesize $\uelt$} \put(-1,1.75){\footnotesize $\selt$}
\end{picture} \end{center}}, the {\em diamond relations} are $\myqX_{\uelt,\selt}\myqY_{\telt,\uelt} = \myqY_{\relt,\selt}\myqX_{\telt,\relt}$ and $\myqX_{\uelt,\telt}\myqY_{\selt,\uelt} = \myqY_{\relt,\telt}\myqX_{\selt,\relt}$. 
The color $i$ {\em crossing relation} at each $\telt \in L$ is the relation 
\[\sum_{\selt: \selt \myarrow{i} \telt} \myqX_{\telt,\selt}\myqY_{\selt,\telt} - \sum_{\uelt: \telt \myarrow{i} \uelt} \myqX_{\uelt,\telt}\myqY_{\telt,\uelt} = \mym_{i}(\telt).\]
Taken all together, these relations comprise the {\em DC relations}. 
Let $V[L]$ be the complex vector space freely generated by $\{v_{\telt}\}_{\telt \in L}$. 
For each $i \in I$ and $\telt \in L$, declare that
\vspace*{0.1in}

\noindent
$\displaystyle (\blacklozenge)\hspace*{1.25in} \myqx_{i}.v_{\telt} \stackrel{\mbox{\tiny def}}{:=} \sum_{\uelt \in L}\myqX_{\uelt,\telt}v_{\uelt}\ \ \mbox{ and }\ \ \myqy_{i}.v_{\telt} \stackrel{\mbox{\tiny def}}{:=} \sum_{\relt \in L}\myqY_{\relt,\telt}v_{\relt}.$

Given an edge-tagged diamond-colored modular or distributive lattice $L$ with an edge $\selt \myarrow{i} \telt$, we let $\myqP_{\selt,\telt} := \myqX_{\telt,\selt}\myqY_{\selt,\telt}$ be the {\em edge product} and let $\myqQ_{\selt,\telt} := \sqrt{\myqP_{\selt,\telt}}$ be the principal square root of $\myqP_{\selt,\telt}$. 
Let $L' := \{\xelt'\}_{\xelt \in L}$ be the edge-colored directed graph with $\selt' \myarrow{i} \telt'$ in $L'$ if and only if edge $\selt \myarrow{i} \telt$ in $L$ has nonzero edge product. 
We make $L'$ an edge-tagged poset by assigning $(\myqX'_{\telt',\selt'} := \myqQ_{\selt,\telt},\myqY'_{\selt',\telt'} := \myqQ_{\selt,\telt})$ to each edge $\selt' \myarrow{i} \telt'$ of $L'$. 

Part {\sl (1)} of the next result is a mild modification of Proposition 3.4 from \cite{DonSupp} and Lemma 3.1 of \cite{DLP2}. 
Part {\sl (2)} follows easily from {\sl (1)}.

\noindent 
{\bf \RepDiagramTheorem}\ \ {\sl Keep the notation of the preceding paragraphs.  
Suppose each edge $\selt \rightarrow \telt$ is assigned a pair of scalars $(\myqX_{\telt,\selt},\myqY_{\selt,\telt})$ wherein at least one scalar of each pair is nonzero. 
(1) Then, the lattice $L$ is $A_{\Phi}$-structured and the scalars satisfy all DC relations if and only if (i) the action of the generators of $\mathfrak{g}(\Phi)$ defined by the formulas} ($\blacklozenge$) {\sl above is well-defined, (ii) $\{v_{\relt}\}_{\relt \in L}$ is a weight basis for the $\mathfrak{g}(\Phi)$-module $V[L]$, and (iii) the lattice $L$ together with the set of scalar pairs $\{(\myqX_{\telt,\selt},\myqY_{\selt,\telt})\}_{\selt \rightarrow \telt \mbox{\scriptsize \ in } L}$ is its representation diagram. 
(2) Moreover, if the equivalent condition of (1) holds and if each edge product in $L$ is nonzero, then $L'$ coincides with $L$ as a diamond-colored modular or distributive lattice and the biconditional statement of (1) holds when we replace $L$ with the edge-tagged lattice $L'$.}

\vspace*{0.1in}
\noindent {\bf \S \CombinatorialFeatures\ Salient combinatorial features of connected splitting posets and supporting graphs.} 
The fact that splitting posets for Weyl bialternants are rank symmetric and rank unimodal is a direct result of ideas originally due to Dynkin \cite{Dyn}, cf.\ \SpecializationTheorem\ above.  
That connected supporting graphs are rank symmetric, rank unimodal, and strongly Sperner is an application of Proctor's `Peck Poset Theorem' \cite{PrPeck} (see also Proposition 3.11 of \cite{DonSupp}). 
The content of the following proposition is borrowed from these various sources. 

\noindent
{\bf \MainCombinatorialTheorem}\ \ {\sl Suppose $R$ is a connected splitting poset for a $W(\Phi)$-symmetric function $\chi = \sum_{i=1}^{k}\chi_{_{\lambda^{(i)}}}$. 
There is a dominant weight $\lambda \in \{\lambda^{(1)},\ldots,\lambda^{(k)}\}$ such that} $\wt(\melt) = \lambda$ {\sl for some maximal element $\melt \in R$ and such that} $\langle \wt(\xelt)+\lambda,\varrho^{\vee} \rangle \leq \langle \wt(\melt)+\lambda,\varrho^{\vee} \rangle$ {\sl for all $\xelt \in R$. 
Let $\rho: R \longrightarrow \mathbb{Z}$ be the function}  $\xelt \stackrel{\rho}{\mapsto} \langle \wt(\xelt)+\lambda,\varrho^{\vee} \rangle$.  
{\sl Then $\rho(R) = \{0,1,\ldots,2\langle\lambda,\varrho^{\vee}\rangle\}$, $\rho$ is the unique rank function for $R$, and $R$ is rank symmetric and rank unimodal.
In addition to the identities} $\displaystyle \WGF(R)  =  \sum_{i=1}^{k}\chi_{_{\lambda^{(i)}}}$ {\sl and} $\myLENGTH(R) = 2\langle\lambda,\varrho^{\vee}\rangle$, {\sl we have:}
\begin{eqnarray*}
\RGF(R) & = & \mbox{\Large $\displaystyle \sum_{i = 1}^{k}$}\ \left(q^{\langle \lambda-\lambda^{(i)},\varrho^{\vee} \rangle} \mbox{\Large $\displaystyle \prod_{\alpha \in 
\Phi^{+}}$} 
\frac{1-q^{\langle \lambda^{(i)}+\varrho,\alpha^{\vee} 
\rangle}}{1-q^{\langle \varrho,\alpha^{\vee} \rangle}}\right)\\
\myCARD(R) & = & \mbox{\Large $\displaystyle \sum_{i = 1}^{k}$}\ \left(\mbox{\Large $\displaystyle \prod_{\alpha \in 
\Phi^{+}}$} 
\frac{{\langle \lambda^{(i)}+\varrho,\alpha^{\vee} 
\rangle}}{{\langle \varrho,\alpha^{\vee} \rangle}}\right)
\end{eqnarray*}
{\sl Further, suppose $R$ is a supporting graph for the f.d.\ $\mathfrak{g}(\Phi)$-module $V \cong 
V(\lambda^{(1)}) \oplus \cdots \oplus V(\lambda^{(k)})$, where each $V(\lambda^{(i)})$ is a highest-weight  $\lambda^{(i)}$ irreducible f.d.\ $\mathfrak{g}(\Phi)$-module.  
Then $R$ is strongly Sperner.}

When $\Phi = \myA_{n-1}$ and $R$ is a splitting poset for an $\myA_{n-1}$-Weyl bialternant, then the formulas for 
$\RGF(R)$, $\myCARD(R)$, and $\myLENGTH(R)$ can be concretized as in \TypeAFormulas\ below, cf.\ \cite{StanGLn}. 
Some notation: For any positive integer $m$, let $[m]_{q}$ denote the $q$-integer $q^{m-1} + \cdots + q^{1} + 1$; for a type $\myA_{n-1}$-dominant weight $\displaystyle \lambda := \sum_{k=1}^{n-1}\lambda_{k}\omega_{k}$, set $\displaystyle \lambda_{i}^{j} := \sum_{k=i}^{j}\lambda_{k}$. 

\noindent
{\bf \TypeAFormulas}\ \ {\sl If $R$ is a splitting poset for the} $\myA_{n-1}${\sl -Weyl bialternant $\chi_{_{\lambda}}$, then:}\\ 
{\small $\displaystyle \RGF(R) = \prod_{i=1}^{n-1}\prod_{j=i}^{n-1}\frac{[\lambda_{i}^{j}+j+1-i]_{q}}{[j+1-i]_{q}}$}, {\small $\displaystyle \myCARD(R) = \prod_{i=1}^{n-1}\prod_{j=i}^{n-1}\frac{\lambda_{i}^{j}+j+1-i}{j+1-i}$}, {\sl and} {\small $\displaystyle \myLENGTH(R) = \sum_{i=1}^{n-1}\sum_{j=i}^{n-1}\, \lambda_{i}^{j}$.}

\vspace*{0.1in}
\noindent {\bf \S \GTLatticesSection\ Skew-tabular lattices.}
The purpose of this section is to set the combinatorial environment for our $\mathfrak{g}(\myA_{n-1})$ representation constructions and to fix our notation. 
For \S\S \GTLatticesSection--\GTStatementsSection,  fix an integer $n \geq 2$ and assume that $\Phi = \myA_{n-1}$ is the reference root system.    
Also, fix an integer $m \geq n$ and two weakly decreasing $m$-tuples of nonnegative integers $\mysmallP = (\myscriptsizeP_{1},\ldots,\myscriptsizeP_{m})$ and $\mysmallQ = (\myscriptsizeQ_{1},\ldots,\myscriptsizeQ_{m})$ wherein $\myscriptsizeQ_{i} \leq \myscriptsizeP_{i}$ for $i \in \{1,2,\ldots,m\}$. 
We think of $\mysmallP$ and $\mysmallQ$ as integer partitions.  
The shape identified with the partition $\mysmallP$ is its {\em partition diagram} (aka Ferrers diagram), which is a left-justified collection of rows of empty boxes, with $\myscriptsizeP_{1}$ boxes on the top row, $\myscriptsizeP_{2}$ boxes on the 2nd row, etc, and similarly for $\mysmallQ$. 
Henceforth, we refer to partition diagram boxes as {\em cells}. 
The positions of partition diagram cells are indexed by (row,column) pairs of integers, as with matrices. 
The {\em skew shape} $\mysmallP/\mysmallQ$ is the shape obtained by removing the shape $\mysmallQ$ from the upper left portion of the shape $\mysmallP$. 
The position indexing of the cells of $\mysmallP/\mysmallQ$ is inherited from $\mysmallP$. 
From here on, we require that $\mysmallP$ and $\mysmallQ$ are {\em skew-compatible} with respect to $n$ in the sense that $n \geq \max\{k-j\, |\, \myscriptsizeP_{j} = \myscriptsizeP_{j+1} = \cdots = \myscriptsizeP_{k}\ \mbox{and}\ \myscriptsizeQ_{i}<\myscriptsizeP_{i}\ \mbox{for}\ j \leq i \leq k\}$. 
This latter condition is compatible with our use of semistandard tableaux of skew shape $\mysmallP/\mysmallQ$ whose entries are from $\{1,2,\ldots,n\}$, as no column of $\mysmallP/\mysmallQ$ has more than $n$ cells. 

\begin{figure}[t] 
\begin{center}
\GTParallelogramFigure: Positions in a GT $n$-parallelogram.

\vspace*{0.25in}
\begin{tabular}{cccccccccc}
$g_{0,-(m-1)}$ & & & & & & & & & \\
 & $g_{1,1-(m-1)}$ & & & & & & & & \\
 $g_{0,-(m-2)}$ & & $g_{2,2-(m-1)}$ & & & & & & & \\
 & $g_{1,1-(m-2)}$ & & & & & & & & \\
  $g_{0,-(m-3)}$ & & & & & & & & $g_{n-1,n-1-(m-1)}$ & \\
 & & {\tiny $\bullet$} & & & & & & & $g_{n,n-(m-1)}$\\
 & {\tiny $\bullet$} & & & &{\tiny $\bullet$} & & & $g_{n-1,n-1-(m-2)}$ & \\
 {\tiny $\bullet$} & & {\tiny $\bullet$} & & & & & & & $g_{n,n-(m-2)}$\\
 & {\tiny $\bullet$} & & & & {\tiny $\bullet$} & & & {\tiny $\bullet$} & \\
 {\tiny $\bullet$} & & {\tiny $\bullet$} & & & & & & & {\tiny $\bullet$}\\
 & {\tiny $\bullet$} & & & &{\tiny $\bullet$} & & & {\tiny $\bullet$} & \\
 {\tiny $\bullet$} & & & & & & & & & {\tiny $\bullet$} \\
 & & $g_{2,0}$ & & & {\tiny $\bullet$} & & & {\tiny $\bullet$} & \\
 & $g_{1,0}$ & & & & & & & & {\tiny $\bullet$} \\
 $g_{0,0}$ & & $g_{2,1}$ & & & {\tiny $\bullet$} & & & {\tiny $\bullet$} & \\
 & $g_{1,1}$ & & & & & & & & {\tiny $\bullet$}\\
 & & $g_{2,2}$ & & & & & & $g_{n-1,n-2}$ & \\
 & & & & & & & &  & $g_{n,n-1}$ \\
 & & & & & & & & $g_{n-1,n-1}$ & \\
 & & & & & & & & & $g_{n,n}$
\end{tabular}
\end{center}

\vspace*{-0.2in}
\end{figure} 

A {\em semistandard $n$-tableau $T$ of shape} $\mysmallP/\mysmallQ$ is a filling of the cells of shape $\mysmallP/\mysmallQ$ with entries from the set $\{1,2,\ldots,n\}$ wherein the entries weakly increase left-to-right across each row and strictly increase top-to-bottom in each column. 
We use $L_{\mytinyA_{n-1}}^{\mbox{\tiny skew}}(\mysmallP/\mysmallQ)$ to denote the collection of semistandard $n$-tableaux of shape $\mysmallP/\mysmallQ$. 
For any $i \in \{1,2,\ldots,n\}$, let $\#_{i}(T)$ denote the number of appearances of the number $i$ within the cells of $T$. 
As we are working with $\mathfrak{g}(A_{n-1}) \cong \mathfrak{sl}(n,\mathbb{C})$ rather than $\mathfrak{gl}(n,\mathbb{C})$, the following weighting of tableaux is efficacious.  
We let $\wt(T) \in \Lambda$ be given by\footnote{Observe that the definition of $\wt(T)$ also applies when $T$ is any multisubset of $\{1,2,\ldots,n\}$.} 
\[\wt(T) {:=} \sum_{i=1}^{n-1}\left(\rule[-1.75mm]{0mm}{4mm}\#_{i}(T) - \#_{i+1}(T)\right)\omega_{i}.\] 
We define the skew Schur function $\vartheta_{_{\mytinyP/\mytinyQ}}$ to be the following Laurent polynomial: 
\[\vartheta_{_{\mytinyP/\mytinyQ}} {:=} \sum_{T \in L_{\mytinyA_{n-1}}^{\mbox{\tiny skew}}(\mytinyP/\mytinyQ)}\myvarZ^{\smallwt(T)}.\] 
Our $\vartheta_{_{\mytinyP/\mytinyQ}}$ differs mildly from the classical skew Schur function $\mys_{_{\mytinyP/\mytinyQ}}$, which is defined to be the following polynomial in $\mathbb{Z}[x_{1},x_{2},\ldots,x_{n}]$:  
\[\mys_{_{\mytinyP/\mytinyQ}} := \sum_{T \in L_{\mytinyA_{n-1}}^{\mbox{\tiny skew}}(\mytinyP/\mytinyQ)}x_{1}^{\#_1(T)}x_{2}^{\#_2(T)}\cdots x_{n}^{\#_n(T)}.\]
Observe that $\vartheta_{_{\mytinyP/\mytinyQ}}$ is the image of $\mys_{_{\mytinyP/\mytinyQ}}$ under the change of variables $x_{i} = z_{i-1}^{-1}z_{i}$ taking $z_{0} := 1 =: z_{n}$. 

\begin{figure}[t] 
\begin{center}
\GTTriangleFigure: Non-inert positions in a GT $n$-triangle.

\vspace*{0.25in}
\begin{tabular}{cccccccccc}
 & & & & & & & & & $g_{n,1}$\\
 & & & & & & & & $g_{n-1,1}$ & \\
 & & & & & & & & {\tiny $\bullet$} & $g_{n,2}$\\
 & & & & & {\tiny $\bullet$} & & & & {\tiny $\bullet$} \\
 & & $g_{3,1}$ & & & & & & {\tiny $\bullet$} & \\
 & $g_{2,1}$ & & & & {\tiny $\bullet$} & & & & {\tiny $\bullet$} \\
 $g_{1,1}$ & & $g_{3,2}$ & & & & & & {\tiny $\bullet$} & \\
 & $g_{2,2}$ & & & & {\tiny $\bullet$} & & & & {\tiny $\bullet$}\\
 & & $g_{3,3}$ & & & & & & {\tiny $\bullet$} & \\
 & & & & & {\tiny $\bullet$} & & & & {\tiny $\bullet$}\\
 & & & & & & & & {\tiny $\bullet$} & $g_{n,n-1}$ \\
 & & & & & & & & $g_{n-1,n-1}$ & \\
 & & & & & & & & & $g_{n,n}$
\end{tabular}
\end{center}
\end{figure} 

Now, there is a natural partial ordering of the elements of $L_{\mytinyA_{n-1}}^{\mbox{\tiny skew}}(\mysmallP/\mysmallQ)$ that, we shall see, is effective for constructing an explicit weight basis for the corresponding representation of  $\mathfrak{g}(A_{n-1})$. 
For $S$ and $T$ in $L_{\mytinyA_{n-1}}^{\mbox{\tiny skew}}(\mysmallP/\mysmallQ)$, say $S \leq T$ if $S_{r,c} \geq T_{r,c}$ for each position $(r,c)$ of $\mysmallP/\mysmallQ$. 
With respect to this `reverse componentwise' partial ordering, one can easily see that the tableau $T$ covers $S$ if there is some position $(r,c)$ with $S_{r,c} = T_{r,c}+1$ but $S_{q,b} = T_{q,b}$ for each position $(q,b) \ne (r,c)$, in which case we give edge $S \rightarrow T$ the color $i := T_{r,c}$ and write $S \myarrow{i} T$. 

To further analyze the structure of $L_{\mytinyA_{n-1}}^{\mbox{\tiny skew}}(\mysmallP/\mysmallQ)$ as an edge-colored poset, it will help to have in hand a certain larger object to be denoted $\widetilde{L_{\mytinyA_{n-1}}^{\mbox{\tiny skew}}}(\mysmallP/\mysmallQ)$. 
The cells in a non-empty column $c$ of the skew shape $\mysmallP/\mysmallQ$ form a single-column shape of some length not exceeding $n$. 
If there are $N$ such nonempty columns, number them consecutively from left to right. 
For $1 \leq c \leq N$, denote by $\mysmallS^{(c)}$ the partition associated with the cells of (nonempty) column $c$ of our skew shape. 
Let $L_{c}$ be the diamond-colored Gelfand-Tsetlin distributive lattice denoted $L_{\mytinyA}^{\mbox{\tiny GT-left}}(\mysmallS^{(c)})$ in \cite{DonSupp} or $GT(\mysmallS^{(c)})$ in \cite{HL}. 
Elements of $L_{c}$ are columnar tableaux with shape $\mysmallS^{(c)}$ and with strictly-increasing entries (reading from top to bottom) from the set $\{1,2,\ldots,n\}$; the partial ordering and edge-coloring rules for $L_{c}$ exactly coincide with our partially ordering and edge-coloring rules for skew-shaped tableaux. 
Therefore, in the language of Section 3 of \cite{DonDiamond}, $L_{\mytinyA_{n-1}}^{\mbox{\tiny skew}}(\mysmallP/\mysmallQ)$ can be viewed as a vertex subset of the diamond-colored distributive lattice $\widetilde{L_{\mytinyA_{n-1}}^{\mbox{\tiny skew}}}(\mysmallP/\mysmallQ) := L_{1} \times \cdots \times L_{c} \times \cdots \times L_{N}$. 
It is not hard to see that $L_{\mytinyA_{n-1}}^{\mbox{\tiny skew}}(\mysmallP/\mysmallQ)$ meets the criteria of Proposition 3.5.2 of \cite{DonDiamond} (i.e.\ there is a path from the min of $\widetilde{L_{\mytinyA_{n-1}}^{\mbox{\tiny skew}}}(\mysmallP/\mysmallQ)$ to its max that only uses vertices of $L_{\mytinyA_{n-1}}^{\mbox{\tiny skew}}(\mysmallP/\mysmallQ)$, and $L_{\mytinyA_{n-1}}^{\mbox{\tiny skew}}(\mysmallP/\mysmallQ)$ is `closed under componentwise joins and meets'), from which it follows that $L_{\mytinyA_{n-1}}^{\mbox{\tiny skew}}(\mysmallP/\mysmallQ)$ is a diamond-colored distributive lattice with edges colored by $\{1,2,\ldots,n-1\}$. 
We formalize the preceding claim in \FirstGTLatticeResult\ below. 
We call $L_{\mytinyA_{n-1}}^{\mbox{\tiny skew}}(\mysmallP/\mysmallQ)$ a {\em skew-tabular lattice}.

Fix a color $i \in \{1,2,\ldots,n-1\}$ and a tableau $T$ from our edge-colored skew-tabular (distributive) lattice $L_{\mytinyA_{n-1}}^{\mbox{\tiny skew}}(\mysmallP/\mysmallQ)$. 
The definitions from \S \SetupSectionOne\ concerning $i$-components apply to $L_{\mytinyA_{n-1}}^{\mbox{\tiny skew}}(\mysmallP/\mysmallQ)$, so we can consider the quantities $\mym_{i}(T)$, $\rho_{i}(T)$, $\delta_{i}(T)$, etc.  
It is an easy exercise to check that $\mym_{i}(T) = \rho_{i}(T) - \delta_{i}(T) = \#_{i}(T)-\#_{i+1}(T)$.  
Thus, $\WGF(L_{\mytinyA_{n-1}}^{\mbox{\tiny skew}}(\mysmallP/\mysmallQ)) = \vartheta_{_{\mytinyP/\mytinyQ}}$. 

Consider an array $\telt$ of integers, to be viewed as a sequence of column vectors arranged in a staggered, parallelogram-like pattern. 
These parallelogram patterns generalize certain triangular arrays utilized in the classical constructions of \cite{GT}. 
Positions of the array are indexed by pairs $(i,j)$ with $i \in \{0,1,\ldots,n-1,n\}$ and $j \in \{i-(m-1), i-(m-2), \ldots, i-1, i\} =: C_{i}$, where the latter is to be thought of as the indexing set for the $i^{\mbox{\tiny th}}$ `column' of $\telt$; 
the entry in the $(i,j)$ position of array $\telt$ is denoted $g_{i,j}(\telt)$. 
See \GTParallelogramFigure\ for a depiction. 
We say $\telt$ is a GT $n$-{\em parallelogram framed by} $\mysmallP/\mysmallQ$ if (1) $g_{0,j}(\telt) = \myscriptsizeQ_{1-j}$ for each $j \in C_{0}$ and $g_{n,j}(\telt) = \myscriptsizeP_{n+1-j}$ for each $j \in C_{n}$ and (2) for any array entry $g_{i,j}(\telt)$ we have 
\begin{center}
{\small
\begin{tabular}{ccccc}
 $g_{i-1,j-1}(\telt)$ & & & & $g_{i+1,j}(\telt)$\\
 & \begin{rotate}{315}$\!\!\!\!\leq$\end{rotate} & & \begin{rotate}{45}$\geq$\end{rotate} & \\
 & & $g_{i,j}(\telt)$ & & \\
 & \begin{rotate}{45}$\geq$\end{rotate} & & \begin{rotate}{315}$\!\!\!\!\leq$\end{rotate} & \\
 $g_{i-1,j}(\telt)$ & & & & $g_{i+1,j+1}(\telt)$
\end{tabular}
}
\end{center}
where an inequality is ignored if $(p,q)$ is not a position within the array. 
Formally, the $i^{\mbox{\tiny th}}$ {\em column} of $\telt$ is the $m$-tuple/partition $\telt^{(i)} := (g_{i,i}(\telt),g_{i,i-1}(\telt),\ldots,g_{i,i-(m-1)}(\telt))$. 

There is a well-known one-to-one correspondence between semistandard $n$-tableaux of shape $\mysmallP/\mysmallQ$ and GT $n$-parallelograms framed by $\mysmallP/\mysmallQ$, see for example \cite{Alexandersson}. 
From a GT $n$-parallelogram $\telt$ framed by $\mysmallP/\mysmallQ$, we form a semistandard $n$-tableau $T$ of shape $\mysmallP/\mysmallQ$ by placing an entry $i$ into each cell of the partition $\telt^{(i)}/\telt^{(i-1)}$, as $i$ ranges from $1$ to $n$; reverse this process to obtain a GT $n$-parallelogram framed by $\mysmallP/\mysmallQ$ from a given semistandard $n$-tableau of shape $\mysmallP/\mysmallQ$. 
We view elements of $L_{\mytinyA_{n-1}}^{\mbox{\tiny skew}}(\mysmallP/\mysmallQ)$ simultaneously as semistandard $n$-tableaux and as GT $n$-parallelograms (\ParallelogramFig). 
We use italicized capital letters $S$, $T$, $U$, etc when referring to skew-tabular lattice elements as tableaux and bold-faced lower-case letters $\selt$, $\telt$, $\uelt$, etc when referring to such elements as parallelograms. 
In the remainder of the paper, we generally prefer GT parallelograms over tableaux, as the former better facilitate our constructions here and afford certain applications to be considered in \S \OrthogonalSection\ below and also in \cite{DD}. 

\newcommand{\VertexParallelogram}[8]{\parbox[c]{1.2cm}{\setlength{\unitlength}{0.3cm}
\begin{picture}(3,5)
\put(0,2.2){\circle*{0.5}} 
\put(#7,#8){
\begin{picture}(0,0)
\put(0,-2){\begin{picture}(0,0)
\put(0.25,3.5){\tiny 0} \put(1.25,3){\tiny #1} \put(2.25,2.5){\tiny #4} \put(3.25,2){\tiny 0} 
\put(0.25,2.5){\tiny 0} \put(1.25,2){\tiny #2} \put(2.25,1.5){\tiny #5} \put(3.25,1){\tiny 3} 
\put(0.25,1.5){\tiny 2} \put(1.25,1){\tiny #3} \put(2.25,0.5){\tiny #6} \put(3.25,0){\tiny 3} 
\end{picture}}
\end{picture}} \end{picture}}
}

\begin{figure}[th]
\begin{center}
\ParallelogramFig: The skew-tabular lattice $L_{\mytinyA_{2}}^{\mbox{\tiny skew}}({\setlength{\unitlength}{0.125cm}\begin{picture}(3,0)\put(0,0){\line(0,1){1}} \put(1,0){\line(0,1){1}} \put(2,0){\line(0,1){2}} \put(3,0){\line(0,1){2}} \put(0,0){\line(1,0){3}} \put(0,1){\line(1,0){3}} \put(2,2){\line(1,0){1}}\end{picture}})$ of \IntroFig, re-built using GT $3$-parallelograms.  

\setlength{\unitlength}{1.5cm}
\begin{picture}(4,6.5)
\put(1,0){\line(-1,1){1}}
\put(1,0){\line(1,1){2}}
\put(0,1){\line(1,1){2}}
\put(2,1){\line(-1,1){1}}
\put(2,1){\line(0,1){1}}
\put(1,2){\line(0,1){1}}
\put(2,2){\line(-1,1){2}}
\put(2,2){\line(1,1){2}}
\put(3,2){\line(-1,1){1}}
\put(3,2){\line(0,1){1}}
\put(1,3){\line(1,1){2}}
\put(2,3){\line(0,1){1}}
\put(3,3){\line(-1,1){2}}
\put(0,4){\line(1,1){2}}
\put(4,4){\line(-1,1){2}}
\put(2,6){\VertexParallelogram{0}{2}{3}{0}{3}{3}{0.0}{2.5}}
\put(1,5){\VertexParallelogram{0}{1}{3}{0}{3}{3}{-4.75}{4}}
\put(3,5){\VertexParallelogram{0}{2}{3}{0}{2}{3}{0.0}{2.5}}
\put(0,4){\VertexParallelogram{0}{0}{3}{0}{3}{3}{-4.75}{4}}
\put(2,4){\VertexParallelogram{0}{1}{3}{0}{2}{3}{-6.75}{3}}
\put(1.55,4){\vector(1,0){0.35}}
\put(4,4){\VertexParallelogram{0}{2}{2}{0}{2}{3}{0.0}{2.5}}
\put(1,3){\VertexParallelogram{0}{0}{3}{0}{2}{3}{-6.75}{3}}
\put(0.55,3){\vector(1,0){0.35}}
\put(2,3){\VertexParallelogram{0}{1}{3}{0}{2}{3}{2}{6.5}}
\put(2.5,3.75){\vector(-2,-3){0.45}}
\put(3,3){\VertexParallelogram{0}{1}{2}{0}{2}{3}{0.25}{0.75}}
\put(1,2){\VertexParallelogram{0}{0}{3}{0}{1}{3}{-6.75}{3}}
\put(0.55,2){\vector(1,0){0.35}}
\put(2,2){\VertexParallelogram{0}{0}{2}{0}{2}{3}{-6.5}{-2}}
\put(1.55,1.1){\vector(1,2){0.375}}
\put(3,2){\VertexParallelogram{0}{1}{2}{0}{1}{3}{0.25}{0.75}}
\put(0,1){\VertexParallelogram{0}{0}{3}{0}{0}{3}{-4.75}{4}}
\put(2,1){\VertexParallelogram{0}{0}{2}{0}{1}{3}{0.25}{0.75}}
\put(1,0){\VertexParallelogram{0}{0}{2}{0}{0}{3}{2.5}{1.75}}
\put(1.55,0){\vector(-1,0){0.35}}
\put(1,5){\NEEdgeLabelForLatticeI{{\em 1}}}
\put(3,5){\NWEdgeLabelForLatticeI{{\em 2}}}
\put(0,4){\NEEdgeLabelForLatticeI{{\em 1}}}
\put(2,4){\NWEdgeLabelForLatticeI{{\em 2}}}
\put(2,4){\NEEdgeLabelForLatticeI{{\em 1}}}
\put(4,4){\NWEdgeLabelForLatticeI{{\em 1}}}
\put(1,3){\NEEdgeLabelForLatticeI{{\em 1}}}
\put(1,3){\NWEdgeLabelForLatticeI{{\em 2}}}
\put(2,3){\VerticalEdgeLabelForLatticeI{{\em 2}}}
\put(3,3){\NWEdgeLabelForLatticeI{{\em 1}}}
\put(3,3){\NEEdgeLabelForLatticeI{{\em 1}}}
\put(1,2){\VerticalEdgeLabelForLatticeI{{\em 2}}}
\put(1.25,2.25){\NEEdgeLabelForLatticeI{{\em 1}}}
\put(2.2,1.8){\NWEdgeLabelForLatticeI{{\em 1}}}
\put(1.8,1.8){\NEEdgeLabelForLatticeI{{\em 1}}}
\put(3,2){\VerticalEdgeLabelForLatticeI{{\em 2}}}
\put(2.75,2.25){\NWEdgeLabelForLatticeI{{\em 1}}}
\put(0,1){\NEEdgeLabelForLatticeI{{\em 2}}}
\put(2,1){\VerticalEdgeLabelForLatticeI{{\em 2}}}
\put(2,1){\NWEdgeLabelForLatticeI{{\em 1}}}
\put(2,1){\NEEdgeLabelForLatticeI{{\em 1}}}
\put(1,0){\NWEdgeLabelForLatticeI{{\em 1}}}
\put(1,0){\NEEdgeLabelForLatticeI{{\em 2}}}
\end{picture}
\end{center}
\end{figure}

One can see that $S \myarrow {i} T$ for tableaux $S$ and $T$ in $L_{\mytinyA_{n-1}}^{\mbox{\tiny skew}}(\mysmallP/\mysmallQ)$ if and only if the corresponding parallelograms $\selt$ and $\telt$ are identical except in some position $(i,j)$ wherein $g_{i,j}(\selt)+1=g_{i,j}(\telt)$. 
It is easily checked, then, that within the $i$-component $\comp_{i}(\telt)$ we have
\begin{eqnarray*}
\mym_{i}(\telt) & = & 2\sum_{q \in C_{i}}g_{i,q}(\telt) - \sum_{q \in C_{i+1}}g_{i+1,q}(\telt) - \sum_{q \in C_{i-1}}g_{i-1,q}(\telt)\\
 & = & \sum_{q=0}^{m-1}\left(\rule[-2.5mm]{0mm}{6.5mm}\left[\rule[-1.25mm]{0mm}{4.5mm}g_{i,i-q}(\telt)-g_{i+1,i+1-q}(\telt)\right] + \left[\rule[-1.25mm]{0mm}{4.5mm}(g_{i,i-q}(\telt)-g_{i-1,i-1-q}(\telt)\right]\right).
\end{eqnarray*}
We describe the unique maximal element $\melt$ and unique minimal element $\nelt$ of  $L_{\mytinyA_{n-1}}^{\mbox{\tiny skew}}(\mysmallP/\mysmallQ)$ as follows. 
Momentarily set $\myscriptsizeP_{r} := \myscriptsizeQ_{m}$ if $r > m$ and $\myscriptsizeQ_{r} := \myscriptsizeP_{1}$ if $r < 1$. 
Then, $g_{i,j}(\melt) = \min\{\myscriptsizeQ_{1-j},\myscriptsizeP_{1+i-j}\}$ and $g_{i,j}(\nelt) = \max\{\myscriptsizeQ_{1+i-j},\myscriptsizeP_{n+1-j}\}$.  
For any GT $n$-parallelogram $\xelt$ framed by $\mysmallP/\mysmallQ$, we have $\rho(\xelt) = \sum \left(\rule[-2.5mm]{0mm}{6.5mm}g_{i,j}(\xelt)-g_{i,j}(\nelt)\right)$ and $\delta(\xelt) = \sum \left(\rule[-2.5mm]{0mm}{6.5mm}g_{i,j}(\melt)-g_{i,j}(\xelt)\right)$, where each sum is taken over all positions $(i,j)$ within the defining parallelogram.  

When $\mysmallQ = (0,\ldots,0)$ so that the shape $\mysmallP/\mysmallQ$ is non-skew, then the lattice $L_{\mytinyA_{n-1}}^{\mbox{\tiny GT}}(\mysmallP) := L_{\mytinyA_{n-1}}^{\mbox{\tiny skew}}(\mysmallP/\mysmallQ)$ is what we call a {\em classical Gelfand--Tsetlin lattice}, cf.\ \cite{PrGZ}, \cite{DonSupp}, \cite{HL}.
For any $\telt \in L_{\mytinyA_{n-1}}^{\mbox{\tiny GT}}(\mysmallP)$, we have $g_{i,j}(\telt) = 0$ if $j \leq 0$. 
That is, all of the information for any such array is contained within the triangle of {\em non-inert positions} given by $\{(i,j)\, |\, 1 \leq j \leq i \leq n\}$. 
Therefore we can refer to the arrays comprising $L_{\mytinyA_{n-1}}^{\mbox{\tiny GT}}(\mysmallP)$ as {\em GT $n$-triangles framed by} $\mysmallP$. 
For a depiction of these positions, see \GTTriangleFigure. 
Versions of GT $n$-triangles framed by non-skew shapes were the objects used in the weight basis constructions of the irreducible $\mathfrak{g}(\myA_{n-1})$-modules in \cite{GT}, so such objects are often called Gelfand--Tsetlin patterns. 
Those GT patterns are the same as our GT triangles but are reflected across the line $y=x$; we prefer our visual orientation of these arrays in part because it seems to  afford a more direct connection with posets of join irreducibles of the skew-tabular lattices, an idea that is further explored in \cite{DD}.  

Next we construct an isomorphism between a given single-color component of $L_{\mytinyA_{n-1}}^{\mbox{\tiny skew}}(\mysmallP/\mysmallQ)$ and a product of chains. 
This furnishes a demonstration of one part of \FirstGTLatticeResult, which is below. 
For any $\xelt \in L_{\mytinyA_{n-1}}^{\mbox{\tiny skew}}(\mysmallP/\mysmallQ)$ and color $k \in \{1,2,\ldots,n-1\}$, the $k$-component $\comp_{k}(\xelt)$ is naturally isomorphic to a product of chains in the following way. 
Note that \[\max\{g_{k-1,k-1-r}(\xelt),g_{k+1,k-r}(\xelt)\} \leq \min\{g_{k-1,k-r}(\xelt),g_{k+1,k+1-r}(\xelt)\},\] where $g_{p,q}(\xelt)$ is to be ignored in the preceding comparisons if $(p,q)$ is not a valid position. 
So for $i \in \{0,1,\ldots,m-1\}$ let $b_{i} = b_{i}(\xelt) := \max\{g_{k-1,k-1-i}(\xelt),g_{k+1,k-i}(\xelt)\}$ and $t_{i} = t_{i}(\xelt) := \min\{g_{k-1,k-i}(\xelt),g_{k+1,k+1-i}(\xelt)\}$. 
In this case, $\yelt \in \comp_{k}(\xelt)$ if and only if $g_{p,q}(\yelt) = g_{p,q}(\xelt)$ for all positions with $p \ne k$ and $b_{i}(\xelt) \leq g_{k,k-q_{i}(\xelt)}(\yelt) \leq t_{i}(\xelt)$ for all $i  \in \{0,1,\ldots,m-1\}$. 
That is, $\comp_{k}(\xelt)$ is isomorphic to the chain product $\mathcal{C}^{(k)}(\xelt) := \{(y_{0},\ldots,y_{m-1}) \in \mathbb{Z}^{m}\, |\, b_{i} \leq y_{i} \leq t_{i} \mbox{ for } 0 \leq i \leq m-1\}$ via the mapping $\phi: \comp_{k}(\xelt) \longrightarrow \mathcal{C}^{(k)}(\xelt)$ given by $\phi(\yelt) = (g_{k,k}(\yelt),\ldots,g_{k,k-(m-1)}(\yelt))$.

\noindent
{\bf \FirstGTLatticeResult}\ \ {\sl The skew-tabular lattice} $L:=L_{\mytinyA_{n-1}}^{\mbox{\tiny skew}}(\mysmallP/\mysmallQ)$ {\sl is an} $\myA_{n-1}$-{\sl structured diamond-colored distributive lattice, and, for each $\xelt \in L$ and color $k \in \{1,2,\ldots,n-1\}$, the $k$-component} $\comp_{k}(\xelt)$ {\sl is isomorphic to a product of chains. 
In particular,} $\WGF(L) = \vartheta_{_{\mytinyP/\mytinyQ}}$ {\sl is an} $\myA_{n-1}$-{\sl symmetric function. 
Moreover, $\rho(\xelt) = \sum_{i=1}^{n-1}\big(\mym_{i}(\xelt)+\mym_{i}(\melt)\big)$, where $\melt$ is the unique maximal element of $L$.}

{\em Proof:} That $K := \widetilde{L_{\mytinyA_{n-1}}^{\mbox{\tiny skew}}}(\mysmallP/\mysmallQ)$ is a diamond-colored distributive lattice follows from Proposition 3.5.1 of \cite{DonDiamond}. 
To see that $L$ is a diamond-colored distributive lattice, it suffices to show that, as a vertex subset of $K$, $L$ satisfies the hypotheses of Proposition 3.5.2 of \cite{DonDiamond}. 
Of course, the maximal tableau of $L$ is also the maximal tableau of $K$, and each is obtained by placing the smallest possibly entry within each cell of each column of the skew shape $\mysmallP/\mysmallQ$; the minimal tableaux of $L$ and $K$ also coincide, and in this case each cell of each column has the largest possible entry. 
By minimizing the entries of the minimal tableau one cell at a time, reading cells from top to bottom within each column and reading columns from left-to-right, it is clear that we do not at any point violate any of the defining inequalities for semistandard-ness of $\mysmallP/\mysmallQ$-shaped tableaux. 
In this way, we obtain a path from the min of $K$ to its max that only uses vertices of $L$. 
It is an easy exercise to verify that when $S$ and $T$ are tableaux in $L$, then their component-wise join, i.e.\ the tableau in $K$ whose entry in its $(r,c)$ cell is $S_{r,c} \vee T_{r,c} = \min\{S_{r,c},T_{r,c}\}$, satisfies the  defining inequalities for semistandard-ness of $\mysmallP/\mysmallQ$-shaped tableaux and is therefore in $L$ as well. 
A similar statement holds for the component-wise meet of $S$ and $T$. 
Proposition 3.5.2 of \cite{DonDiamond} now applies. 

It is shown in the paragraph preceding the proposition statement that each $k$-component is a product of chains.  
So, to apply \WInvariantLemma, it suffices to show that $L$ is $\myA_{n-1}$-structured. 
To do so, we check that when $\selt \myarrow{i} \telt$ in $L$, then $\mym_{j}(\selt) + \langle \alpha_{i},\alpha_{j}^{\vee} \rangle = \mym_{j}(\telt)$ for all $j \not= i$. 
Let us say that $\selt \myarrow{i} \telt$ with $g_{i,i-q}(\telt) = g_{i,i-q}(\selt)+1$.  
First we analyze the case that $1 \leq i < n-2$ with $j=i+1$.  
In particular, $\langle \alpha_{i},\alpha_{j}^{\vee} \rangle = -1$.  
Then $\mym_{j}(\telt) - \mym_{j}(\selt) = -g_{j+1,j+1-q}(\telt)+g_{j+1,j+1-q}(\selt) = -1$.
Analysis of the case $1 < i \leq n-1$ and $j=i-1$ is entirely similar. 
When $i$ and $j$ satisfy $|i-j| \geq 2$ so that $\langle \alpha_{i},\alpha_{j}^{\vee} \rangle = 0$, then the terms of the sum for $\mym_{j}(\telt)$ are exactly the same as the those for $\mym_{j}(\selt)$, hence $\mym_{j}(\telt) = \mym_{j}(\selt)$.
It follows from \MainCombinatorialTheorem\ that the formula for $\rho(\xelt)$ is as claimed.\hfill\QED

Continue with $L:=L_{\mytinyA_{n-1}}^{\mbox{\tiny skew}}(\mysmallP/\mysmallQ)$. 
By expressing the $\myA_{n-1}$-symmetric function $\WGF(L) = \vartheta_{_{\mytinyP/\mytinyQ}}$ in terms of classical Schur functions, we simultaneously obtain the decomposition of the $\mathfrak{g}(\myA_{n-1})$-module $V[L]$ as a direct sum of irreducible modules. 
Of course, these decomposition results are consequences of the Littlewood--Richardson rule. 
Below, we state a more general version of this result obtained by Zelevinsky and Stembridge, which we call the `ZS Rule' and formulate the statement in terms of GT parallelograms. 
Fix a semistandard $n$-tableau $X$ of skew shape $\mysmallP/\mysmallQ$ and its corresponding GT $n$-parallelogram $\xelt$. 
For a fixed $i \in \{1,\ldots,n\}$, an $i$-{\em string} of $X$ is the set of all cells on a given row whose entries are $i$. 
A {\em right-to-left (RTL) subtableau} $X'$ of $X$ is a collection of some of the cells / cell entries of $X$ obtained by ({\sl i}) reading rows from top to bottom, ({\sl ii}) reading cells from right to left across rows, beginning at rightmost cell on the top row, and ({\sl iii}) including in $X'$ the entire $i$-string of any row $r$ of $X$ whenever $X'$ contains a cell from row $r$ with entry $i$. 
The {\em length} of an RTL subtableau is its number of cells. 
Our tableau $X$ and parallelogram $\xelt$ are {\em ballot admissible}\footnote{The more common term is {\em Littlewood--Richardson (LR) tableau}; see Chapter 7 Appendix 1.3 of \cite{StanText2} for some discussion relating to such terminology.} if, for each $i \in \{1,2,,\ldots,n-1\}$, no RTL subtableau of $X$ has more cells with an `$i+1$' than cells with an `$i$'. 

Next, we consider ballot-admissibility in terms of GT parallelogram entries. 
In the formula that follows, let $g_{i,j} := g_{i,j}(\xelt)$ if $i \in \{0,1,\ldots,n\}$ and $j \in C_{i}$. 
(Ignore $g_{ij}$ in any formulas if $j \not\in C_{i}$.) 
Set 
\[\myd(\xelt;k,k-p) := \sum_{q=0}^{p-1}\left(\rule[-2.5mm]{0mm}{6mm}2g_{k,k-q}-g_{k-1,k-1-q} - g_{k+1,k+1-q}\right) + g_{k,k-p} - g_{k+1,k+1-p},\] 
where it is understood that $k \in \{1,\ldots,n-1\}$ and $p \in \{0,\ldots,m-1\}$.  
It is not hard to see that $\xelt$ is ballot admissible if and only if $\myd(\xelt;k,k-p) \geq 0$ for all $k$ and $p$.  
The preceding notion generalizes as follows. 
Take an $\myA_{n-1}$-dominant weight $\nu = \sum_{i=1}^{n-1}\,\nu_{i}\omega_{i}$. 
Say the GT $n$-parallelogram $\zelt$ or corresponding $n$-tableau $Z$ is $\nu${\em -ballot admissible} if, for all $k$ and $p$, we have $\myd(\zelt;k,k-q) \geq -\nu_{k}$. 
Equivalently, $Z$ and $\zelt$ are $\nu$-ballot admissible if $\wt(Z')+\nu$ is dominant for all RTL subtableaux $Z'$ of $Z$, where $\wt(Z') := \sum_{i=1}^{n-1}\mym_{i}(Z')\omega_{i}$ with $\mym_{i}(Z') := \#_{i}(Z') - \#_{i+1}(Z')$. 

In this notation, we have the following well-known result interpreted as claims about the weight- and rank-generating functions of skew-tabular lattices. 

\noindent 
{\bf \LRCorollaryOne\ (The Zelevinsky-Stembridge Rule, and a Specialization)}\ \ {\sl Keeping the above language and notation concerning} $L:=L_{\mytinyA_{n-1}}^{\mbox{\tiny skew}}(\mysmallP/\mysmallQ)$, {\sl the quantity} $\vartheta_{_{\mytinyP/\mytinyQ}} = \WGF(L)$ {\sl satisfies}  
\[\mathcal{A}(\myvarZ^{\varrho+\nu})\vartheta_{_{\mytinyP/\mytinyQ}} = \sum_{\zelt \in \mathcal{Z}}\mathcal{A}(\myvarZ^{\nu+\smallwt(\zelt)})
\hspace*{0.35in}\mbox{\sl and}\hspace*{0.35in}
\chi_{_{\nu}}\vartheta_{_{\mytinyP/\mytinyQ}} = \sum_{\zelt \in \mathcal{Z}}\chi_{_{\nu+\tinywt(\zelt)}}.\]
{\sl Now take $\nu = 0$.}  
{\sl Let $\melt$ denote the unique maximal element of $L$, and, in accordance with \TypeAFormulas, set $\lambda_{i}^{j}(\xelt) := \sum_{k=i}^{j}\mym_{k}(\xelt)$ for any $\xelt$ in $L$.  
Then,} 
\[\RGF(L) = \sum_{\zelt \in \mathcal{Z}} \left(q^{\rho(\melt)-\rho(\zelt)}\prod_{i=1}^{n-1}\prod_{j=i}^{n-1}\frac{[\lambda_{i}^{j}(\zelt)+j+1-i]_{q}}{[j+1-i]_{q}}\right)\] 
{\sl is a symmetric and unimodal polynomial in the variable $q$.} 

{\em Proof:} See \cite{Z} and \cite{StemLR} for proofs of the first claim of the statement. 
The second claim concerning the rank-generating function of $L$ therefore follows from \CombinatorialPropositions.\hfill\QED

\vspace*{0.1in}
\noindent {\bf \S \GTStatementsSection\ Skew-tabular lattices as representation diagrams and supporting graphs.}
In this section we state and prove most of our main representation theoretic results about skew-tabular lattices, continuing with the notation of the previous section. 
Let $\selt \myarrow{i} \telt$ be an edge from  $L_{\mytinyA_{n-1}}^{\mbox{\tiny skew}}(\mysmallP/\mysmallQ)$ for GT $n$-parallelograms $\selt$ and $\telt$, and let $(i,j)$ be the position wherein $g_{i,j}(\selt) + 1 = g_{i,j}(\telt)$. 
We assign to this edge a scalar pair $(\myqX_{\telt,\selt},\myqY_{\selt,\telt})$. 
In the formulas for $\myqX_{\telt,\selt}$ and $\myqY_{\selt,\telt}$ below, we use $g_{p,q}$ as shorthand for $g_{p,q}(\telt)$ for any position $(p,q)$. 
Declare that 
\begin{eqnarray*}
\rule[-14mm]{0.0mm}{25mm}\myqX_{\telt,\selt} & := & -\ \frac{\mbox{$\displaystyle \prod_{k\in C_{i+1}}$} \left(\rule[-2.25mm]{0mm}{5.75mm}g_{i,j}-g_{i+1,k}+j-k\right)}{\mbox{$\displaystyle \prod_{k\in C_{i} \setminus \{j\}}$} \left(\rule[-2.25mm]{0mm}{5.75mm}g_{i,j}-g_{i,k}+j-k-1\right)}\\
\setlength{\unitlength}{1cm}
\begin{picture}(0,0)
\put(-4.4,1.5){\LARGE $\left(\star\right)$}
\end{picture}\myqY_{\selt,\telt} & := & \frac{\mbox{$\displaystyle \prod_{k\in  C_{i-1}}$} \left(\rule[-2.25mm]{0mm}{5.75mm}g_{i,j}-g_{i-1,k}+j-k-1\right)}{\mbox{$\displaystyle \prod_{k\in C_{i} \setminus \{j\}}$} \left(\rule[-2.25mm]{0mm}{5.75mm}g_{i,j}-g_{i,k}+j-k\right)}
\end{eqnarray*}
We discovered these formulas by naturally modifying the edge coefficient formulas for classical GT lattices as first recorded by Gelfand and Tsetlin in \cite{GT}. 
The following observations are stated for the record; their proofs are routine and therefore omitted. 

\noindent
{\bf \ObservationalLemma}\ \ {\sl In the notation of the preceding paragraph, we have: (1) The numbers} $\myqX_{\telt,\selt}$ {\sl and} $\myqY_{\selt,\telt}$ {\sl are positive and rational. (2) If $g_{i,k} = g_{i+1,k+1}$, then $g_{i,j}-g_{i+1,k+1}+j-(k+1) = g_{i,j}-g_{i,k}+j-k-1$, so removing these respective factors from the numerator and denominator of the formula for} $\myqX_{\telt,\selt}$ {\sl does not alter its value; similarly we have $g_{i,j}-g_{i-1,k-1}+j-(k-1)-1 = g_{i,j}-g_{i,k}+j-k$ when $g_{i,k} = g_{i-1,k-1}$, so removing these respective factors from the numerator and denominator of the formula for} $\myqY_{\selt,\telt}$ {\sl does not alter its value.} 

The classical-case formulas can be expressed as follows. 
Say $\mysmallQ = (0,\ldots,0)$ with $m \geq n$, and consider $L_{\mytinyA_{n-1}}^{\mbox{\tiny skew}}(\mysmallP/\mysmallQ) = L_{\mytinyA_{n-1}}^{\mbox{\tiny GT}}(\mysmallP)$. 
For $1 \leq i \leq n$, let $C^{\triangleleft}_{i} := \{1,2,\ldots,i\}$, a modified set of indices for the non-inert positions in the $i^{\mbox{\tiny th}}$ column of any associated GT $n$-triangle, cf.\ \GTTriangleFigure. Set $C^{\triangleleft}_{0} := \emptyset$.  
From \cite{MoLast} or \cite{HL}, the classical-case formulas for the coefficients on our edge $\selt \myarrow{i} \telt$, which we denote `$\myqX^{\triangleleft}_{\telt,\selt}$' and `$\myqY^{\triangleleft}_{\selt,\telt}$', are the same as $\left(\star\right)$ above, but with each instance of $C_{l}$, for $0 \leq l \leq n$, replaced by $C^{\triangleleft}_{l}$. 
Unfortunately, it is not always the case that $(\myqX_{\telt,\selt},\myqY_{\selt,\telt}) = (\myqX^{\triangleleft}_{\telt,\selt},\myqY^{\triangleleft}_{\selt,\telt})$, as we see with edge $\selt = {\VertexTableauForText{1}{1}{3}{2}{3}} \mylongarrow{2} {\VertexTableauForText{1}{1}{3}{2}{2}} = \telt$ from  
$L_{\mytinyA_{2}}^{\mbox{\tiny skew}}({\setlength{\unitlength}{0.125cm}\begin{picture}(3,0)\put(0,-0.25){\line(0,1){2}} \put(1,-0.25){\line(0,1){2}} \put(2,-0.25){\line(0,1){2}} \put(3,0.75){\line(0,1){1}} \put(0,-0.25){\line(1,0){2}} \put(0,0.75){\line(1,0){3}} \put(0,1.75){\line(1,0){3}}\end{picture}}) = 
L_{\mytinyA_{2}}^{\mbox{\tiny GT}}({\setlength{\unitlength}{0.125cm}\begin{picture}(3,0)\put(0,-0.25){\line(0,1){2}} \put(1,-0.25){\line(0,1){2}} \put(2,-0.25){\line(0,1){2}} \put(3,0.75){\line(0,1){1}} \put(0,-0.25){\line(1,0){2}} \put(0,0.75){\line(1,0){3}} \put(0,1.75){\line(1,0){3}}\end{picture}})$, cf.\ \SolitarySigmaFig. 
Taking $m = n = 3$, the associated GT $3$-parallelograms are $\selt =$ 
\parbox[c]{1.2cm}{\setlength{\unitlength}{0.3cm}
\begin{picture}(3,5)
\put(-0.5,2.25){
\begin{picture}(0,0)
\put(0,-2){\begin{picture}(0,0)
\put(0.25,3.5){\tiny 0} \put(1.25,3){\tiny 0} \put(2.25,2.5){\tiny 0} \put(3.25,2){\tiny 0} 
\put(0.25,2.5){\tiny 0} \put(1.25,2){\tiny 0} \put(2.25,1.5){\tiny 1} \put(3.25,1){\tiny 2} 
\put(0.25,1.5){\tiny 0} \put(1.25,1){\tiny 2} \put(2.25,0.5){\tiny 2} \put(3.25,0){\tiny 3} 
\end{picture}}
\end{picture}
}
\end{picture}
}
and $\telt =$ 
\parbox[c]{1.2cm}{\setlength{\unitlength}{0.3cm}
\begin{picture}(3,5)
\put(-0.5,2.25){
\begin{picture}(0,0)
\put(0,-2){\begin{picture}(0,0)
\put(0.25,3.5){\tiny 0} \put(1.25,3){\tiny 0} \put(2.25,2.5){\tiny 0} \put(3.25,2){\tiny 0} 
\put(0.25,2.5){\tiny 0} \put(1.25,2){\tiny 0} \put(2.25,1.5){\tiny 2} \put(3.25,1){\tiny 2} 
\put(0.25,1.5){\tiny 0} \put(1.25,1){\tiny 2} \put(2.25,0.5){\tiny 2} \put(3.25,0){\tiny 3} 
\end{picture}}
\end{picture}
}
\end{picture}
}. 
We have  $(\myqX_{\telt,\selt},\myqY_{\selt,\telt}) = (\frac{3}{2},2)$ and $(\myqX^{\triangleleft}_{\telt,\selt},\myqY^{\triangleleft}_{\selt,\telt}) = (3,1)$. 

In \cite{HL}, Hersh and Lenart provide edge-coefficient formulas for classical GT lattices that are an alternative to the formulas of \cite{GT}. 
Fix an edge-color $i$ and an edge $\selt \myarrow{i} \telt$ from $L_{\mytinyA_{n-1}}^{\mbox{\tiny skew}}(\mysmallP/\mysmallQ) = L_{\mytinyA_{n-1}}^{\mbox{\tiny GT}}(\mysmallP)$. 
In particular, there is a fixed number $j \in \{1,2,\ldots,i\}$ such that $g_{i,j}(\selt)+1=g_{i,j}(\telt)$ and such that $\selt$ and $\telt$ agree in all other array positions. 
In \cite{HL}, formulas for the coefficients on our fixed edge, here denoted $(\myqX^{\mbox{\tiny tab}}_{\telt,\selt},\myqY^{\mbox{\tiny tab}}_{\selt,\telt})$, are expressed in terms of the associated tableaux. 
Rather than re-produce definitions from \cite{HL} for the tableaux-related quantities used in their formulas, we offer here interpretations of these quantities in terms of the non-inert positions of the GT $n$-triangle associated with $\telt$. 
Say an index $r \in \{1,2,3,\ldots,i\}$ is {\em $i$-significant in} $\telt$ if $r=1$ or else if $g_{i,i+1-r} < g_{i,i+2-r}$.  
Let $p$ be the number of $i$-significant indices in $\telt$; the subsequence $(r_{1},r_{2},\ldots,r_{p})$ of the set $\{1,2,3,\ldots,i\}$ will denote the strictly increasing sequence of $i$-significant indices in $\telt$.  
Since the number $g_{i,j}$ in the array for $\telt$ may be decreased by one, then $j = i+1-r_{k}$ for some $i$-significant index $r_{k}$ in $\telt$. 
Next, we consider certain differences in array values in or near the $i^{\mbox{\tiny th}}$ column of the GT $n$-triangle associated with $\telt$. 
For each $l \in \{1,2,\ldots,p-1\}$, we let $a_{l}' := r_{l+1}+g_{i,i}-g_{i,i+2-r_{l+1}}$ and $b_{l}' := g_{i,i+2-r_{l+1}} - g_{i-1,i+1-r_{l+1}}$.
Declare that $a_{k} := r_{k}+g_{i,i}-g_{i,j}+1$ and $b_{k} := g_{i+1,j+1}-g_{i,j}-1$. 
For $l \in \{1,2,\ldots,p\} \setminus \{k\}$, set $a_{l} := r_{l}+g_{i,i}-g_{i,i+1-r_{l}}$ and $b_{l} := g_{i+1,i+2-r_{l}} - g_{i,i+1-r_{l}}$. 
The quantities defined in the previous two sentences are interpreted in \cite{HL} in terms of tableaux.  
The following are formulas (13) and (14) from \S 3 of \cite{HL}:
\begin{eqnarray*}
\rule[-9mm]{0.0mm}{15mm}\myqX^{\mbox{\tiny tab}}_{\telt,\selt} & := & b_{k} {\mbox{$\displaystyle \prod_{l=1}^{k-1}$} \left(\rule[-2.25mm]{0mm}{5.75mm}1+\frac{b_{l}}{a_{k}-a_{l}}\right)}{\mbox{$\displaystyle \prod_{l=k+1}^{p}$} \left(\rule[-2.25mm]{0mm}{5.75mm}1-\frac{b_{l}}{a_{l}-a_{k}}\right)}\\
\setlength{\unitlength}{1cm}
\begin{picture}(0,0)
\put(-3.9,0.85){\LARGE $\left(\star\star\right)$}
\end{picture}
\myqY^{\mbox{\tiny tab}}_{\selt,\telt} & := & b_{k}' {\mbox{$\displaystyle \prod_{l=1}^{k-1}$} \left(\rule[-2.25mm]{0mm}{5.75mm}1-\frac{b_{l}'}{a_{k}'-a_{l}'}\right)}{\mbox{$\displaystyle \prod_{l=k+1}^{p-1}$} \left(\rule[-2.25mm]{0mm}{5.75mm}1+\frac{b_{l}'}{a_{l}'-a_{k}'}\right)}
\end{eqnarray*}
For example, for the aforementioned edge $\selt = {\VertexTableauForText{1}{1}{3}{2}{3}} \mylongarrow{2} {\VertexTableauForText{1}{1}{3}{2}{2}} = \telt$ from  
$L_{\mytinyA_{2}}^{\mbox{\tiny skew}}({\setlength{\unitlength}{0.125cm}\begin{picture}(3,0)\put(0,-0.25){\line(0,1){2}} \put(1,-0.25){\line(0,1){2}} \put(2,-0.25){\line(0,1){2}} \put(3,0.75){\line(0,1){1}} \put(0,-0.25){\line(1,0){2}} \put(0,0.75){\line(1,0){3}} \put(0,1.75){\line(1,0){3}}\end{picture}}) = 
L_{\mytinyA_{2}}^{\mbox{\tiny GT}}({\setlength{\unitlength}{0.125cm}\begin{picture}(3,0)\put(0,-0.25){\line(0,1){2}} \put(1,-0.25){\line(0,1){2}} \put(2,-0.25){\line(0,1){2}} \put(3,0.75){\line(0,1){1}} \put(0,-0.25){\line(1,0){2}} \put(0,0.75){\line(1,0){3}} \put(0,1.75){\line(1,0){3}}\end{picture}})$, we have 
$(\myqX^{\mbox{\tiny tab}}_{\telt,\selt},\myqY^{\mbox{\tiny tab}}_{\selt,\telt}) = (\frac{3}{2},2) = (\myqX_{\telt,\selt},\myqY_{\selt,\telt})$. 
The following result generalizes the latter observation. 

\noindent
{\bf \ConsistencyLemma}\ \ {\sl Continue with the assumptions of the previous three paragraphs that} $\mysmallQ = (0,\ldots,0)$ {\sl with} $m \geq n$. 
{\sl Suppose $\selt \myarrow{i} \telt$ in} $L_{\mytinyA_{n-1}}^{\mbox{\tiny skew}}(\mysmallP/\mysmallQ) = L_{\mytinyA_{n-1}}^{\mbox{\tiny GT}}(\mysmallP)$. 
{\sl Then} $\myqX^{\mbox{\tiny tab}}_{\telt,\selt} = \myqX_{\telt,\selt}$ {\sl and} $\myqY^{\mbox{\tiny tab}}_{\selt,\telt} = \myqY_{\selt,\telt}$. 

{\em Proof.} We begin by viewing $\selt$ and $\telt$ as GT $n$-triangles; then we rewrite the formulas $\left(\star\star\right)$ for $\myqX^{\mbox{\tiny tab}}_{\telt,\selt}$ and $\myqY^{\mbox{\tiny tab}}_{\selt,\telt}$ in terms of the non-inert positions of $\telt$; and then we write each formula in terms of the entries of $\telt$ when the latter is viewed as a GT $n$-parallelogram.  Now, by routine algebra, \[\myqX^{\mbox{\tiny tab}}_{\telt,\selt} = -\prod_{q \in \{2,\ldots,i+1\}}(g_{i,j}-g_{i+1,q}+j-q) \mbox{\LARGE $/$}\hspace*{-0.1in} \prod_{q \in \{2,\ldots,i+1\}\setminus\{j\}}(g_{i,j}-g_{i,q}+j-q-1).\] Multiply the right-hand side by $\displaystyle \prod_{q=i+1-(m-1)}^{1}(g_{i,j}-g_{i+1,q}+j-q) \mbox{\LARGE $/$}\hspace*{-0.1in} \prod_{q=i-(m-1)}^{0}(g_{i,j}-g_{i,q}+j-q-1)$, which is unity, to get $\myqX_{\telt,\selt}$. Similarly see that
\[\myqY^{\mbox{\tiny tab}}_{\selt,\telt} = \prod_{q \in \{1,\ldots,i-1\}}(g_{i,j}-g_{i-1,q}+j-q-1) \mbox{\LARGE $/$}\hspace*{-0.1in} \prod_{q \in \{2,\ldots,i\}\setminus\{j\}}(g_{i,j}-g_{i,q}+j-q).\] 
Now multiply the right-hand side by $\displaystyle \prod_{q=i-1-(m-1)}^{0}(g_{i,j}-g_{i-1,q}+j-q-1) \mbox{\LARGE $/$}\hspace*{-0.1in} \prod_{q=i-(m-1)}^{1}(g_{i,j}-g_{i,q}+j-q)$, which is unity, to get $\myqY_{\selt,\telt}$.\hfill\QED

The next theorem renders a central result of \cite{HL} -- viz.\ Theorem 4.1 --  
in the language of this paper. 
\GTBasisTheorem\ can be interpreted as saying that for non-skew shapes, the associated classical GT lattice is a supporting graph for some $\mathfrak{g}(\myA_{n-1})$-module. 
We note that the insights of \cite{PrGZ} afforded the crucial connection between the formulas of \cite{GT} and the classical GT lattices as studied in \cite{HL}; this  accounts for the many authorial attributions of the next result. 

\noindent
{\bf \GTBasisTheorem\ (Gelfand--Tsetlin/Proctor/Hersh--Lenart)}\ \ {\sl Assume} $\mysmallQ = (0,\ldots,0)$ {\sl and $m=n$. 
Let} $L := L_{\mytinyA_{n-1}}^{\mbox{\tiny skew}}(\mysmallP/\mysmallQ) = L_{\mytinyA_{n-1}}^{\mbox{\tiny GT}}(\mysmallP)$, {\sl and assign the scalar pair} $(\myqX^{\mbox{\tiny tab}}_{\telt,\selt},\myqY^{\mbox{\tiny tab}}_{\selt,\telt})$ {\sl of $\left(\star\star\right)$ to each edge $\selt \rightarrow \telt$ in $L$. 
Then $L$ is} $\myA_{n-1}${\sl -structured and the scalars satisfy all DC relations. 
Therefore the action of the generators of} $\mathfrak{g}(\myA_{n-1})$ {\sl on the vector space $V[L]$ as defined by the formulas} ($\blacklozenge$) {\sl in} \S \SetupSectionTwo\ {\sl is well-defined, $\{v_{\telt}\}_{\telt \in L}$ is a weight basis for the} $\mathfrak{g}(\myA_{n-1})${\sl -module $V[L]$, and the lattice $L$ together with the set of scalar pairs} $\{(\myqX^{\mbox{\tiny tab}}_{\telt,\selt},\myqY^{\mbox{\tiny tab}}_{\selt,\telt})\}_{\selt \rightarrow \telt \mbox{\scriptsize \ in } L}$ {\sl is its representation diagram.} 

Now we return to the case of general skew shapes. The goal of this section is to demonstrate the following result.

\noindent
{\bf \GTGeneralTheorem}\ \ 
{\sl Let} $L := L_{\mytinyA_{n-1}}^{\mbox{\tiny skew}}(\mysmallP/\mysmallQ)$, {\sl and assign the scalar pairs $\{(\myqX_{\telt,\selt},\myqY_{\selt,\telt})\}_{\selt \rightarrow \telt \mbox{\scriptsize \ in } L}$ of $\left(\star\right)$ above to the edges of $L$. 
Then $L$ is} $\myA_{n-1}${\sl -structured and the scalars satisfy all DC relations. 
Therefore the action of the generators of} $\mathfrak{g}(\myA_{n-1})$ {\sl on the vector space $V[L]$ as defined by the formulas} ($\blacklozenge$) {\sl in} \S \SetupSectionTwo\ {\sl is well-defined, $\{v_{\telt}\}_{\telt \in L}$ is a weight basis for the} $\mathfrak{g}(\myA_{n-1})${\sl -module $V[L]$, and the lattice $L$ together with $\{(\myqX_{\telt,\selt},\myqY_{\selt,\telt})\}_{\selt \rightarrow \telt \mbox{\scriptsize \ in } L}$ is its representation diagram.}

Our first proof of \GTGeneralTheorem\ was a direct verification of the DC relations. 
However, we subsequently realized that $L_{\mytinyA_{n-1}}^{\mbox{\tiny skew}}(\mysmallP/\mysmallQ)$ can be identified as a component within some larger classical GT lattice $L_{\mytinyA_{m+n-1}}^{\mbox{\tiny GT}}(\mysmallP')$ in such a way that scalar pairs on corresponding edges are exactly the same. 
\GTGeneralTheorem\ follows immediately from this observation. 
This is the proof of \GTGeneralTheorem\ that we formally pursue below, using GT parallelograms. 
Before we do so, we note that an `embedding' of our given skew-tabular lattice $K:=L_{\mytinyA_{n-1}}^{\mbox{\tiny skew}}(\mysmallP/\mysmallQ)$ in some classical GT lattice $L:=L_{\mytinyA_{m+n-1}}^{\mbox{\tiny GT}}(\mysmallP')$ is easy to visualize via tableaux: The classical non-skew shape is just the shape associated with $\mysmallP$, tableaux entries are from $\{1,2,\ldots,m+n\}$, and the topmost tableau $T'$ of the $\{m+1,\ldots,m+n-1\}$-component of $L$ that is isomorphic to (a re-coloring of) $K$ has the same entries in its non-skew-shape cells as the topmost tableau $T$ in $L$ but the entry in each cell of the skew-shape part of $T'$ is the entry from $T$ increased by $m$. 
We thank the referee for sharing this latter observation with us. 

Using GT parallelograms, we codify the embedding from the previous paragraph as follows. 
Let $\mysmallP' := (\myscriptsizeP'_{1},\ldots,\myscriptsizeP'_{m+n}) := (\myscriptsizeP_{1},\ldots,\myscriptsizeP_{m},0,\ldots,0)$ and $\mysmallQ' := (\myscriptsizeQ'_{1},\ldots,\myscriptsizeQ'_{m+n}) := (0,\ldots,0)$.
View an $n$-parallelogram $\xelt$ framed by $\mysmallP/\mysmallQ$ as an array $(g_{i,j}(\xelt))$ where $i \in \{0,1,\ldots,n\}$ and $j \in C_{i}$.  
Now create a new $(m+n)$-triangle $\phi(\xelt)$ framed by $\mysmallP'$ by taking 
\[g'_{k,l}(\phi(\xelt)) := \left\{\begin{array}{cc} 
g_{0,l-m}(\xelt) & k \in \{1,\ldots,m-1\}, l \in \{1,\ldots,k\}\\
0 & k \in \{m+1,\ldots,m+n\}, l \in \{1,\ldots,k-m\}\\
g_{k-m,l-m}(\xelt) & k \in \{m,\ldots,m+n\}, l \in \{k-m+1,\ldots,k\} 
\end{array}\right.\] 
for $1 \leq l \leq k \leq m+n$. 
To view $\phi(\xelt)$ as an $(m+n)$-parallelogram framed by $\mysmallP'/\mysmallQ'$, set $g'_{k,l}(\phi(\xelt)) := 0$ when $k \in \{0,\ldots,m+n\}, l \in \{k-(m+n),\ldots,0\}$; in this case the array $\phi(\xelt) = \left(\rule[-2.5mm]{-0.1mm}{6mm}g'_{k,l}(\phi(\xelt))\right)$ is defined for all $k \in \{0,1,\ldots,m+n\}$ and $l \in C'_{k} := \{k,k-1,\ldots,k-(m+n-1)\}$. 
It is clear that we have simply built a larger triangle and parallelogram around the original parallelogram without losing any of the original information. 
So the function $\phi: L_{\mytinyA_{n-1}}^{\mbox{\tiny skew}}(\mysmallP/\mysmallQ) \longrightarrow L_{\mytinyA_{m+n-1}}^{\mbox{\tiny GT}}(\mysmallP')$ is an injection, and it is evident that $\selt \myarrow{i} \telt$ in  $L_{\mytinyA_{n-1}}^{\mbox{\tiny skew}}(\mysmallP/\mysmallQ)$ if and only if $\phi(\selt) \mylongarrow{\mbox{\tiny $i\!+\!m$}} \phi(\telt)$. 

\noindent 
{\bf \StandardizeProp}\ \ {\sl Keep the notation of the preceding paragraph. 
Define $\sigma: \{1,\ldots,n-1\} \longrightarrow \{m+1,\ldots,m+n-1\} =: J$ by $\sigma(i) := m+i$. 
Let} $L := L_{\mytinyA_{n-1}}^{\mbox{\tiny skew}}(\mysmallP/\mysmallQ)^{\sigma}$ {\sl with maximal element denoted $\melt$, let} $K := L_{\mytinyA_{m+n-1}}^{\mbox{\tiny GT}}(\mysmallP')${\sl , and let $\phi: L \longrightarrow K$ be the mapping of sets identified above. 
Then, $\phi$ is an edge-, edge-color-, and weight-preserving injection and} $\phi(L) = \comp_{J}(\phi(\melt))$. 
{\sl Further, suppose that $\selt \myarrow{i} \telt$ in}  $L_{\mytinyA_{n-1}}^{\mbox{\tiny skew}}(\mysmallP/\mysmallQ)${\sl , and set $\selt' := \phi(\selt)$ and $\telt' := \phi(\telt)$ so that} $\phi(\selt) \mylongarrow{\mbox{\tiny $i\!+\!m$}} \phi(\telt)$. {\sl Then,} $\myqX_{\telt',\selt'} = \myqX_{\telt,\selt}$ {\sl and} $\myqY_{\selt',\telt'} = \myqY_{\selt,\telt}$. 

{\em Proof.} Observations at the end of the paragraph preceding the proposition statement make clear that our mapping $\phi$ is edge- and edge-color-preserving (and therefore weight-preserving) and that $\phi(L) = \comp_{J}(\phi(\melt))$. 
Assume now that for $\selt \myarrow{i} \telt$ in $L_{\mytinyA_{n-1}}^{\mbox{\tiny skew}}(\mysmallP/\mysmallQ)$, we have $g_{i,j}(\selt)+1=g_{i,j}(\telt)$. 
Some notation before we confirm the claimed coefficient equalities: For an integer $k$ and integer set $\mathcal{S}$, let $\mathcal{S}+k = \{s+k\, |\, s \in \mathcal{S}\}$. 
From the formula for $\myqX_{\telt',\selt'}$ at $(\star)$, the numerator is indexed by $C'_{i+1+m}$ and the denominator by $C'_{i+m} \setminus \{j+m\}$. 
From the definition of $\phi$, we see that $g'_{i+1+m,k+1+m}=0$ if $k+1+m \in C'_{i+1+m} \setminus (C_{i+1}+m)$. 
But $k+1+m \in C'_{i+1+m} \setminus (C_{i+1}+m)$ if and only if $k+m \in C'_{i+m} \setminus (C_{i}+m)$, in which case $g'_{i+m,k+m}=0$. 
This observation allows us to invoke \ObservationalLemma.2 in the second equality in what follows: 
\begin{eqnarray*}
\rule[-14mm]{0.0mm}{25mm}\myqX_{\telt',\selt'}  & = & -\ \frac{\mbox{$\displaystyle \prod_{k+1+m\in C'_{i+1+m}}$} \left(\rule[-2.25mm]{0mm}{5.75mm}g'_{i+m,j+m}-g'_{i+1+m,k+1+m}+j-k-1\right)}{\mbox{$\displaystyle \prod_{k+m\in C'_{i+m} \setminus \{j+m\}}$} \left(\rule[-2.25mm]{0mm}{5.75mm}g'_{i+m,j+m}-g'_{i+m,k+m}+j-k-1\right)}\\ 
 & = & -\ \frac{\mbox{$\displaystyle \prod_{k+1+m\in (C_{i+1}+m)}$} \left(\rule[-2.25mm]{0mm}{5.75mm}g'_{i+m,j+m}-g'_{i+1+m,k+1+m}+j-k-1\right)}{\mbox{$\displaystyle \prod_{k+m\in (C_{i}+m) \setminus \{j+m\}}$} \left(\rule[-2.25mm]{0mm}{5.75mm}g'_{i+m,j+m}-g'_{i+m,k+m}+j-k-1\right)}\\
 & = & -\ \frac{\mbox{$\displaystyle \prod_{k+1\in C_{i+1}}$} \left(\rule[-2.25mm]{0mm}{5.75mm}g_{i,j}-g_{i+1,k+1}+j-k-1\right)}{\mbox{$\displaystyle \prod_{k\in C_{i} \setminus \{j\}}$} \left(\rule[-2.25mm]{0mm}{5.75mm}g_{i,j}-g_{i,k}+j-k-1\right)},
\end{eqnarray*}
and the latter is just $\myqX_{\telt,\selt}$. Similarly see that $\myqY_{\selt',\telt'} = \myqY_{\selt,\telt}$.\hfill\QED

{\em Proof of \GTGeneralTheorem.} In the notation of \StandardizeProp, let $\mathfrak{g}'$ be the Lie subalgebra of $\mathfrak{g}(\myA_{m+n-1})$ generated by $\{\myqx_{j},\myqy_{j}\}_{j \in J}$, so $\mathfrak{g}'$ is naturally isomorphic to $\mathfrak{g}(\myA_{n-1})$.  
Since, by \GTBasisTheorem, $K := L_{\mytinyA_{m+n-1}}^{\mbox{\tiny GT}}(\mysmallP')$ is a representation diagram for a $\mathfrak{g}(\myA_{m+n-1})$-module, then 
$L \cong \comp_{J}(\phi(\melt))$ is a representation diagram for a $\mathfrak{g}'$-module. 
By \ConsistencyLemma, $\myqX^{\mbox{\tiny tab}}_{\telt',\selt'} = \myqX_{\telt',\selt'}$ and $\myqY^{\mbox{\tiny tab}}_{\selt',\telt'} = \myqY_{\selt',\telt'}$ whenever $\selt' = \phi(\selt) \mylongarrow{\mbox{\tiny $i\!+\!m$}} \phi(\telt) = \telt'$ for some edge $\selt \myarrow{i} \telt$ in $L_{\mytinyA_{n-1}}^{\mbox{\tiny skew}}(\mysmallP/\mysmallQ)$. 
By \StandardizeProp, it now follows that $L_{\mytinyA_{n-1}}^{\mbox{\tiny skew}}(\mysmallP/\mysmallQ)$, with scalar pairs as in the statement of \GTGeneralTheorem, is a representation diagram for a $\mathfrak{g}(\myA_{n-1})$-module.\hfill\QED

To close this section, we record the following corollary, which combines \GTGeneralTheorem\ and \MainCombinatorialTheorem. 
For classical GT lattices, these results are due to Stanley \cite{StanGLn} and Proctor \cite{PrGZ}. 

\noindent
{\bf \GTCombinatorics}\ \ {\sl Skew-tabular lattices are rank symmetric, rank unimodal, and strongly Sperner.}\hfill\QED

\vspace*{0.1in}
\noindent {\bf \S \OrthogonalSection\ Constructions of certain orthogonal Lie algebra representations.} 
 As an application of our perspective on skew-tabular lattices and GT parallelograms, we explicitly construct new weight bases for certain orthogonal Lie algebra representations. 
To do so, we build diamond-colored distributive lattices from GT-like patterns and locate within them large components isomorphic to classical GT lattices.\footnote{In particular, no non-classical skew-tabular lattices are required for the orthogonal representations that follow. However, in \cite{DDW} we use an entirely analogous approach to construct certain irreducible representations of the type $\myE_{6}$ and $\myE_{7}$ simple Lie algebras, and non-classical skew-tabular lattices are required.}  
This will allow us to determine coefficients for lattice edges and so realize each such lattice as a representation diagram. 

The arrays we use to build these distributive lattice representation diagrams are an interpretation of certain combinatorial objects considered by Proctor and Stanley in their proof that all minuscule posets are Gaussian, cf.\ Theorem 6 of \cite{PrEur}. 
Consider a simple Lie algebra $\mathfrak{g}(\myX_{n})$ with $\myX \in \{\myA, \myB, \myC, \myD, \myE, \myF, \myG\}$.  
Let $\omega$ be a dominant minuscule weight, so $\omega$ must be a fundamental weight and $\myX_{n}$ cannot be $\myE_{8}$, $\myF_{4}$, or $\myG_{2}$, cf.\ Exercise 13.4.13 of \cite{Hum}. 
Fix the vertex-colored minuscule poset $P(\omega)$ associated with the pair $(\mathfrak{g}(\myX_{n}),\omega)$, cf.\ \cite{PrEur}, and let $L(\omega) := \mathbf{J}_{\mbox{\tiny color}}(P(\omega))$ be the diamond-colored distributive lattice of order ideals from $P(\omega)$. 
Since, in the notation of \cite{DonPosetModels}, the edge-colored directed graph $L(\omega)$ can be identified with the corresponding weight diagram $\Pi(\omega)$, then any element $\xelt \in L(\omega)$ can be identified with its weight $\wt(\xelt)$. 
Now let $m$ be any positive integer. 
Following \cite{PrEur}, we say an $m\omega$ {\em dual-multichain} $\mathcal{M}$ is a sequence of elements $\mathcal{M} := (\xelt_{1},\ldots,\xelt_{m})$ of $L(\omega)$ such that $\xelt_{1} \geq \cdots \geq \xelt_{m}$. 
The weight $wt(\mathcal{M})$ of $\mathcal{M}$ is declared to be $\wt(\xelt_{1}) + \cdots + \wt(\xelt_{m})$. 
Within this framework, Lemma 6 of \cite{PrEur} interprets the result, due to Seshadri (see \cite{Seshadri}, \cite{LS}), that 
\[\chi_{_{m\omega}} = \stackrel{\mbox{\LARGE $\sum$}}{\parbox[t]{0.75in}{\tiny \begin{center}\vspace*{-0.125in}$m\omega\ \!$ dual-\\ multichains$\ \! \mathcal{M}$\end{center}}}\!\!\!\!\myvarZ^{wt(\mathcal{M})}.\] 

\begin{figure}[t] 
\begin{center}
\OddOrthPatternFigure: Inequalities for an odd orthogonal ideal array.

\vspace*{0.25in}
\begin{tabular}{cccccccccc}
 & & & & & & & & & $c_{n,1}$ \\
 & & & & & & & & \begin{rotate}{45}$\geq$\end{rotate} & \\
 & & & & & & & $c_{n-1,1}$ & & \\
 & & & & & & \begin{rotate}{45}$\geq$\end{rotate} & & \begin{rotate}{315}$\!\!\!\!\leq$\end{rotate} & \\
 & & & & & $c_{n-2,1}$ & & & & $c_{n,2}$ \\
 & & & & \begin{rotate}{45}$\geq$\end{rotate} & & \begin{rotate}{315}$\!\!\!\!\leq$\end{rotate} & & \begin{rotate}{45}$\geq$\end{rotate} & \\
 & & & {\LARGE $\iddots$} & & & {\LARGE $\cdot$} & {\LARGE $\cdot$} & {\LARGE $\cdot$} & \\
 & \begin{rotate}{45}$\geq$\end{rotate} & & & & & & & & \\
$c_{1,1}$ & & & & & & {\LARGE $\cdot$} & {\LARGE $\cdot$} & {\LARGE $\cdot$} & \\
 \vspace*{-0.1in}
 & \begin{rotate}{315}$\!\!\!\!\leq$\end{rotate} & & & & & & & & \\
 & & & {\LARGE $\ddots$} & & & {\LARGE $\cdot$} & {\LARGE $\cdot$} & {\LARGE $\cdot$} & \\
 & & & & \begin{rotate}{315}$\!\!\!\!\leq$\end{rotate} & & \begin{rotate}{45}$\geq$\end{rotate} & & \begin{rotate}{315}$\!\!\!\!\leq$\end{rotate} & \\
 & & & & & $c_{n-2,n-2}$ & & & & $c_{n,n-1}$ \\
 & & & & & & \begin{rotate}{315}$\!\!\!\!\leq$\end{rotate} & & \begin{rotate}{45}$\geq$\end{rotate} & \\
 & & & & & & & $c_{n-1,n-1}$ & & \\
 & & & & & & & & \begin{rotate}{315}$\!\!\!\!\leq$\end{rotate} & \\
 & & & & & & & & & $c_{n,n}$ 
\end{tabular}
\end{center}

\vspace*{-0.3in}
\end{figure} 

\noindent 
Now, any $m\omega$ dual-multichain $\mathcal{M} = (\xelt_{1},\ldots,\xelt_{m})$ corresponds to an array $\big(c_{v}(\mathcal{M})\big)_{v \in P(\omega)}$ where, for any $v \in P(\omega)$, we set $c_{v} := \rule[-2.6mm]{0.2mm}{7.1mm}\left\{i \in \{1,\ldots,m\}\, \rule[-2.6mm]{0.1mm}{7.1mm}\, v \in \xelt_{i}\right\}\rule[-2.6mm]{0.2mm}{7.1mm}$. 
Notice that $0 \leq c_{v}(\mathcal{M}) \leq c_{u}(\mathcal{M}) \leq m$ whenever $u \leq v$ in $P(\omega)$. 
This is how we obtained our descriptions of the GT-like orthogonal ideal arrays we use in this section to define our distributive lattice representation diagrams. 

From here, the details of our set-up mainly follow \S 8 of \cite{DonPosetModels}.
Fix positive integers $m$ and $n$, with $n \geq 3$. 
An {\em odd orthogonal ideal array} of size $n$ and bounded by $m$ is an array $(c_{i,j})_{1 \leq j \leq i \leq n}$ of nonnegative integers with $c_{n,n} \leq m$ and satisfying the inequalities pictured in \OddOrthPatternFigure. 
As described in the second paragraph above, the poset that frames each such array is the minuscule poset associated with the `spin-node' fundamental weight $\omega_n$ associated with the type $\myB_n$ root system. 
Let $L_{\mytinyB_n}(m\omega_n)$ be the collection of odd orthogonal ideal arrays of size $n$ and bounded by $m$. 
For $\telt \in L_{\mytinyB_n}(m\omega_n)$, we use the notation $c_{i,j}(\telt)$ to refer to the $(i,j)$-entry of the corresponding array. 
Partially order the elements of $L_{\mytinyB_n}(m\omega_n)$ by componentwise comparison, so $\selt \leq \telt$ if and only if $c_{i,j}(\selt) \leq c_{i,j}(\telt)$ for all $1 \leq i \leq j \leq n$. 
It is routine to check that the partially ordered set $L_{\mytinyB_n}(m\omega_n)$ is a distributive lattice. 
\hfill Then\ \ $\selt \rightarrow \telt$\ \ in the Hasse diagram for\ \ $L_{\mytinyB_n}(m\omega_n)$\ \ if and only if

\newpage
\begin{figure}[ht] 
\begin{center}
\EvenOrthPatternFigure.E: Inequalities for an even orthogonal ideal array, with $n$ even.

\vspace*{0.25in}
\begin{tabular}{ccccccccccccc}
 & & & & & & & & & & & & $c_{n-1,1}$ \\
 & & & & & & & & & & & \begin{rotate}{45}$\geq$\end{rotate} & \\
 & & & & & & & & & & $c_{n-2,1}$ & & \\
 & & & & & & & & & \begin{rotate}{45}$\geq$\end{rotate} & & \begin{rotate}{315}$\!\!\!\!\leq$\end{rotate} & \\
 & & & & & & & & $c_{n-3,1}$ & & & & $c_{n,1}$ \\
 & & & & & & & \begin{rotate}{45}$\geq$\end{rotate} & & \begin{rotate}{315}$\!\!\!\!\leq$\end{rotate} & & \begin{rotate}{45}$\geq$\end{rotate} & \\
 & & & & & & $c_{n-4,1}$ & & & & $c_{n-2,2}$ & & \\
 & & & & & \begin{rotate}{45}$\geq$\end{rotate} & & \begin{rotate}{315}$\!\!\!\!\leq$\end{rotate} & & \begin{rotate}{45}$\geq$\end{rotate} & & \begin{rotate}{315}$\!\!\!\!\leq$\end{rotate} & \\
 & & & & $c_{n-5,1}$ & & & & $c_{n-3,2}$ & & & & $c_{n-1,2}$ \\
 & & & \begin{rotate}{45}$\geq$\end{rotate} & & \begin{rotate}{315}$\!\!\!\!\leq$\end{rotate} & & \begin{rotate}{45}$\geq$\end{rotate} & & \begin{rotate}{315}$\!\!\!\!\leq$\end{rotate} & & \begin{rotate}{45}$\geq$\end{rotate} & \\
 & & $\ \iddots\ \ $ & & & & & $\cdot$ & $\cdot$ & $\cdot$ & & &  \\
 & \begin{rotate}{45}$\geq$\end{rotate} & & & & & & & & & & & \\
$\ c_{1,1}\ \ $ & & & & & & & $\cdot$ & $\cdot$ & $\cdot$ & & & \\
 & \begin{rotate}{315}$\!\!\!\!\leq$\end{rotate} & & & & & & & & & & \\
 & & $\ \ddots\ \ $ & & & & &  $\cdot$ & $\cdot$ & $\cdot$ & & & \\
 & & & \begin{rotate}{315}$\!\!\!\!\leq$\end{rotate} & & \begin{rotate}{45}$\geq$\end{rotate} & & \begin{rotate}{315}$\!\!\!\!\leq$\end{rotate} & & \begin{rotate}{45}$\geq$\end{rotate} & & \begin{rotate}{315}$\!\!\!\!\leq$\end{rotate} & \\
 & & & & $c_{n-5,n-5}$ & & & & $c_{n-3,n-4}$ & & & & $c_{n-1,\frac{n}{2}-1}$ \\
 & & & & & \begin{rotate}{315}$\!\!\!\!\leq$\end{rotate} & & \begin{rotate}{45}$\geq$\end{rotate} & & \begin{rotate}{315}$\!\!\!\!\leq$\end{rotate} & & \begin{rotate}{45}$\geq$\end{rotate} & \\
 & & & & & & $c_{n-4,n-4}$ & & & & $c_{n-2,n-3}$ & & \\
 & & & & & & & \begin{rotate}{315}$\!\!\!\!\leq$\end{rotate} & & \begin{rotate}{45}$\geq$\end{rotate} & & \begin{rotate}{315}$\!\!\!\!\leq$\end{rotate} & \\
 & & & & & & & & $c_{n-3,n-3}$ & & & & $c_{n,\frac{n}{2}-1}$ \\
 & & & & & & & & & \begin{rotate}{315}$\!\!\!\!\leq$\end{rotate} & & \begin{rotate}{45}$\geq$\end{rotate} & \\
 & & & & & & & & & & $c_{n-2,n-2}$ & & \\
 & & & & & & & & & & & \begin{rotate}{315}$\!\!\!\!\leq$\end{rotate} & \\
 & & & & & & & & & & & & $c_{n-1,\frac{n}{2}}$ 
\end{tabular}
\end{center}
\end{figure} 

\noindent
there is a pair $(i,j)$ such that $c_{p,q}(\selt) = c_{p,q}(\telt)$ when $(p,q) \not= (i,j)$ and $c_{i,j}(\selt) + 1 = c_{i,j}(\telt)$. 
In this case, we attach the color $i$ to this edge of the Hasse diagram and write $\selt \myarrow{i} \telt$. 
In this way we realize $L_{\mytinyB_n}(m\omega_n)$ as a diamond-colored distributive lattice.

Now let $n \geq 4$. 
When $n$ is even, we define an {\em even orthogonal ideal array} of size $n$ and bounded by $m$ to be an array $(c_{i,j})_{1 \leq j \leq i \leq n-2} \cup (c_{n-1,j})_{1 \leq j \leq \frac{n}{2}} \cup (c_{n,j})_{1 \leq j \leq \frac{n}{2}-1}$ of nonnegative integers with $c_{n-1,\frac{n}{2}} \leq m$ and satisfying the inequalities pictured in \EvenOrthPatternFigure.E.   
When $n$ is odd, we define an {\em even orthogonal ideal array} of size $n$ and bounded by $m$ to be an array $(c_{i,j})_{1 \leq j \leq i \leq n-2} \cup (c_{n-1,j})_{1 \leq j \leq \frac{n-1}{2}} \cup (c_{n,j})_{1 \leq j \leq \frac{n-1}{2}}$ of nonnegative integers with $c_{n,\frac{n-1}{2}} \leq m$ and satisfying the inequalities pictured in \EvenOrthPatternFigure.O. 
We use $L_{\mytinyD_n}(m\omega_{n-1})$ to denote the collection of even orthogonal ideal patterns of size

\newpage
\begin{figure}[ht] 
\begin{center}
\EvenOrthPatternFigure.O: Inequalities for an even orthogonal ideal pattern, with $n$ odd.

\vspace*{0.25in}
\begin{tabular}{ccccccccccccc}
 & & & & & & & & & & & & $c_{n-1,1}$ \\
 & & & & & & & & & & & \begin{rotate}{45}$\geq$\end{rotate} & \\
 & & & & & & & & & & $c_{n-2,1}$ & & \\
 & & & & & & & & & \begin{rotate}{45}$\geq$\end{rotate} & & \begin{rotate}{315}$\!\!\!\!\leq$\end{rotate} & \\
 & & & & & & & & $c_{n-3,1}$ & & & & $c_{n,1}$ \\
 & & & & & & & \begin{rotate}{45}$\geq$\end{rotate} & & \begin{rotate}{315}$\!\!\!\!\leq$\end{rotate} & & \begin{rotate}{45}$\geq$\end{rotate} & \\
 & & & & & & $c_{n-4,1}$ & & & & $c_{n-2,2}$ & & \\
 & & & & & \begin{rotate}{45}$\geq$\end{rotate} & & \begin{rotate}{315}$\!\!\!\!\leq$\end{rotate} & & \begin{rotate}{45}$\geq$\end{rotate} & & \begin{rotate}{315}$\!\!\!\!\leq$\end{rotate} & \\
 & & & & $c_{n-5,1}$ & & & & $c_{n-3,2}$ & & & & $c_{n-1,2}$ \\
 & & & \begin{rotate}{45}$\geq$\end{rotate} & & \begin{rotate}{315}$\!\!\!\!\leq$\end{rotate} & & \begin{rotate}{45}$\geq$\end{rotate} & & \begin{rotate}{315}$\!\!\!\!\leq$\end{rotate} & & \begin{rotate}{45}$\geq$\end{rotate} & \\
 & & $\ \iddots\ \ $ & & & & & $\cdot$ & $\cdot$ & $\cdot$ & & &  \\
 & \begin{rotate}{45}$\geq$\end{rotate} & & & & & & & & & & & \\
$\ c_{1,1}\ \ $ & & & & & & & $\cdot$ & $\cdot$ & $\cdot$ & & & \\
 & \begin{rotate}{315}$\!\!\!\!\leq$\end{rotate} & & & & & & & & & & \\
 & & $\ \ddots\ \ $ & & & & &  $\cdot$ & $\cdot$ & $\cdot$ & & & \\
 & & & \begin{rotate}{315}$\!\!\!\!\leq$\end{rotate} & & \begin{rotate}{45}$\geq$\end{rotate} & & \begin{rotate}{315}$\!\!\!\!\leq$\end{rotate} & & \begin{rotate}{45}$\geq$\end{rotate} & & \begin{rotate}{315}$\!\!\!\!\leq$\end{rotate} & \\
 & & & & $c_{n-5,n-5}$ & & & & $c_{n-3,n-4}$ & & & & $c_{n,\frac{n-1}{2}-1}$ \\
 & & & & & \begin{rotate}{315}$\!\!\!\!\leq$\end{rotate} & & \begin{rotate}{45}$\geq$\end{rotate} & & \begin{rotate}{315}$\!\!\!\!\leq$\end{rotate} & & \begin{rotate}{45}$\geq$\end{rotate} & \\
 & & & & & & $c_{n-4,n-4}$ & & & & $c_{n-2,n-3}$ & & \\
 & & & & & & & \begin{rotate}{315}$\!\!\!\!\leq$\end{rotate} & & \begin{rotate}{45}$\geq$\end{rotate} & & \begin{rotate}{315}$\!\!\!\!\leq$\end{rotate} & \\
 & & & & & & & & $c_{n-3,n-3}$ & & & & $c_{n-1,\frac{n-1}{2}}$ \\
 & & & & & & & & & \begin{rotate}{315}$\!\!\!\!\leq$\end{rotate} & & \begin{rotate}{45}$\geq$\end{rotate} & \\
 & & & & & & & & & & $c_{n-2,n-2}$ & & \\
 & & & & & & & & & & & \begin{rotate}{315}$\!\!\!\!\leq$\end{rotate} & \\
 & & & & & & & & & & & & $c_{n,\frac{n-1}{2}}$ 
\end{tabular}
\end{center}
\end{figure} 
\noindent
$n$ and bounded by $m$ and partially ordered as in the odd orthogonal case. 
The result is a diamond-colored distributive lattice. 
An alternative pattern of inequalities for even orthogonal ideal arrays replaces each $c_{n-1,k}$ with $c_{n,k}$ and vice-versa;  the resulting diamond-colored distributive lattice of such arrays is denoted $L_{\mytinyD_n}(m\omega_{n})$. 
Another way to realize $L_{\mytinyD_n}(m\omega_{n})$ is by recoloring the edges of $L_{\mytinyD_n}(m\omega_{n-1})$ by exchanging colors $n-1$ and $n$. 
In the $\myD_n$ case, we refer to $\omega_{n-1}$ and $\omega_{n}$ as the `spin-node fundamental weights', for obvious reasons. 

Our goal is to explicitly construct even and odd orthogonal Lie algebra representations whose highest weight is a multiple of a spin-node fundamental weight, cf.\ \NewOrthogonalConstructions/\NewOrthogonalCorollary\ below. 
From Theorem 8.10 in \S 8 of \cite{DonSupp}, we know that $L_{\mytinyB_n}(m\omega_{n})$, $L_{\mytinyD_n}(m\omega_{n-1})$, and $L_{\mytinyD_n}(m\omega_{n})$ are splitting distributive lattices for $\chi_{_{m\omega_{n}}}^{\mytinyB_{n}}$, $\chi_{_{m\omega_{n-1}}}^{\mytinyD_{n}}$, and $\chi_{_{m\omega_{n}}}^{\mytinyD_{n}}$ respectively. 
Alternatively, if, as in the second paragraph of this section, we view orthogonal ideal arrays as dual-multichains, then the claims of the preceding sentence can be deduced from the Seshadri result stated as Lemma 6 in \cite{PrEur}. 
In view of \RepDiagramTheorem, it therefore only remains to supply coefficients to the edges of these lattices and check the DC relations. 
In fact, we only need to supply edge coefficients for $L_{\mytinyD_n}(m\omega_{n-1})$, as the lattices $L_{\mytinyB_{n-1}}(m\omega_{n-1})$ and $L_{\mytinyD_n}(m\omega_{n})$ are simply re-colorings of $L_{\mytinyD_n}(m\omega_{n-1})$. 

For the moment, let $J := \{1,2,\ldots,n-1\}$, assume $n$ is even, and fix some ideal array $\xelt$ in $L_{\mytinyD_n}(m\omega_{n-1})$. 
In particular, all elements of $\comp_{J}(\xelt)$ have the same entries in positions $(n,1)$, $(n,2)$, $\ldots$, $(n,\frac{n}{2}-1)$, namely $c_{n,1}(\xelt), c_{n,2}(\xelt), \ldots, c_{n,\frac{n}{2}-1}(\xelt)$. 
Set $g_{n,1} := 0$, $g_{n,n} := m$, and $g_{n,q} := c_{n,\lfloor q/2 \rfloor}(\xelt)$ for $1 < q < n$, and let $\mysmallP := (g_{n,n}, \ldots, g_{n,1})$. 
Now we identify $\comp_{J}(\xelt)$ with the classical GT lattice $L_{\mytinyA_{n-1}}^{\mbox{\tiny GT}}(\mysmallP)$. 
For any $\telt \in \comp_{J}(\xelt)$, the corresponding element of said GT lattice has the following entries in its non-inert positions: $g_{p,q}(\telt) := c_{p,q}(\telt)$ for $1 \leq p \leq n-2$ and $g_{n-1,q}(\telt) := c_{n-(q\, \mbox{\tiny mod}\, 2),\lfloor (q+1)/2 \rfloor}(\telt)$. 
See \OrthogonalBijectionFig\ for an example when $n=6$. 
Suppose now that $\selt \myarrow{i} \telt$ is an edge in $\comp_{J}(\xelt)$, wherein $c_{i,j'}(\selt) + 1 = c_{i,j'}(\telt)$. 
If $i = n-1$, set $j := 2j'-1$, otherwise take $j := j'$. 
Thus $g_{i,j}(\selt)+1=g_{i,j}(\telt)$, while $g_{p,q}(\selt) = g_{p,q}(\telt)$ when $(p,q) \ne (i,j)$. 
With the foregoing notation, let $\myqP_{\selt,\telt}$ be the product of the positive and rational coefficients defined by $\left(\star\right)$ in \S \GTStatementsSection. 
Associate to the edge $\selt \myarrow{i} \telt$ in $L_{\mytinyD_n}(m\omega_{n-1})$ the two coefficients $\myqX^{\mbox{\tiny orth}}_{\telt,\selt} := \sqrt{\myqP_{\selt,\telt}} =: \myqY^{\mbox{\tiny orth}}_{\selt,\telt}$. 
Now assume $n$ is odd. 
All elements of $\comp_{J}(\xelt)$ share the same list of numbers $c_{n,1}(\xelt), c_{n,2}(\xelt), \ldots, c_{n,\frac{n-1}{2}}(\xelt)$; we take $g_{n,n} := c_{n,\lfloor (n-1)/2 \rfloor}(\xelt)$, but all other quantities are the same as in the `$n$ is even' case. 
With these modifications, we define the edge coefficients $\myqX^{\mbox{\tiny orth}}_{\telt,\selt}$ and $\myqY^{\mbox{\tiny orth}}_{\telt,\selt}$ as in the `$n$ is even' case. 
Thus, for all $n$, we have defined edge coefficients $\myqX^{\mbox{\tiny orth}}_{\telt,\selt}$ and $\myqY^{\mbox{\tiny orth}}_{\telt,\selt}$ on any edge $\selt \myarrow{i} \telt$ in $L_{\mytinyD_n}(m\omega_{n-1})$ whenever $i \in J$.

Next, set $J' := \{1,2,\ldots,n-2,n\}$, and assume $n$ is even. 
Note now that all elements of $\comp_{J'}(\xelt)$ share the same list of numbers $c_{n-1,1}(\xelt), c_{n-1,2}(\xelt), \ldots, c_{n-1,\frac{n}{2}}(\xelt)$. 
This time we set $g'_{n,1} := c_{n-1,1}(\xelt)$, $g'_{n,n} := c_{n-1,n/2}(\xelt)$, and $g'_{n,k} := c_{n-1,\lfloor (k+1)/2 \rfloor}(\xelt)$ for $1 < k < n$, and let $\mysmallP' := (g'_{n,n}, \ldots, g'_{n,1})$. 
For any $\telt \in \comp_{J'}(\xelt)$, the corresponding element of $L_{\mytinyA_{n-1}}^{\mbox{\tiny GT}}(\mysmallP')$ has the following entries in its non-inert positions: $g'_{p,q}(\telt) := c_{p,q}(\telt)$ for $1 \leq p \leq n-2$ and $g_{n-1,q}'(\telt) := c_{n-(q\, \mbox{\tiny mod}\, 2),\lfloor (q+1)/2 \rfloor}(\telt)$. 
Suppose now that $\selt \myarrow{i} \telt$ is an edge in $\comp_{J'}(\xelt)$, wherein $c_{i,j'}(\selt) + 1 = c_{i,j'}(\telt)$. 
If $i = n$, set $j := 2j'$, otherwise take $j := j'$. 
Thus $g'_{i,j}(\selt)+1=g'_{i,j}(\telt)$ with $g'_{p,q}(\selt) = g'_{p,q}(\telt)$ when $(p,q) \ne (i,j)$. 
In the formulas for the positive and rational coefficients defined by $\left(\star\right)$ in \S \GTStatementsSection, replace each $g_{r,s}$ with $g'_{r,s}$, and let $\myqP'_{\selt,\telt}$ be their product. 
Since  $g'_{i,q}(\telt) := c_{i,q}(\telt) =: g_{i,q}(\telt)$ $1 \leq i \leq n-2$, then $\myqP'_{\selt,\telt} = \myqP_{\selt,\telt}$, so each of $\myqX^{\mbox{\tiny orth}}_{\telt,\selt}$ and $\myqY^{\mbox{\tiny orth}}_{\selt,\telt}$ as defined in the previous paragraph can also be realized as $\sqrt{\myqP'_{\selt,\telt}}$. 
When $i=n$, we associate to the edge $\selt \myarrow{n} \telt$ in $L_{\mytinyD_n}(m\omega_{n-1})$ the two coefficients $\myqX^{\mbox{\tiny orth}}_{\telt,\selt} := \sqrt{\myqP'_{\selt,\telt}} =: \myqY^{\mbox{\tiny orth}}_{\selt,\telt}$. 
Now assume $n$ is odd. 
All elements of $\comp_{J'}(\xelt)$ share the same list of numbers $c_{n,1}(\xelt), c_{n,2}(\xelt), \ldots, c_{n,\frac{n-1}{2}}(\xelt)$; we take $g_{n,n} := m$, but all other quantities remain the same. 
With these modifications, we define the edge coefficients $\myqX^{\mbox{\tiny orth}}_{\telt,\selt}$ and $\myqY^{\mbox{\tiny orth}}_{\telt,\selt}$ on an edge $\selt \myarrow{n} \telt$ as in the `$n$ is even' case. 

\begin{figure}[tb]
\begin{center}
\OrthogonalBijectionFig: For the even orthogonal ideal array $\xelt \in L_{\mytinyD_6}(9\omega_{5})$ indicated below, we see a natural correspondence between elements of the $J$-component $\comp_{J}(\xelt)$ and elements of the skew-tabular lattice $L_{\mytinyA_{5}}^{\mbox{\tiny GT}}({\setlength{\unitlength}{0.125cm}
\begin{picture}(9,0)
\put(0,-0.4){\line(1,0){2}} \put(0,0.15){\line(1,0){2}} \put(0,0.7){\line(1,0){7}} \put(0,1.25){\line(1,0){9}} \put(0,1.8){\line(1,0){9}}
\put(0,-0.45){\line(0,1){2.3}} \put(1,-0.45){\line(0,1){2.3}} \put(2,-0.45){\line(0,1){2.3}} \put(3,0.65){\line(0,1){1.2}} \put(4,0.65){\line(0,1){1.2}} \put(5,0.65){\line(0,1){1.2}} \put(6,0.65){\line(0,1){1.2}} \put(7,0.65){\line(0,1){1.2}} \put(8,1.2){\line(0,1){0.65}} \put(9,1.2){\line(0,1){0.65}}
\end{picture}}
)$, 
where $J$ is the set of colors $\{1,2,3,4,5\}$.  

\setlength{\unitlength}{1cm}
\begin{picture}(15,7)
\put(0,2.5){$\xelt =$}
\put(0,-1.5){\begin{picture}(5,0)
\put(0.85,6){\vector(0,-1){1.5}} \put(0.55,6.3){\tiny color} \put(0.775,6.1){\tiny $1$}
\put(0.8,4){\scriptsize $4$}
\put(1.25,6.3){\vector(0,-1){1.5}}
\put(0.95,6.6){\tiny color}
\put(1.175,6.4){\tiny $2$}
\put(1.2,4.3){\scriptsize $3$}
\put(1.2,3.7){\scriptsize $5$}
\put(1.65,6.6){\vector(0,-1){1.5}}
\put(1.35,6.9){\tiny color}
\put(1.575,6.7){\tiny $3$}
\put(1.6,4.6){\scriptsize $2$}
\put(1.6,4){\scriptsize $5$}
\put(1.6,3.4){\scriptsize $6$}
\put(2.05,6.9){\vector(0,-1){1.5}}
\put(1.75,7.2){\tiny color}
\put(1.975,7){\tiny $4$}
\put(2,4.9){\scriptsize $1$}
\put(2,4.3){\scriptsize $4$}
\put(2,3.7){\scriptsize $5$}
\put(2,3.1){\scriptsize $7$}
\put(2.4,5.2){\scriptsize $1$}
\put(3.7,5.3){\vector(-1,0){1}} \put(3.8,5.35){\tiny color} \put(4.025,5.15){\tiny $5$}
\put(2.4,4.6){\scriptsize $2$}
\put(3.7,4.7){\vector(-1,0){1}} \put(3.8,4.75){\tiny color} \put(4.025,4.55){\tiny $6$}
\put(2.4,4){\scriptsize $4$}
\put(3.7,4.1){\vector(-1,0){1}} \put(3.8,4.15){\tiny color} \put(4.025,3.95){\tiny $5$}
\put(2.4,3.4){\scriptsize $7$}
\put(3.7,3.5){\vector(-1,0){1}} \put(3.8,3.55){\tiny color} \put(4.025,3.35){\tiny $6$}
\put(2.4,2.8){\scriptsize $8$}
\put(3.7,2.9){\vector(-1,0){1}} \put(3.8,2.95){\tiny color} \put(4.025,2.75){\tiny $5$}
\end{picture}}
\put(10,-1.5){\begin{picture}(5,0)
\put(0.2,6.9){\vector(1,-1){1}} 
\put(-1.8,7.5){\parbox{4cm}{\small \begin{center}Generic element\end{center}}}
\put(-1.8,7.1){\parbox{4cm}{\small \begin{center}of $L_{\mytinyA_{5}}^{\mbox{\tiny GT}}({\setlength{\unitlength}{0.125cm}
\begin{picture}(9,0)
\put(0,-0.4){\line(1,0){2}} \put(0,0.15){\line(1,0){2}} \put(0,0.7){\line(1,0){7}} \put(0,1.25){\line(1,0){9}} \put(0,1.8){\line(1,0){9}}
\put(0,-0.45){\line(0,1){2.3}} \put(1,-0.45){\line(0,1){2.3}} \put(2,-0.45){\line(0,1){2.3}} \put(3,0.65){\line(0,1){1.2}} \put(4,0.65){\line(0,1){1.2}} \put(5,0.65){\line(0,1){1.2}} \put(6,0.65){\line(0,1){1.2}} \put(7,0.65){\line(0,1){1.2}} \put(8,1.2){\line(0,1){0.65}} \put(9,1.2){\line(0,1){0.65}}
\end{picture}}
)$\end{center}}}
\put(0,4){$g_{1,1}$}
\put(0.6,4.5){$g_{2,1}$}
\put(0.6,3.5){$g_{2,2}$}
\put(1.2,5){$g_{3,1}$}
\put(1.2,4){$g_{3,2}$}
\put(1.2,3){$g_{3,3}$}
\put(1.8,5.5){$g_{4,1}$}
\put(1.8,4.5){$g_{4,2}$}
\put(1.8,3.5){$g_{4,3}$}
\put(1.8,2.5){$g_{4,4}$}
\put(2.4,6){$g_{5,1}$}
\put(2.6,5){\small $2$}
\put(2.4,4){$g_{5,3}$}
\put(2.6,3){\small $7$}
\put(2.4,2){$g_{5,5}$}
\put(3,6.5){{\small $0$} $=g_{6,1}$}
\put(3,5.5){{\small $2$} $=g_{6,2}$}
\put(3,4.5){{\small $2$} $=g_{6,3}$}
\put(3,3.5){{\small $7$} $=g_{6,4}$}
\put(3,2.5){{\small $7$} $=g_{6,5}$}
\put(3,1.5){{\small $9$} $=g_{6,6}$}
\end{picture}}
\put(6,-1.5){\begin{picture}(5,0)
\put(-0.2,6.4){\vector(1,-1){1}} 
\put(-2,7){\parbox{4cm}{\small \begin{center}Generic element\end{center}}}
\put(-2,6.6){\parbox{4cm}{\small \begin{center}of $\comp_{J}(\xelt)$\end{center}}}
\put(0,4){$c_{1,1}$}
\put(0.6,4.5){$c_{2,1}$}
\put(0.6,3.5){$c_{2,2}$}
\put(1.2,5){$c_{3,1}$}
\put(1.2,4){$c_{3,2}$}
\put(1.2,3){$c_{3,3}$}
\put(1.8,5.5){$c_{4,1}$}
\put(1.8,4.5){$c_{4,2}$}
\put(1.8,3.5){$c_{4,3}$}
\put(1.8,2.5){$c_{4,4}$}
\put(2.4,6){$c_{5,1}$}
\put(2.6,5){\small $2$}
\put(2.4,4){$c_{5,2}$}
\put(2.6,3){\small $7$}
\put(2.4,2){$c_{5,3}$}
\end{picture}}
\end{picture}
\end{center}
\end{figure}

Taking the two previous paragraphs together, we have defined, for all $n$, edge coefficients $\myqX^{\mbox{\tiny orth}}_{\telt,\selt}$ and $\myqY^{\mbox{\tiny orth}}_{\telt,\selt}$ on any edge $\selt \myarrow{i} \telt$ in $L_{\mytinyD_n}(m\omega_{n-1})$ whenever $i \in \{1,2,\dots,n\}$. 
The next result is an application of \RepDiagramTheorem\ and \GTBasisTheorem.

\noindent
{\bf \NewOrthogonalConstructions}\ \ {\sl Let} $L := L_{\mytinyD_{n}}(m\omega_{n-1})$, {\sl and assign the above scalars} $\{(\myqX^{\mbox{\tiny orth}}_{\telt,\selt},\myqY^{\mbox{\tiny orth}}_{\selt,\telt})\}_{\selt \rightarrow \telt \mbox{\scriptsize \ in } L}$ {\sl to the edges of $L$. 
Then $L$ is} $\myD_{n}${\sl -structured and the scalars satisfy all DC relations. 
Therefore the action of the generators of} $\mathfrak{g}(\myD_{n})$ {\sl on the vector space $V[L]$ as defined by the formulas} ($\blacklozenge$) {\sl in} \S \SetupSectionTwo\ {\sl is well-defined, $\{v_{\telt}\}_{\telt \in L}$ is a weight basis for the irreducible} $\mathfrak{g}(\myD_{n})${\sl -module $V[L]$ with highest weight $m\omega_{n-1}$, and the lattice $L$ together with the set of scalar pairs $\{(\myqX_{\telt,\selt},\myqY_{\selt,\telt})\}_{\selt \rightarrow \telt \mbox{\scriptsize \ in } L}$ is its representation diagram.}

{\em Proof.} That $L$ is $\myD_{n}$-structured is demonstrated in Proposition 8.7 of \S 8 of \cite{DonPosetModels}. 
Once we establish that $L$ is a representation diagram for some $\mathfrak{g}(\myD_{n})$-module as claimed, then we can identify $V[L]$ as irreducible with highest weight $m\omega_{n-1}$, since $L$ is a splitting distributive lattice for $\chi_{_{m\omega_{n-1}}}^{\mytinyD_{n}}$ by Theorem 8.10 of \cite{DonPosetModels}.  
As above, we take $J = \{1,2,\ldots,n-1\}$ and $J' = \{1,2,\ldots,n-2,n\}$. 
In the paragraphs preceding the theorem statement, we saw that each $J$-component of $L$ is a classical GT lattice. 
By definition, the coefficients assigned to each of these component edges have the same product (nonzero, by \ObservationalLemma.1) as the coefficients from the corresponding edge in the appropriate GT lattice. 
So, \RepDiagramTheorem.2 applies. 
It follows that all color $i$ crossing relations are satisfied when $i \in J$. 
This same reasoning shows that any diamond relation is satisfied if the edge colors for the diamond are within the set $J$. 
Similarly, we can use $J'$-components to establish all color $n$ crossing relations and to see that any diamond relation is satisfied if the diamond edge colors are within $J'$. 
It only remains to check diamond relations for diamonds with two edges of color $n-1$ and two of color $n$, as in \parbox{1.4cm}{\begin{center}
\setlength{\unitlength}{0.2cm}
\begin{picture}(5,3)
\put(2.5,0){\circle*{0.5}} \put(0.5,2){\circle*{0.5}}
\put(2.5,4){\circle*{0.5}} \put(4.5,2){\circle*{0.5}}
\put(0.5,2){\line(1,1){2}} \put(2.5,0){\line(-1,1){2}}
\put(4.5,2){\line(-1,1){2}} \put(2.5,0){\line(1,1){2}}
\put(1,0.7){\em \small n} \put(2.5,0.7){\em \small n-1}
\put(-0.3,2.55){\em \small n-1} \put(3.25,2.55){\em \small n}
\put(3,-0.75){\footnotesize $\relt$} \put(5.25,1.75){\footnotesize $\telt$}
\put(3,4){\footnotesize $\uelt$} \put(-1,1.75){\footnotesize $\selt$}
\end{picture} \end{center}}. 
A routine application of the definitions shows that $\myqP_{\selt,\uelt} \myqP'_{\telt,\uelt} = \myqP'_{\relt,\selt} \myqP_{\relt,\telt}$, which suffices to establish the two diamond relations in question. 
The claims within the last sentence of the theorem statement now follow from \RepDiagramTheorem.1.\hfill\QED

With respect to any given ordering of the weight basis constructed for $V[L]$ in \NewOrthogonalConstructions, the representing matrices for the $\mathfrak{g}(\myD_{n})$-generators $\myqx_{i}$ and $\myqy_{i}$ will be transposes of one another. 
However, this transpose property is obtained at the expense of not always having rational matrix entries. 
At this time, we have not yet obtained a set of exclusively rational solutions to the DC relations for $L_{\mytinyD_{n}}(m\omega_{n-1})$. 

\noindent
{\bf \NewOrthogonalCorollary}\ \ {\sl (1) For $n \geq 3$, view} $K := L_{\mytinyB_n}(m\omega_n)$ {\sl to be a re-coloring of} $L_{\mytinyD_{n+1}}(m\omega_{n}) := M$ {\sl where color $n+1$ edges in $M$ are given color $n$ in $K$.  Assign coefficients from the edges of $M$ to the corresponding edges of $K$. 
Then $K$ is} $\myB_{n}${\sl -structured and the edge coefficients satisfy all DC relations. 
Therefore the action of the generators of} $\mathfrak{g}(\myB_{n})$ {\sl on the vector space $V[K]$ as defined by the formulas} ($\blacklozenge$) {\sl in} \S \SetupSectionTwo\ {\sl is well-defined, $\{v_{\telt}\}_{\telt \in K}$ is a weight basis for the irreducible} $\mathfrak{g}(\myB_{n})${\sl -module $V[K]$ with highest weight $m\omega_{n}$, and the lattice $K$ together with the set of scalar pairs assigned to its edges is its representation diagram.}\\ 
{\sl (2) For $n \geq 4$, view} $L := L_{\mytinyD_n}(m\omega_n)$ {\sl to be a re-coloring of} $L_{\mytinyD_{n}}(m\omega_{n-1}) =: N$ {\sl wherein color $n$ (resp.\ color $n-1$) edges in $N$ are given color $n-1$ (resp.\ color $n$) in $L$.  Assign coefficients from the edges of $N$ to the corresponding edges of $L$. 
Then $L$ is} $\myD_{n}${\sl -structured and the edge coefficients satisfy all DC relations. 
Therefore the action of the generators of} $\mathfrak{g}(\myD_{n})$ {\sl on the vector space $V[L]$ as defined by the formulas} ($\blacklozenge$) {\sl in} \S \SetupSectionTwo\ {\sl is well-defined, $\{v_{\telt}\}_{\telt \in L}$ is a weight basis for the irreducible} $\mathfrak{g}(\myD_{n})${\sl -module $V[L]$ with highest weight $m\omega_{n}$, and the lattice $L$ together with the set of scalar pairs assigned to its edges is its representation diagram.} 

{\em Proof.} 
From Proposition 8.7 and Theorem 8.10 of \cite{DonPosetModels}, it follows that $K$ is $\myB_{n}$-structured, $\WGF(K) = \chi_{_{m\omega_{n}}}^{\mytinyB_{n}}$, $L$ is $\myD_{n}$-structured, and $\WGF(L) = \chi_{_{m\omega_{n}}}^{\mytinyD_{n}}$. 
We use induced actions from inclusions certain $\mathfrak{g}(\myB_{n}) \hookrightarrow \mathfrak{g}(\myD_{n+1})$ and $\mathfrak{g}(\myD_{n}) \hookrightarrow \mathfrak{g}(\myD_{n})$ to get remaining conclusions. 
In particular, let's say  $\{\myqx^{\mytinyX_{m}}_{i},\myqy^{\mytinyX_{m}}_{i}\}_{i=1}^{n}$ are generators for $\mathfrak{g}(\myX_{m})$. 
Consider the monomorphism $\mathfrak{g}(\myB_{n}) \longrightarrow \mathfrak{g}(\myD_{n+1})$ induced by $\myqx^{\mytinyB_{n}}_{i} \mapsto \myqx^{\mytinyD_{n+1}}_{i}$ and $\myqy^{\mytinyB_{n}}_{i} \mapsto \myqy^{\mytinyD_{n+1}}_{i}$ for $1 \leq i \leq n-1$ with $\myqx^{\mytinyB_{n}}_{n} \mapsto \myqx^{\mytinyD_{n+1}}_{n}+\myqx^{\mytinyD_{n+1}}_{n+1}$ and $\myqy^{\mytinyB_{n}}_{i} \mapsto \myqy^{\mytinyD_{n+1}}_{n}+\myqy^{\mytinyD_{n+1}}_{n+1}$. 
View $V[M]$ as a $\mathfrak{g}(\myB_{n})$-module under the action induced by the foregoing inclusion $\mathfrak{g}(\myB_{n}) \hookrightarrow \mathfrak{g}(\myD_{n+1})$. 
For the proof of part {\sl (1)}, it suffices to observe that $K$ is the representation diagram for the weight basis $\{v_{\xelt}\}_{\xelt \in M}$ of this $\mathfrak{g}(\myB_{n})$-module. 
Now consider the monomorphism $\mathfrak{h} := \mathfrak{g}(\myD_{n}) \longrightarrow \mathfrak{g}(\myD_{n}) =: \mathfrak{g}$ induced by $\myqx^{\mytinyD_{n}}_{i} \mapsto \myqx^{\mytinyD_{n}}_{i}$ and $\myqy^{\mytinyD_{n}}_{i} \mapsto \myqy^{\mytinyD_{n}}_{i}$ for $1 \leq i \leq n-2$ with $\myqx^{\mytinyD_{n}}_{n-1} \mapsto \myqx^{\mytinyD_{n}}_{n}$, $\myqy^{\mytinyD_{n}}_{n-1} \mapsto \myqy^{\mytinyD_{n}}_{n}$, $\myqx^{\mytinyD_{n}}_{n} \mapsto \myqx^{\mytinyD_{n}}_{n-1}$,  and $\myqy^{\mytinyD_{n}}_{n} \mapsto \myqy^{\mytinyD_{n}}_{n-1}$. 
View $V[N]$ as an $\mathfrak{h}$-module under the action induced by the foregoing inclusion $\mathfrak{h} \hookrightarrow \mathfrak{g}$. 
For the proof of part {\sl (2)}, it suffices to observe that $L$ is the representation diagram for the weight basis $\{v_{\xelt}\}_{\xelt \in N}$ of this $\mathfrak{h}$-module. 
\hfill\QED

Our constructions appear to be new. 
In particular, they are distinct from Molev's explicit constructions cf.\ \S 3 of \cite{MoSurvey}, as ours possess a restriction property that Molev's bases do not, and vice-versa. 
To be precise, let $\mathfrak{g} := \mathfrak{g}(\Phi)$ be the Lie algebra associated with some root system $\Phi$, and let $\mathfrak{g}_{1} \subset \mathfrak{g}_{2} \subset \cdots \mathfrak{g}_{k-1} \subset \mathfrak{g}_{k} = \mathfrak{g}$ be a sequence of Levi subalgebras of $\mathfrak{g}$ wherein $J_{k} := I$ and each $\mathfrak{g}_{i}$ resides within $\mathfrak{g}_{i+1}$ and is built from generators corresponding to a proper subset $J_{i} \subset J_{i+1}$. 
Following \S 3.3 of \cite{DonSupp}, a supporting graph $R$ or (its weight basis $\{v_{\xelt}\}_{\xelt \in R}$ for) a $\mathfrak{g}$-module $V[R]$ {\em restricts irreducibly for this chain of Levi subalgebras} if, for each $i \in \{1,\ldots,k-1\}$, each $J_{i}$-component of $R$ realizes an irreducible representation of $\mathfrak{g}_{i}$. 
A distinguishing feature of Molev's weight bases for the irreducible representations of $\mathfrak{g}(\myD_{n})$ is that they restrict irreducibly for the chain of Levi subalgebras $\mathfrak{g}(\myA_{1} \oplus \myA_{1}) \subset \mathfrak{g}(\myA_{3}) \subset \mathfrak{g}(\myD_{4}) \subset \mathfrak{g}(\myD_{5}) \subset \cdots \subset \mathfrak{g}(\myD_{n})$, where the inclusions within this chain are determined by the obvious inclusion of Dynkin diagrams.  
Now let $L$ be one of $L_{\mytinyD_{n}}(m\omega_{n-1})$ or $L_{\mytinyD_{n}}(m\omega_{n})$.  
It is not hard to see that $L$ does not always restrict irreducibly for the chain $\mathfrak{g}(\myA_{1} \oplus \myA_{1}) \subset \mathfrak{g}(\myA_{3}) \subset \mathfrak{g}(\myD_{4}) \subset \mathfrak{g}(\myD_{5}) \subset \cdots \subset \mathfrak{g}(\myD_{n})$. 
Similarly, Molev's weight bases for the irreducible $\mathfrak{g}(\myB_{n})$-modules restrict irreducibly for the chain of Levi subalgebras  $\mathfrak{g}(\myA_{1}) \subset \mathfrak{g}(\myB_{2}) \subset \mathfrak{g}(\myB_{3}) \subset \mathfrak{g}(\myB_{4}) \subset \cdots \subset \mathfrak{g}(\myB_{n})$, but our weight bases do not thusly restrict in general.  
The restriction properties our weight bases enjoy are presented in \OrthogonalSolitaryTheorem\ below. 

Rank symmetry and rank unimodality of the foregoing orthogonal lattices is well-known, cf.\ \cite{PrEur}. 
However, the Sperner aspect of the following result appears to be new.

\noindent
{\bf \NewOrthogonalCombinatorics}\ \ {\sl Let $L$ be any one of} $L_{\mytinyB_n}(m\omega_n)$, $L_{\mytinyD_{n+1}}(m\omega_{n})$, {\sl or} $L_{\mytinyD_{n+1}}(m\omega_{n+1})$.  {\sl Then $L$ is rank symmetric, rank unimodal, and strongly Sperner with} $\displaystyle \RGF(L) = \prod_{i=1}^{n}\prod_{j=i+1}^{n+1}\frac{[m+2n+2-i-j]_{q}}{[2n+2-i-j]_{q}}$. 

{\em Proof.}  In view of \NewOrthogonalConstructions/\NewOrthogonalCorollary, apply \MainCombinatorialTheorem. 
The expression given for  $\displaystyle \RGF(L)$ is a straightforward simplification of the formula in \MainCombinatorialTheorem; alternatively, use Theorem 10.6.3 of \cite{DonDiamond}.\hfill\QED

\vspace*{0.1in}
\noindent {\bf \S \SolitarySection\ Some extremal properties of the foregoing representation constructions.} 
We close this paper with a discussion of so-called `extremal properties' (cf.\ \cite{DonSupp}; particularly the solitary and edge-minimal properties summarized in \S \SetupSectionTwo\ above) enjoyed by the weight bases we have constructed above for special linear and orthogonal Lie algebra representations. 
In particular, we state conditions on a given skew shape $\mysmallP/\mysmallQ$ that are sufficient to guarantee that the associated skew-tabular lattice $L_{\mytinyA_{n-1}}^{\mbox{\tiny skew}}(\mysmallP/\mysmallQ)$ is solitary or edge-minimal. 
We also prove that each of $L_{\mytinyB_{n-1}}(m\omega_{n-1})$, $L_{\mytinyD_{n}}(m\omega_{n-1})$,  and $L_{\mytinyD_{n}}(m\omega_{n})$ is solitary and edge-minimal and restricts irreducibly for a distinctive chain of subalgebras. 

\begin{figure}[htb]
\begin{center}
\SolitarySigmaFig:  The skew-tabular lattice $L_{\mytinyA_{2}}^{\mbox{\tiny skew}}({\setlength{\unitlength}{0.125cm}\begin{picture}(3,0)\put(0,-0.25){\line(0,1){2}} \put(1,-0.25){\line(0,1){2}} \put(2,-0.25){\line(0,1){2}} \put(3,0.75){\line(0,1){1}} \put(0,-0.25){\line(1,0){2}} \put(0,0.75){\line(1,0){3}} \put(0,1.75){\line(1,0){3}}\end{picture}}) = L_{\mytinyA_{2}}^{\mbox{\tiny skew}}({\setlength{\unitlength}{0.125cm}\begin{picture}(3,0)\put(0,0){\line(0,1){1}} \put(1,0){\line(0,1){1}} \put(2,0){\line(0,1){2}} \put(3,0){\line(0,1){2}} \put(0,0){\line(1,0){3}} \put(0,1){\line(1,0){3}} \put(2,2){\line(1,0){1}}\end{picture}})^{\sigma_{0}}$.  

\setlength{\unitlength}{1.5cm}
\begin{picture}(4,6.5)
\put(1,0){\line(-1,1){1}}
\put(1,0){\line(1,1){2}}
\put(0,1){\line(1,1){2}}
\put(2,1){\line(-1,1){1}}
\put(2,1){\line(0,1){1}}
\put(1,2){\line(0,1){1}}
\put(2,2){\line(-1,1){2}}
\put(2,2){\line(1,1){2}}
\put(3,2){\line(-1,1){1}}
\put(3,2){\line(0,1){1}}
\put(1,3){\line(1,1){2}}
\put(2,3){\line(0,1){1}}
\put(3,3){\line(-1,1){2}}
\put(0,4){\line(1,1){2}}
\put(4,4){\line(-1,1){2}}
\put(2,6){\VertexTableauTwo{1}{1}{1}{2}{2}{-0.65}{-0.05}}
\put(1,5){\VertexTableauTwo{1}{1}{1}{2}{3}{-0.65}{-0.05}}
\put(3,5){\VertexTableauTwo{1}{1}{2}{2}{2}{0.15}{-0.05}}
\put(0,4){\VertexTableauTwo{1}{1}{1}{3}{3}{-0.75}{-0.05}}
\put(2,4){\VertexTableauTwo{1}{1}{2}{2}{3}{0.35}{-0.15}}
\put(4,4){\VertexTableauTwo{1}{1}{3}{2}{2}{0.15}{-0.05}}
\put(1,3){\VertexTableauTwo{1}{1}{2}{3}{3}{-0.8}{-0.3}}
\put(2,3){\VertexTableauTwo{1}{2}{2}{2}{3}{-0.75}{-0.2}}
\put(3,3){\VertexTableauTwo{1}{1}{3}{2}{3}{0.15}{-0.3}}
\put(1,2){\VertexTableauTwo{1}{2}{2}{3}{3}{-0.8}{-0.2}}
\put(2,2){\VertexTableauTwo{1}{1}{3}{3}{3}{0.15}{-0.3}}
\put(3,2){\VertexTableauTwo{1}{2}{3}{2}{3}{0.15}{-0.2}}
\put(0,1){\VertexTableauTwo{2}{2}{2}{3}{3}{-0.75}{-0.25}}
\put(2,1){\VertexTableauTwo{1}{2}{3}{3}{3}{0.25}{-0.25}}
\put(1,0){\VertexTableauTwo{2}{2}{3}{3}{3}{0.25}{-0.25}}
\put(1,5){\NEEdgeLabelForLatticeI{{\em 2}}}
\put(3,5){\NWEdgeLabelForLatticeI{{\em 1}}}
\put(0,4){\NEEdgeLabelForLatticeI{{\em 2}}}
\put(2,4){\NWEdgeLabelForLatticeI{{\em 1}}}
\put(2,4){\NEEdgeLabelForLatticeI{{\em 2}}}
\put(4,4){\NWEdgeLabelForLatticeI{{\em 2}}}
\put(1,3){\NEEdgeLabelForLatticeI{{\em 2}}}
\put(1,3){\NWEdgeLabelForLatticeI{{\em 1}}}
\put(2,3){\VerticalEdgeLabelForLatticeI{{\em 1}}}
\put(3,3){\NWEdgeLabelForLatticeI{{\em 2}}}
\put(3,3){\NEEdgeLabelForLatticeI{{\em 2}}}
\put(1,2){\VerticalEdgeLabelForLatticeI{{\em 1}}}
\put(1.25,2.25){\NEEdgeLabelForLatticeI{{\em 2}}}
\put(2.2,1.8){\NWEdgeLabelForLatticeI{{\em 2}}}
\put(1.8,1.8){\NEEdgeLabelForLatticeI{{\em 2}}}
\put(3,2){\VerticalEdgeLabelForLatticeI{{\em 1}}}
\put(2.75,2.25){\NWEdgeLabelForLatticeI{{\em 2}}}
\put(0,1){\NEEdgeLabelForLatticeI{{\em 1}}}
\put(2,1){\VerticalEdgeLabelForLatticeI{{\em 1}}}
\put(2,1){\NWEdgeLabelForLatticeI{{\em 2}}}
\put(2,1){\NEEdgeLabelForLatticeI{{\em 2}}}
\put(1,0){\NWEdgeLabelForLatticeI{{\em 2}}}
\put(1,0){\NEEdgeLabelForLatticeI{{\em 1}}}
\end{picture}
\end{center}
\end{figure}

For skew-tabular lattices, we begin by defining some operations on skew shapes. 
Our focal skew shape is $\mysmallP/\mysmallQ$, where $\mysmallP$ and $\mysmallQ$ are $m$-tuples of nonnegative integers and are skew-compatible with respect to $n$ for some fixed integers $m$ and $n$ satisfying $m \geq n \geq 2$. 
To form $(\mysmallP/\mysmallQ)^{\sigma_{0}}$, add cells directly below each column of $\mysmallP/\mysmallQ$ so that each column now has exactly $n$ cells, and then remove $\mysmallP/\mysmallQ$. 
(The operation $(\mysmallP/\mysmallQ)^{\sigma_{0}}$ is so denoted because of a connection we will make with the involution $\sigma_{0}$ on the set of colors $I = \{1,2,\ldots,n-1\}$ induced by the longest element $w_{0}$ of the Weyl group $W(\myA_{n-1})$.) 
To form $(\mysmallP/\mysmallQ)^{*}$, rotate $(\mysmallP/\mysmallQ)^{\sigma_{0}}$ by $180^{\circ}$. 
To form $(\mysmallP/\mysmallQ)^{\bowtie}$, rotate $\mysmallP/\mysmallQ$ by $180^{\circ}$. 
It is useful to set $(\mysmallP/\mysmallQ)^{\varepsilon} := \mysmallP/\mysmallQ$.
Notice that we have a natural action of the Klein four-group $\mathbb{Z}_{2} \times \mathbb{Z}_{2} = \{(0,0), (1,0), (0,1), (1,1)\}$ on such skew shapes via the identifications 
\begin{center}
\begin{tabular}{ccccccc}
$(0,0).\mysmallP/\mysmallQ = (\mysmallP/\mysmallQ)^{\varepsilon}$ & \hspace*{0.05in} & $(1,0).\mysmallP/\mysmallQ = (\mysmallP/\mysmallQ)^{\sigma_{0}}$ & \hspace*{0.05in} & $(0,1).\mysmallP/\mysmallQ = (\mysmallP/\mysmallQ)^{*}$ & \hspace*{0.05in} & $(1,1).\mysmallP/\mysmallQ = (\mysmallP/\mysmallQ)^{\bowtie}$
\end{tabular}
\end{center}

Some concepts from the general setting are apropos here; these are mostly standard. 
Let $R$ be a ranked poset with edges colored by our simple-root-indexing set $I$ associated with a given root system $\Phi$, and assume there are two coefficients, denoted $\myqX_{\telt,\selt}$ and $\myqY_{\selt,\telt}$, attached to each edge $\selt \myarrow{i} \telt$; so $R$ is an edge-tagged poset. 
As in \cite{ADLPOne} and \cite{ADLMPPW}, we define edge-tagged posets $R^{\sigma_{0}}$, $R^{*}$, and $R^{\bowtie}$ as follows\footnote{The notation `$R^{\triangle}$' from \cite{ADLPOne} and \cite{ADLMPPW} is replaced here by `$R^{\bowtie}$', which we call the `bowtie' of $R$.}. 
Declare $R^{\sigma_{0}}$ to be the set $R$ whose edges are of the form $\selt \mylongarrow{\mbox{\tiny $\sigma_{0}(i)$}} \telt$, whenever $\selt \myarrow{i} \telt$ is an edge in $R$, and are assigned the corresponding edge coefficients from $R$. 
As a set, we have $R^{*} := \{\telt^{*}\}_{\telt \in R}$, and we have an edge $\telt^{*} \myarrow{i} \selt^{*}$ in $R^{*}$ with coefficients $\myqX_{\telt^{*},\selt^{*}}$ and $\myqY_{\selt^{*},\telt^{*}}$ if and only if $\selt \myarrow{i} \telt$ is an edge in $R$ whose coefficients are $\myqX_{\telt,\selt} = \myqY_{\selt^{*},\telt^{*}}$ and $\myqY_{\selt,\telt} = \myqX_{\telt^{*},\selt^{*}}$.  
Let $R^{\bowtie}$ be the edge-tagged poset $(R^{\sigma_{0}})^{*}$.  
For completeness, let $R^{\varepsilon}$ be the edge-tagged poset $R$. 
The notations $^{\varepsilon}$, $^{*}$, $^{\sigma_{0}}$, and $^{\bowtie}$ are understood to apply to ranked posets with edges colored by $I$ but with no edge coefficients.

\begin{figure}[htb]
\begin{center}
\SolitaryDualFig:  The skew-tabular lattice $L_{\mytinyA_{2}}^{\mbox{\tiny skew}}({\setlength{\unitlength}{0.125cm}\begin{picture}(3,0)\put(0,-0.25){\line(0,1){1}} \put(1,-0.25){\line(0,1){2}} \put(2,-0.25){\line(0,1){2}} \put(3,-0.25){\line(0,1){2}} \put(0,-0.25){\line(1,0){3}} \put(0,0.75){\line(1,0){3}} \put(1,1.75){\line(1,0){2}}\end{picture}}) = L_{\mytinyA_{2}}^{\mbox{\tiny skew}}({\setlength{\unitlength}{0.125cm}\begin{picture}(3,0)\put(0,0){\line(0,1){1}} \put(1,0){\line(0,1){1}} \put(2,0){\line(0,1){2}} \put(3,0){\line(0,1){2}} \put(0,0){\line(1,0){3}} \put(0,1){\line(1,0){3}} \put(2,2){\line(1,0){1}}\end{picture}})^{*}$.  

\setlength{\unitlength}{1.5cm}
\begin{picture}(4,6.5)
\put(2,0){\line(-1,1){2}}
\put(2,0){\line(1,1){2}}
\put(1,1){\line(1,1){2}}
\put(3,1){\line(-1,1){2}}
\put(0,2){\line(1,1){2}}
\put(4,2){\line(-1,1){2}}
\put(2,2){\line(0,1){1}}
\put(1,3){\line(0,1){1}}
\put(3,3){\line(0,1){1}}
\put(2,4){\line(0,1){1}}
\put(2,3){\line(-1,1){2}}
\put(3,4){\line(-1,1){2}}
\put(1,4){\line(1,1){1}}
\put(0,5){\line(1,1){1}}
\put(2,3){\line(1,1){1}}
\put(1,6){\VertexTableauThree{1}{1}{2}{1}{2}{-0.75}{-0.05}}
\put(0,5){\VertexTableauThree{2}{1}{2}{1}{2}{-0.75}{-0.15}}
\put(2,5){\VertexTableauThree{1}{1}{2}{1}{3}{0.2}{-0.15}}
\put(1,4){\VertexTableauThree{2}{1}{2}{1}{3}{-0.75}{-0.3}}
\put(2,4){\VertexTableauThree{1}{1}{3}{1}{3}{-0.7}{-0.05}}
\put(3,4){\VertexTableauThree{1}{1}{2}{2}{3}{0.1}{-0.15}}
\put(1,3){\VertexTableauThree{2}{1}{3}{1}{3}{-0.75}{-0.05}}
\put(2,3){\VertexTableauThree{2}{1}{2}{2}{3}{0.1}{-0.2}}
\put(3,3){\VertexTableauThree{1}{1}{3}{2}{3}{0.1}{-0.05}}
\put(0,2){\VertexTableauThree{3}{1}{3}{1}{3}{-0.75}{-0.15}}
\put(2,2){\VertexTableauThree{2}{1}{3}{2}{3}{0.2}{-0.15}}
\put(4,2){\VertexTableauThree{1}{2}{3}{2}{3}{0.1}{-0.15}}
\put(3,1){\VertexTableauThree{2}{2}{3}{2}{3}{0.1}{-0.25}}
\put(1,1){\VertexTableauThree{3}{1}{3}{2}{3}{-0.80}{-0.25}}
\put(2,0){\VertexTableauThree{3}{2}{3}{2}{3}{0.1}{-0.25}}
\put(2,2){\VerticalEdgeLabelForLatticeI{{\em 2}}}
\put(1,3){\VerticalEdgeLabelForLatticeI{{\em 2}}}
\put(3,3){\VerticalEdgeLabelForLatticeI{{\em 2}}}
\put(2,4){\VerticalEdgeLabelForLatticeI{{\em 2}}}
\put(0,5){\NEEdgeLabelForLatticeI{{\em 1}}}
\put(1,4){\NEEdgeLabelForLatticeI{{\em 1}}}
\put(3,1){\NEEdgeLabelForLatticeI{{\em 1}}}
\put(2,2){\NEEdgeLabelForLatticeI{{\em 1}}}
\put(1.25,3.25){\NEEdgeLabelForLatticeI{{\em 1}}}
\put(1.75,2.75){\NEEdgeLabelForLatticeI{{\em 1}}}
\put(2,0){\NEEdgeLabelForLatticeI{{\em 2}}}
\put(1,1){\NEEdgeLabelForLatticeI{{\em 2}}}
\put(0,2){\NEEdgeLabelForLatticeI{{\em 2}}}
\put(1,1){\NWEdgeLabelForLatticeI{{\em 1}}}
\put(3,4){\NWEdgeLabelForLatticeI{{\em 1}}}
\put(2,2){\NWEdgeLabelForLatticeI{{\em 1}}}
\put(2.75,3.25){\NWEdgeLabelForLatticeI{{\em 1}}}
\put(2.25,2.75){\NWEdgeLabelForLatticeI{{\em 1}}}
\put(2,0){\NWEdgeLabelForLatticeI{{\em 1}}}
\put(3,1){\NWEdgeLabelForLatticeI{{\em 1}}}
\put(4,2){\NWEdgeLabelForLatticeI{{\em 1}}}
\put(1,4){\NWEdgeLabelForLatticeI{{\em 2}}}
\put(2,5){\NWEdgeLabelForLatticeI{{\em 2}}}
\end{picture}
\end{center}
\end{figure}

\begin{figure}[htb]
\begin{center}
\SolitaryBowtieFig:  The skew-tabular lattice $L_{\mytinyA_{2}}^{\mbox{\tiny skew}}({\setlength{\unitlength}{0.125cm}\begin{picture}(3,0)\put(0,-0.25){\line(0,1){2}} \put(1,-0.25){\line(0,1){2}} \put(2,0.75){\line(0,1){1}} \put(3,0.75){\line(0,1){1}} \put(0,-0.25){\line(1,0){1}} \put(0,0.75){\line(1,0){3}} \put(0,1.75){\line(1,0){3}}\end{picture}}) = L_{\mytinyA_{2}}^{\mbox{\tiny skew}}({\setlength{\unitlength}{0.125cm}\begin{picture}(3,0)\put(0,0){\line(0,1){1}} \put(1,0){\line(0,1){1}} \put(2,0){\line(0,1){2}} \put(3,0){\line(0,1){2}} \put(0,0){\line(1,0){3}} \put(0,1){\line(1,0){3}} \put(2,2){\line(1,0){1}}\end{picture}})^{\bowtie}$.  

\setlength{\unitlength}{1.5cm}
\begin{picture}(4,6.5)
\put(2,0){\line(-1,1){2}}
\put(2,0){\line(1,1){2}}
\put(1,1){\line(1,1){2}}
\put(3,1){\line(-1,1){2}}
\put(0,2){\line(1,1){2}}
\put(4,2){\line(-1,1){2}}
\put(2,2){\line(0,1){1}}
\put(1,3){\line(0,1){1}}
\put(3,3){\line(0,1){1}}
\put(2,4){\line(0,1){1}}
\put(2,3){\line(-1,1){2}}
\put(3,4){\line(-1,1){2}}
\put(1,4){\line(1,1){1}}
\put(0,5){\line(1,1){1}}
\put(2,3){\line(1,1){1}}
\put(1,6){\VertexTableauFour{1}{2}{1}{1}{-0.65}{-0.05}}
\put(0,5){\VertexTableauFour{1}{3}{1}{1}{-0.7}{-0.15}}
\put(2,5){\VertexTableauFour{1}{2}{1}{2}{0.2}{-0.15}}
\put(1,4){\VertexTableauFour{1}{3}{1}{2}{-0.75}{-0.3}}
\put(2,4){\VertexTableauFour{1}{2}{2}{2}{-0.7}{-0.2}}
\put(3,4){\VertexTableauFour{1}{2}{1}{3}{0.1}{-0.15}}
\put(1,3){\VertexTableauFour{1}{3}{2}{2}{-0.75}{-0.10}}
\put(2,3){\VertexTableauFour{1}{3}{1}{3}{0.1}{-0.4}}
\put(3,3){\VertexTableauFour{1}{2}{2}{3}{0.1}{-0.05}}
\put(0,2){\VertexTableauFour{2}{3}{2}{2}{-0.7}{-0.15}}
\put(2,2){\VertexTableauFour{1}{3}{2}{3}{0.25}{-0.2}}
\put(4,2){\VertexTableauFour{1}{2}{3}{3}{0.1}{-0.15}}
\put(3,1){\VertexTableauFour{1}{3}{3}{3}{0.2}{-0.25}}
\put(1,1){\VertexTableauFour{2}{3}{2}{3}{-0.80}{-0.25}}
\put(2,0){\VertexTableauFour{2}{3}{3}{3}{0.25}{-0.25}}
\put(2,2){\VerticalEdgeLabelForLatticeI{{\em 1}}}
\put(1,3){\VerticalEdgeLabelForLatticeI{{\em 1}}}
\put(3,3){\VerticalEdgeLabelForLatticeI{{\em 1}}}
\put(2,4){\VerticalEdgeLabelForLatticeI{{\em 1}}}
\put(0,5){\NEEdgeLabelForLatticeI{{\em 2}}}
\put(1,4){\NEEdgeLabelForLatticeI{{\em 2}}}
\put(3,1){\NEEdgeLabelForLatticeI{{\em 2}}}
\put(2,2){\NEEdgeLabelForLatticeI{{\em 2}}}
\put(1.25,3.25){\NEEdgeLabelForLatticeI{{\em 2}}}
\put(1.75,2.75){\NEEdgeLabelForLatticeI{{\em 2}}}
\put(2,0){\NEEdgeLabelForLatticeI{{\em 1}}}
\put(1,1){\NEEdgeLabelForLatticeI{{\em 1}}}
\put(0,2){\NEEdgeLabelForLatticeI{{\em 1}}}
\put(1,1){\NWEdgeLabelForLatticeI{{\em 2}}}
\put(3,4){\NWEdgeLabelForLatticeI{{\em 2}}}
\put(2,2){\NWEdgeLabelForLatticeI{{\em 2}}}
\put(2.75,3.25){\NWEdgeLabelForLatticeI{{\em 2}}}
\put(2.25,2.75){\NWEdgeLabelForLatticeI{{\em 2}}}
\put(2,0){\NWEdgeLabelForLatticeI{{\em 2}}}
\put(3,1){\NWEdgeLabelForLatticeI{{\em 2}}}
\put(4,2){\NWEdgeLabelForLatticeI{{\em 2}}}
\put(1,4){\NWEdgeLabelForLatticeI{{\em 1}}}
\put(2,5){\NWEdgeLabelForLatticeI{{\em 1}}}
\end{picture}
\end{center}
\end{figure}

We can naturally define related objects $V^{\sigma_{0}}$, $V^{*}$, and $V^{\bowtie}$ for any $\mathfrak{g}(\Phi)$-module $V =: V^{\varepsilon}$. 
The dual (or contragredient) of $V$ is $V^{*}$.  
The $\mathfrak{g}(\Phi)$-module $V^{\sigma_{0}}$ is defined by the homomorphism $\eta: \mathfrak{g}(\Phi) \longrightarrow \mathfrak{gl}(V)$ wherein $\eta(\myqx_{i})(v) := \myqx_{\sigma_{0}(i)}.v$ and $\eta(\myqy_{i})(v) := \myqy_{\sigma_{0}(i)}.v$. 
We set $V^{\bowtie} := (V^{*})^{\sigma_{0}}$. 
One can see that $V^{\sigma_{0}} \cong V^{*}$ and $V^{\bowtie} \cong V^{\varepsilon} \cong (V^{\sigma_{0}})^{*}$. 
There are naturally related actions on the group ring $\mathbb{Z}[\Lambda]$ induced by $(\myvarZ^{\mu})^{*} := \myvarZ^{-\mu}$, $(\myvarZ^{\mu})^{\sigma_{0}} := \myvarZ^{-w_{0}.\mu}$, and $(\myvarZ^{\mu})^{\bowtie} := ((\myvarZ^{\mu})^{*})^{\sigma_{0}} = \myvarZ^{w_{0}.\mu}$; observe that these actions preserve the ring of $W(\Phi)$-symmetric functions. 
Moreover, as in Proposition 2.16 of \cite{DonPosetModels}, one can see that $\chi^{\sigma_{0}} = \chi^{*}$ and $\chi^{\bowtie} = \chi^{\varepsilon} = (\chi^{\sigma_{0}})^{*}$. 

\noindent 
{\bf \StarProposition}\ \ {\sl  
(1) Let $R$ be a representation diagram for a $\mathfrak{g}(\Phi)$-module $V$, so $V[R] \cong V$. Then $V[R^{*}] \cong V[R]^{*} \cong V[R]^{\sigma_{0}} \cong V[R^{\sigma_{0}}]$ and $V[R^{\bowtie}] \cong V[R]^{\bowtie} \cong V[R]^{\varepsilon} \cong V[R^{\varepsilon}]$. 
(2) Let $R$ be a splitting poset for a $W(\Phi)$-symmetric function $\chi$, so} $\WGF(R) = \chi$.  {\sl Then} $\WGF(R^{*}) = \WGF(R)^{*} = \WGF(R)^{\sigma_{0}} = \WGF(R^{\sigma_{0}})$ {\sl and} $\WGF(R^{\bowtie}) = \WGF(R)^{\bowtie} = \WGF(R)^{\varepsilon} = \WGF(R^{\varepsilon})$.

{\em Proof.} Part {\sl (1)} follows from Lemmas 2.3/2.4 and Proposition 2.5 of \cite{ADLPOne}.  Part {\sl (2)} follows from Lemmas 2.1/2.2 of \cite{ADLMPPW}.\hfill\QED

Returning now to the specific setting of skew-tabular lattices, the edge-color-preserving poset isomorphisms of the next proposition are readily apparent. 
See \SolitaryFigs\ for illustrations without edge coefficients.   

\noindent
{\bf \StarIsomorphisms}\ \ {\sl Let $\star \in \{\varepsilon,\sigma_{0},*,\bowtie\}$.  Then} $L_{\mytinyA_{n-1}}^{\mbox{\tiny skew}}((\mysmallP/\mysmallQ)^{\star}) \cong L_{\mytinyA_{n-1}}^{\mbox{\tiny skew}}(\mysmallP/\mysmallQ)^{\star}$.

\noindent
{\bf \GTSolitaryMinimal}\ \ {\sl Let $V$ be the} $\mathfrak{g}(\myA_{n-1})${\sl -module $V[L]$, where} $L := L_{\mytinyA_{n-1}}^{\mbox{\tiny skew}}(\mysmallP/\mysmallQ)$. 
{\sl (0) For each $\star \in \{\varepsilon,\sigma_{0},*,\bowtie\}$, $V[L^{\star}] \cong V^{\star}$ and} $\WGF(L^{\star}) = \WGF(L)^{\star}$.
{\sl (1) If, for some $\star \in \{\varepsilon,\sigma_{0},*,\bowtie\}$,} $L_{\mytinyA_{n-1}}^{\mbox{\tiny skew}}((\mysmallP/\mysmallQ)^{\star})$ {\sl is a solitary (respectively, edge-minimal) supporting graph for $V^{\star}$, then} $L_{\mytinyA_{n-1}}^{\mbox{\tiny skew}}((\mysmallP/\mysmallQ)^{\lozenge})$ {\sl is a solitary (respectively, edge-minimal) supporting graph for $V^{\lozenge}$ for all $\lozenge \in \{\varepsilon,\sigma_{0},*,\bowtie\}$.} 
{\sl (2) If} $(\mysmallP/\mysmallQ)^{\star}$ {\sl is non-skew for some $\star \in \{\varepsilon,\sigma_{0},*,\bowtie\}$, then for all $\lozenge \in \{\varepsilon,\sigma_{0},*,\bowtie\}$,} $L_{\mytinyA_{n-1}}^{\mbox{\tiny skew}}((\mysmallP/\mysmallQ)^{\lozenge})$ {\sl is a solitary and edge-minimal supporting graph for the irreducible} $\mathfrak{g}(\myA_{n-1})${\sl -module $V^{\lozenge}$.}

{\em Proof.} Part {\sl (0)} follows from \StarProps. 
For {\sl (1)}, see Lemmas 2.3/2.4 and Proposition 2.5 of \cite{ADLPOne}.   
Part {\sl (2)} follows from Theorem 4.4 of \cite{DonSupp} or Theorem 5.5 of \cite{HL}.\hfill\QED

Note, however, that there exist solitary and edge-minimal skew-tabular lattices that fall outside the purview of \GTSolitaryMinimal.2. 
The skew-tabular lattice $L_{\mytinyA_{2}}^{\mbox{\tiny skew}}({\setlength{\unitlength}{0.125cm}\begin{picture}(3,0)\put(0,-0.75){\line(0,1){2}} \put(1,-0.75){\line(0,1){2}} \put(2,0.25){\line(0,1){2}} \put(3,0.25){\line(0,1){2}} \put(0,-0.75){\line(1,0){1}} \put(0,0.25){\line(1,0){3}} \put(0,1.25){\line(1,0){3}} \put(2,2.25){\line(1,0){1}}\end{picture}}) = L_{\mytinyA_{2}}^{\mbox{\tiny skew}}({\setlength{\unitlength}{0.125cm}\begin{picture}(3,0)\put(0,-0.75){\line(0,1){2}} \put(1,-0.75){\line(0,1){2}} \put(2,0.25){\line(0,1){2}} \put(3,0.25){\line(0,1){2}} \put(0,-0.75){\line(1,0){1}} \put(0,0.25){\line(1,0){3}} \put(0,1.25){\line(1,0){3}} \put(2,2.25){\line(1,0){1}}\end{picture}})^{\bowtie}$ of \IntroFigTwo\ as well as $L_{\mytinyA_{2}}^{\mbox{\tiny skew}}({\setlength{\unitlength}{0.125cm}\begin{picture}(3,0)\put(0,-0.25){\line(0,1){1}} \put(1,-0.25){\line(0,1){2}} \put(2,-0.25){\line(0,1){2}} \put(3,0.75){\line(0,1){1}} \put(0,-0.25){\line(1,0){2}} \put(0,0.75){\line(1,0){3}} \put(1,1.75){\line(1,0){2}}\end{picture}}) = L_{\mytinyA_{2}}^{\mbox{\tiny skew}}({\setlength{\unitlength}{0.125cm}\begin{picture}(3,0)\put(0,-0.75){\line(0,1){2}} \put(1,-0.75){\line(0,1){2}} \put(2,0.25){\line(0,1){2}} \put(3,0.25){\line(0,1){2}} \put(0,-0.75){\line(1,0){1}} \put(0,0.25){\line(1,0){3}} \put(0,1.25){\line(1,0){3}} \put(2,2.25){\line(1,0){1}}\end{picture}})^{\sigma_{0}} = L_{\mytinyA_{2}}^{\mbox{\tiny skew}}({\setlength{\unitlength}{0.125cm}\begin{picture}(3,0)\put(0,-0.75){\line(0,1){2}} \put(1,-0.75){\line(0,1){2}} \put(2,0.25){\line(0,1){2}} \put(3,0.25){\line(0,1){2}} \put(0,-0.75){\line(1,0){1}} \put(0,0.25){\line(1,0){3}} \put(0,1.25){\line(1,0){3}} \put(2,2.25){\line(1,0){1}}\end{picture}})^{*}$ are examples. 
Moreover, there exist skew-tabular lattice supporting graphs that are neither solitary nor edge-minimal, such as $L_{\mytinyA_{2}}^{\mbox{\tiny skew}}({\setlength{\unitlength}{0.125cm}\begin{picture}(2,0)\put(0,-0.75){\line(0,1){2}} \put(1,-0.75){\line(0,1){3}} \put(2,0.25){\line(0,1){2}} \put(0,-0.75){\line(1,0){1}} \put(0,0.25){\line(1,0){2}} \put(0,1.25){\line(1,0){2}} \put(1,2.25){\line(1,0){1}}\end{picture}}) = L_{\mytinyA_{2}}^{\mbox{\tiny skew}}({\setlength{\unitlength}{0.125cm}\begin{picture}(2,0)\put(0,-0.25){\line(0,1){1}} \put(1,-0.25){\line(0,1){2}} \put(2,0.75){\line(0,1){1}} \put(0,-0.25){\line(1,0){1}} \put(0,0.75){\line(1,0){2}} \put(1,1.75){\line(1,0){1}}\end{picture}})^{*}$. 
We leave verification of the claims in the preceding two sentences to the reader as pleasant exercises. 

We can make more definitive claims about our orthogonal Lie algebra representation constructions. 
We observe for the record the following congruences of edge-tagged posets: $L_{\mytinyB_n}(m\omega_n) \cong L_{\mytinyB_n}(m\omega_n)^{*} \cong L_{\mytinyB_n}(m\omega_n)^{\sigma_{0}} \cong L_{\mytinyB_n}(m\omega_n)^{\bowtie}$ for all $n \geq 2$; we have $L_{\mytinyD_n}(m\omega_{n-1}) \cong L_{\mytinyD_n}(m\omega_{n-1})^{*} \cong L_{\mytinyD_n}(m\omega_{n})^{\sigma_{0}} \cong  L_{\mytinyD_n}(m\omega_{n})^{\bowtie}$ and $L_{\mytinyD_n}(m\omega_{n}) \cong L_{\mytinyD_n}(m\omega_{n})^{*} \cong L_{\mytinyD_n}(m\omega_{n-1})^{\sigma_{0}} \cong L_{\mytinyD_n}(m\omega_{n-1})^{\bowtie}$ when $n \geq 4$ is even; and we have $L_{\mytinyD_n}(m\omega_{n-1}) \cong L_{\mytinyD_n}(m\omega_{n-1})^{\sigma_{0}} \cong L_{\mytinyD_n}(m\omega_{n})^{*}  \cong  L_{\mytinyD_n}(m\omega_{n})^{\bowtie}$ and $L_{\mytinyD_n}(m\omega_{n}) \cong L_{\mytinyD_n}(m\omega_{n})^{\sigma_{0}} \cong L_{\mytinyD_n}(m\omega_{n-1})^{*} \cong L_{\mytinyD_n}(m\omega_{n-1})^{\bowtie}$ when $n \geq 4$ is odd. 
The following theorem statement and proof borrow language and results from \S 3 of \cite{DonSupp}. 

\noindent 
{\bf \OrthogonalSolitaryTheorem}\ \ {\sl Let} $\myX \in \{\myB,\myD\}$, {\sl and let} $\mathfrak{g}(\myA_{1}) \subset \mathfrak{g}(\myA_{2}) \subset \cdots \subset \mathfrak{g}(\myA_{n-1}) \subset \mathfrak{g}(\myX_{n})$ {\sl be the chain of Levi subalgebras induced by the obvious inclusions of Dynkin diagrams.  
For} $\myX_{n} = \myB_{n}$, {\sl set} $L := L_{\mytinyB_n}(m\omega_n)$ {\sl with $V = V[L]$ as the corresponding} $\mathfrak{g}(\myB_{n})${\sl -module $V(m\omega_{n})$; and for} $\myX_{n} = \myD_{n}$, {\sl let $L$ be one of} $L_{\mytinyD_n}(m\omega_{n-1})$ {\sl or} $L_{\mytinyD_n}(m\omega_{n})$ {\sl with $V=V[L]$ as the appropriate corresponding} $\mathfrak{g}(\myD_{n})${\sl -module $V(m\omega_{n-1})$ or $V(m\omega_{n})$.} 
{\sl Then $L$ is solitary and edge-minimal as a supporting graph for $V$, and a weight basis for $V$ restricts irreducibly for the given chain of Levi subalgebras if and only if that weight basis has $L$ as its supporting graph.} 

{\em Proof.} We begin with $L :=  L_{\mytinyB_n}(m\omega_n)$, and set $J := \{1,2,\ldots,n-1\}$. 
For each $\xelt \in L$, observe that $\comp_{J}(\xelt) \cong L_{\mytinyA_{n-1}}^{\mbox{\tiny GT}}(\mysmallP)$ for the partition $\mysmallP = (\myscriptsizeP_{1},\ldots,\myscriptsizeP_{n})$, where $\myscriptsizeP_{i} := c_{n,n+1-i}(\xelt) - c_{n,1}(\xelt)$. 
By Theorem 4.4 of \cite{DonSupp}, $L_{\mytinyA_{n-1}}^{\mbox{\tiny GT}}(\mysmallP)$ restricts irreducibly for the chain of subalgebras $\mathfrak{g}(\myA_{1}) \subset \mathfrak{g}(\myA_{2}) \subset \cdots \subset \mathfrak{g}(\myA_{n-1})$ and, by Lemma 4.3 of \cite{DonSupp}, meets the criteria of Lemma 3.7.A of \cite{DonSupp}. 
Since each $J$-component of $L$ is a classical GT lattice, then $L$ restricts irreducibly for the chain of subalgebras $\mathfrak{g}(\myA_{n-1}) \subset \mathfrak{g}(\myB_{n})$. 
All claims about $L$ from our theorem statement above will now follow from Lemma 3.7.A of \cite{DonSupp}, if we can prove that when $\xelt$ and $\yelt$ are distinct and maximal within their respective $J$-components, then $\wt(\xelt) \ne \wt(\yelt)$. 
Contrapositively, we suppose $\wt(\xelt) = \wt(\yelt)$ for such maximal elements $\xelt$ and $\yelt$, and we demonstrate that $\xelt = \yelt$. 
Now, $\comp_{J}(\yelt) \cong L_{\mytinyA_{n-1}}^{\mbox{\tiny GT}}(\mysmallP')$, where $\mysmallP' = (\myscriptsizeP'_{1},\ldots,\myscriptsizeP'_{n})$ with each $\myscriptsizeP'_{i} := c_{n,n+1-i}(\yelt) - c_{n,1}(\yelt)$. 
Maximality of $\xelt$ means that $c_{n-1,i}(\xelt) = c_{n,i+1}(\xelt)$ for each $i \in \{1,2,\ldots,n-1\}$; the remaining array entries for $\xelt$ are similarly determined.  
The same can be said for $\yelt$. 
So, it suffices to show that $\xelt$ and $\yelt$ agree on elements in positions $(n,n+1-i)$ for $i \in \{1,2,\ldots,n\}$.  
With $\wt(\xelt) = \wt(\yelt)$, it follows that $\mysmallP' = \mysmallP$ and $\mym_{n}(\xelt) = \mym_{n}(\yelt)$. 
Now, for any $\telt \in L$, it is easy to see that $\mym_{n}(\telt) = \sum_{i=1}^{n}[2c_{n,i}(\telt) - c_{n-1,i}(\telt) - c_{n-1,i-1}(\telt)]$. 
Therefore, $2c_{n,1}(\xelt) = \mym_{n}(\xelt) = \mym_{n}(\yelt) = 2c_{n,1}(\yelt)$, and hence $c_{n,1}(\xelt) = c_{n,1}(\yelt)$. 
Thus, for each $i$ we have $c_{n,n+1-i}(\xelt) = \myscriptsizeP_{i}+c_{n,1}(\xelt) = \myscriptsizeP'_{i}+c_{n,1}(\yelt) = c_{n,n+1-i}(\yelt)$. 
This shows that $\xelt$ and $\yelt$ agree on all entries in positions $(n,n+1-i)$, and completes the argument.  
Entirely similar arguments demonstrate the theorem claims about $L_{\mytinyD_n}(m\omega_{n-1})$ and $L_{\mytinyD_n}(m\omega_{n})$.\hfill\QED

%
\vspace*{-0.1in}
\renewcommand{\refname}{\Large \bf References}
\renewcommand{\baselinestretch}{1.1}

\end{document}